%% file: agt-1-5.tex
\documentclass{gtart}   
\input agtout

\lognumber{5}
\volumenumber{1}
\volumeyear{2001}
\papernumber{5}
\published{24 February 2001}
\pagenumbers{73}{114}
\received{6 February 2001}
\accepted{12 February 2001}

\def\semidirect{\mathop
 {\rlap{\kern 7pt\vrule height 5pt width .4pt depth 0pt} \times}}

\newdimen\lengtharrow \newbox\uppertext
\newbox\lowertext \newbox\arrowbox
\def\dimmax #1#2{\ifdim #1<#2 #2\else #1\fi}
\def\arrow #1#2
{\setbox\uppertext=\hbox{$\scriptstyle #1$}
\setbox\lowertext=\hbox{$\scriptstyle #2$}
\lengtharrow=\dimmax{\wd\lowertext}{\wd\uppertext}
\ifdim \lengtharrow<1pt \setbox\arrowbox=\hbox{$\rightarrow$}
\else \advance\lengtharrow by 8pt
\setbox\arrowbox=\hbox to\lengtharrow{\rightarrowfill}
\ht\arrowbox=3.66875 pt \dp\arrowbox=-1.5pt
\setbox\arrowbox=\hbox{$\mathop{\box\arrowbox}
\limits_{\box\lowertext}^{\box\uppertext}$}
\fi
\mathrel{\box\arrowbox}}

\def\PP{{\cal P}}
\def\HH{{\cal H}}
\def\MM{{\cal M}}
\def\e{\varepsilon}
\def\PM{{\cal PM}}
\def\f{\varphi}
\def\BZ{{\bf Z}}
\def\prod{{\rm prod}}

\begin{document} 
\title{Presentations for the punctured mapping\\class groups
in terms of Artin groups}
\shorttitle{Presentations for punctured mapping class groups}
\authors{Catherine Labru\`ere\\Luis Paris}
\asciiauthors{Catherine Labruere and Luis Paris}
\coverauthors{Catherine Labru\noexpand\`ere and Luis Paris}

\address{
Laboratoire de Topologie, UMR 5584 du CNRS\\ 
Universit\'e de Bourgogne, BP 47870 21078 Dijon Cedex, France}
\asciiaddress{
Laboratoire de Topologie, UMR 5584 du CNRS, 
Universite de Bourgogne, BP 47870 21078 Dijon Cedex, France}

\email{clabruer@u-bourgogne.fr, lparis@u-bourgogne.fr}

\begin{abstract}
Consider an oriented compact surface $F$ of positive genus, 
possibly with boundary, and a finite set $\PP$ of punctures in 
the interior of $F$, and define the 
punctured mapping class group of $F$ relatively to $\PP$ 
to be the group of isotopy classes of 
orientation-preserving
homeomorphisms $h: 
F\to F$ which pointwise fix the boundary of $F$ and such 
that $h(\PP) = \PP$. In this paper, we calculate 
presentations for all punctured mapping class groups.
More precisely, we show that these groups are isomorphic with
quotients of Artin groups by some relations involving 
fundamental elements of parabolic subgroups.
\end{abstract}

\asciiabstract{Consider an oriented compact surface F of positive
genus, possibly with boundary, and a finite set P of punctures in the
interior of F, and define the punctured mapping class group of F
relatively to P to be the group of isotopy classes of
orientation-preserving homeomorphisms h: F-->F which pointwise fix the
boundary of F and such that h(P) = P. In this paper, we calculate
presentations for all punctured mapping class groups.  More precisely,
we show that these groups are isomorphic with quotients of Artin
groups by some relations involving fundamental elements of parabolic
subgroups.}

\makeshorttitle

\section{Introduction}

Throughout the paper $F=F_{g,r}$ will denote a compact
oriented surface of genus $g$ with $r$ boundary
components, and $\PP=\PP_n=\{P_1,\dots,P_n\}$ a finite set
of points in the interior of $F$, called {\it punctures}.
We denote by $\HH(F,\PP)$ the group of 
orientation-preserving homeomorphisms $h:F\to F$ that
pointwise fix
the boundary of $F$ and such that $h(\PP)=\PP$. The {\it
punctured mapping class group} $\MM(F,\PP)$ of $F$
relatively to $\PP$ is defined to be the group of isotopy
classes of elements of $\HH(F,\PP)$. Note that the group
$\MM(F,\PP)$ only depends up to isomorphism on the genus
$g$, on the number $r$ of boundary components, and on the
cardinality $n$ of $\PP$. If $\PP$ is empty, then we write
$\MM(F)=\MM(F,\emptyset)$, and call $\MM(F)$ the {\it
mapping class group} of $F$.

The {\it pure mapping class group} of $F$ relatively to
$\PP$ is defined to be the subgroup $\PM(F,\PP)$ of
isotopy classes of elements of $\HH(F,\PP)$ that pointwise
fix $\PP$. Let $\Sigma_n$ denote the symmetric group of
$\{1,\dots,n\}$. Then the punctured mapping class group
and the pure mapping class group are related by the
following exact sequence.
$$
1\to \PM(F,\PP_n)\to\MM(F,\PP_n)\to\Sigma_n\to 1\ .
$$

A {\it Coxeter matrix} is a matrix
$M=(m_{i,j})_{i,j=1,\dots,l}$ satisfying:

\smallskip
$\bullet$\qua $m_{i,i}=1$ for all $i=1,\dots, l$;

\smallskip
$\bullet$\qua $m_{i,j}=m_{j,i}\in\{2,3,4,\dots,\infty\}$, for
$i\neq j$.

\smallskip\noindent
A Coxeter matrix $M=(m_{i,j})$ is usually represented by
its {\it Coxeter graph} $\Gamma$. This is defined by the
following data:

\smallskip
$\bullet$\qua $\Gamma$ has $l$ vertices: $x_1,\dots, x_{l}$;

\smallskip
$\bullet$\qua two vertices $x_i$ and $x_j$ are joined by an
edge if $m_{i,j}\ge 3$;

\smallskip
$\bullet$\qua the edge joining two vertices $x_i$ and $x_j$ is
labelled by $m_{i,j}$ if $m_{i,j}\ge 4$.

\smallskip
For $i,j\in\{1,\dots,l\}$, we write:
$$
\prod(x_i,x_j,m_{i,j})=\left\{\matrix{
(x_ix_j)^{m_{i,j}/2}\quad\hfill&{\rm if}\ m_{i,j}\ {\rm
is\ even},\hfill\cr
(x_ix_j)^{(m_{i,j}-1)/2}x_i\quad\hfill&{\rm if}\ m_{i,j}\
{\rm is\ odd}.\hfill\cr}
\right.
$$
The {\it Artin group} $A(\Gamma)$ associated with $\Gamma$
(or with $M$) is the group given by the presentation:
$$
A(\Gamma)=\langle x_1,\dots,x_{l}\,|\,\prod(x_i,x_j,m_{i,j})\!
=\!\prod(x_j,x_i,m_{i,j})\ {\rm if}\ i\neq j\ {\rm and}\
m_{i,j}<\infty\rangle.
$$
The {\it Coxeter group} $W(\Gamma)$ associated with
$\Gamma$ is the quotient of $A(\Gamma)$ by the relations
$x_i^2=1$, $i=1,\dots, l$. We say that $\Gamma$ or
$A(\Gamma)$ is of {\it finite type} if $W(\Gamma)$ is
finite.

For a subset $X$ of the set $\{x_1,\dots,x_{l}\}$ of
vertices of $\Gamma$, we denote by $\Gamma_X$ the Coxeter
subgraph of $\Gamma$ generated by $X$, by $W_X$ the
subgroup of $W(\Gamma)$ generated by $X$, and by $A_X$ the
subgroup of $A(\Gamma)$ generated by $X$. It is a non-trivial 
but well known fact that $W_X$ is the Coxeter
group associated with $\Gamma_X$ (see \cite{Bo}), and
$A_X$ is the Artin group associated with $\Gamma_X$ (see
\cite{Lek}, \cite{Pa}). Both $W_X$ and $A_X$ are called 
{\it parabolic subgroups} of $W(\Gamma)$ and of $A(\Gamma)$,
respectively.

Define the {\it quasi-center} of an Artin group $A(\Gamma)$ 
to be the subgroup of elements $\alpha$ in $A(\Gamma)$ 
satisfying $\alpha X \alpha^{-1} = X$, where $X$ is the 
natural generating set of $A(\Gamma)$. If $\Gamma$ is of 
finite type and connected, then the quasi-center is an 
infinite cyclic group generated by a special element of 
$A(\Gamma)$, called {\it fundamental element}, and denoted 
by $\Delta(\Gamma)$ (see \cite{Del}, \cite{BrSa}).

\medskip
The most significant work on presentations for mapping class 
groups is certainly the paper \cite{HT} of Hatcher and 
Thurston. In this paper, the authors introduced a simply 
connected complex on which the mapping class group $\MM 
(F_{g,0})$ acts, and, using this action and following a 
method due to Brown \cite{Bro}, they obtained a 
presentation for $\MM(F_{g,0})$. However, as pointed out by 
Wajnryb \cite{Waj}, this presentation is rather complicated 
and requires many generators and relations. Wajnryb 
\cite{Waj} used this presentation of Hatcher and Thurston 
to calculate new presentations for $\MM (F_{g,1})$ and for 
$\MM (F_{g,0})$. He actually presented $\MM (F_{g,1})$ as 
the quotient of an Artin group by two relations, and 
presented $\MM (F_{g,0})$ as the quotient of the same Artin 
group by the same two relations plus another one. In 
\cite{Mat}, Matsumoto showed that these three relations are 
nothing else than equalities among powers of fundamental 
elements of parabolic subgroups. Moreover, he showed how to 
interpret these powers of fundamental elements inside the 
mapping class group. Once this interpretation is known, the 
relations in Matsumoto's presentations become trivial. At 
this point, one has ``good'' presentations for 
$\MM(F_{g,1})$ and for $\MM(F_{g,0})$, in the sence that 
one can remember them. Of course, the definition of a 
``good'' presentation depends on the memory of the reader 
and on the time he spends working on the presentation.

One can find in \cite{Loo} another presentation for $\MM 
(F_{g,1})$ as the quotient of an Artin group by relations 
involving fundamental elements of parabolic subgroups. 
Recently, Gervais \cite{Ger} found another ``good'' 
presentation for $\MM (F_{g,r})$ with many generators but 
simple relations.

In the present paper, starting from Matsumoto's 
presentations, we calculate presentations for all punctured 
mapping class groups $\MM(F_{g,r},\PP_n)$ as quotients of 
Artin groups by some relations which involve fundamental 
elements of parabolic subgroups. In particular, $\MM( 
F_{g,0}, \PP_n)$ is presented as the quotient of an Artin 
group by five relations, all of them being equalities among 
powers of fundamental elements of parabolic subgroups. 

The generators in our presentations are Dehn twists and 
braid twists. We define them in Subsection 2.1, and we show 
that they verify some ``braid'' relations that allow us to 
define homomorphisms from Artin groups to punctured mapping 
class groups. The main algebraic tool we use is Lemma 2.5, 
stated in Subsection 2.2, which says how to find a 
presentation for a group $G$ from an exact sequence $1\to 
K\to G\to H\to 1$ and from presentations of $K$ and $H$. We 
also state in Subsection 2.2 some exact sequences involving 
punctured mapping class groups on which Lemma 2.5 will be 
applied. In order to find our presentations, we first need 
to investigate some homomorphisms from finite type Artin 
groups to punctured mapping class groups, and to calculate 
the images under these homomorphisms of some powers of 
fundamental elements. This is the object of Subsection 2.3. 
Once these images are known, one can easily verify that the 
relations in our presentations hold. Of course, it remains 
to prove that no other relation is needed. We state our 
presentation for $\MM (F_{g,r+1}, \PP_n)$ (where $g\ge 1$, 
and $r,n\ge 0$) in Theorem 3.1, and we state our presentation for 
$\MM (F_{g,0}, \PP_n)$ (where $g,n\ge 1$) in Theorem 3.2. 
Then, Subsection 3.1 is dedicated to the proof of Theorem 
3.1, and Subsection 3.2 is dedicated to the proof of Theorem 3.2.


\section{Preliminaries}

\subsection{Dehn twists and braid twists}

We introduce in this subsection some elements of the
punctured mapping class group, the Dehn twists and the
braid twists, which will play a prominent r\^ole throughout
the paper. In particular, the generators for the punctured
mapping class group will be chosen among them.

By an {\it essential circle} in $F\setminus \PP$ we mean an
embedding $s:S^1\to F\setminus\PP$ of the circle whose
image 
is in the interior of $F\setminus\PP$ and
does not bound a disk in $F\setminus\PP$. Two
essential circles $s,s'$ are called {\it isotopic} if there
exists $h\in\HH(F,\PP)$ which represents the identity in
$\MM(F,\PP)$ and such that $h\circ s=s'$. Isotopy of
circles is an equivalence relation which we denote by
$s\simeq s'$. Let $s:S^1\to F\setminus\PP$ be an essential
circle. We choose an embedding $A:[0,1]\times S^1\to
F\setminus\PP$ of the annulus such that $A({1\over 2},z)=s(z)$ for
all $z\in S^1$, and we consider the homeomorphism
$T\in\HH(F,\PP)$ defined by
$$
(T\circ A)(t,z)=A(t,e^{2i\pi t}z),\quad t\in[0,1],\ z\in 
S^1,
$$
and $T$ is the identity on the exterior of the image of $A$
(see Figure \ref{f1}). The {\it Dehn twist} along $s$ is defined
to be the element $\sigma\in\MM(F,\PP)$ represented by $T$. Note
that:

\smallskip
$\bullet$\qua the definition of $\sigma$ does not depend on the
choice of $A$;

\smallskip
$\bullet$\qua the element $\sigma$ does not depend on the orientation of
$s$;

\smallskip
$\bullet$\qua if $s$ and $s'$ are isotopic, then their corresponding
Dehn twists are equal;

\smallskip
$\bullet$\qua if $s$ bounds a disk in $F$ which contains exactly one
puncture,then $\sigma=1$; otherwise, $\sigma$ is of
infinite order;

\smallskip
$\bullet$\qua if 
$\xi\in\MM(F,\PP)$ is represented by $f\in\HH(F,\PP)$, then
$\xi\sigma\xi^{-1}$ is the Dehn twist along $f(s)$.

\begin{figure}[ht]
\centerline{\input{pmc_1}}
\caption{\label{f1} Dehn twist along $s$}
\end{figure}

\medskip
By an {\it arc} we mean an embedding $a:[0,1]\to F$ 
of the segment 
whose image is in the interior of $F$,
such
that $a((0,1))\cap\PP=\emptyset$, and such that both $a(0)$ and $a(1)$
are punctures. Two arcs $a,a'$ are called {\it isotopic} if
there exists $h\in\HH(F,\PP)$ which represents the identity
in $\MM(F,\PP)$ and such that $h\circ a=a'$. Note that
$a(0)=a'(0)$ and $a(1)=a'(1)$ if $a$ and $a'$ are isotopic.
Isotopy of arcs is an equivalence relation which we denote
by $a\simeq a'$. Let $a$ be an arc. We choose an embedding
$A:D^2\to F$ of the unit disk satisfying:

\smallskip
$\bullet$\qua $a(t)=A(t-{1\over 2})$ for all $t\in [0,1]$,

\smallskip
$\bullet$\qua $A(D^2)\cap\PP=\{a(0),a(1)\}$,

\smallskip\noindent
and we consider the homeomorphism $T\in\HH(F,\PP)$ defined
by
$$
(T\circ A)(z)=A(e^{2i\pi|z|}z),\quad z\in D^2,
$$
and $T$ is the identity on the exterior of the image of $A$
(see Figure \ref{f2}). The {\it braid twist} along $a$ is defined
to be the element $\tau\in\MM(F,\PP)$ represented by $T$. Note
that:

\smallskip
$\bullet$\qua the definition of $\tau$ does not depend on the
choice of $A$;

\smallskip
$\bullet$\qua if $a$ and $a'$ are isotopic, then their corresponding
braid twists are equal;

\smallskip
$\bullet$\qua if 
$\xi\in\MM(F,\PP)$ is represented by $f\in\HH(F,\PP)$, then
$\xi\tau\xi^{-1}$ is the braid twist along $f(a)$;

\smallskip
$\bullet$\qua if $s:S^1\to F\setminus\PP$ is the essential circle defined
by $s(z)=A(z)$ (see Figure \ref{f2}), then $\tau^2$ is the Dehn
twist along $s$.

\begin{figure}[ht]
\centerline{\input{pmc_2}}
\caption{\label{f2} Braid twist along a}
\end{figure}

\medskip
We turn now to describe some relations among Dehn twists and
braid twists which will be essential to define
homomorphisms from Artin groups to punctured
mapping class groups.

The first family of relations are known as ``braid
relations'' for Dehn twists (see \cite{Bi3}).

\medskip\noindent
{\bf Lemma 2.1}\qua {\sl Let $s$ and $s'$ be two essential
circles which intersect transversely, and let $\sigma$ and
$\sigma'$ be the Dehn twists along $s$ and $s'$,
respectively. Then:
$$
\matrix{
&\sigma\sigma'=\sigma'\sigma\quad\hfill&{\hbox{\sl if}}\ s\cap s'
=\emptyset,\hfill\cr
&\sigma\sigma'\sigma=\sigma'\sigma\sigma'\quad\hfill&{\hbox{\sl
if}}\ |s\cap s'|=1.}\eqno{\lower .5em\hbox{\qed}}
$$}

The next family of relations are simply the usual braid
relations viewed inside the punctured mapping class group.

\medskip\noindent
{\bf Lemma 2.2}\qua {\sl Let $a$ and $a'$ be two arcs, and let
$\tau$ and $\tau'$ be be the braid twists along $a$ and
$a'$, respectively. Then:
$$
\matrix{
&\tau\tau'=\tau'\tau\quad\hfill&{\hbox{\sl if}}\ a\cap a'=
\emptyset,\hfill\cr
&\tau\tau'\tau=\tau'\tau\tau'\quad\hfill&{\hbox{\sl if}}\ a(0)=
a'(1)\ {\hbox{\sl and}}\ a\cap a'=\{a(0)\}.}\eqno{\lower .5em\hbox{\qed}}
$$}

To our knowledge, the last family of relations does not appear in the
literature. However, their proofs are easy and are
left to the reader.

\medskip\noindent
{\bf Lemma 2.3}\qua {\sl Let $s$ be an essential circle, and let $a$
be an arc which intersects $s$ transversely. Let $\sigma$ be
the Dehn twist along $s$, and let $\tau$ be the braid twist along
$a$. Then:
$$
\matrix{
&\sigma\tau=\tau\sigma\quad\hfill&{\hbox{\sl if}}\ s\cap a=
\emptyset,\hfill\cr
&\sigma\tau\sigma\tau=\tau\sigma\tau\sigma\quad\hfill&
{\hbox{\sl if}}\ |s\cap a|=1.
\cr}\eqno{\lower .5em\hbox{\qed}}
$$}

We finish this subsection by recalling another relation called
{\it lantern relation} (see \cite{Jo}) which is not used
to define homomorphisms between Artin
groups and punctured mapping class groups, but which will
be useful in the remainder.

We point out first that we use the convention in figures
that a letter which appears over a circle or an arc denotes
the corresponding Dehn twist or braid twist, and not the
circle or the arc itself.

\medskip\noindent
{\bf Lemma 2.4}\qua {\sl Consider an embedding of $F_{0,4}$ in
$F\setminus\PP$ and the Dehn twists $e_1,e_2,e_3,e_4,
a,b,c$ represented in Figure \ref{f3}. Then
$$
e_1e_2e_3e_4=abc.\eqno{\qed}
$$}

\begin{figure}[ht]
\centerline{\input{pmc_3}}
\caption{\label{f3} Lantern relation}
\end{figure}


\subsection{Exact sequences}

Now, we introduce in Lemma 2.5 our main tool to obtain
presentations for the punctured mapping class groups.
Briefly, this lemma says how to find a presentation for a
group $G$ from an exact sequence $1\to K\to G\to H\ \to 1$
and from presentations of $H$ and $K$. This lemma will be
applied to the exact sequences (2.1), (2.2), and (2.3) given after
Lemma 2.5.

Consider an exact sequence
$$
1\to K\to G\arrow{\rho}{} H \to 1
$$
and presentations $H=\langle S_H|R_H\rangle$, $K=\langle S_K|
R_K\rangle$ for $H$ and $K$, respectively. For all $x\in
S_H$, we fix some $\tilde x\in G$ such that $\rho(\tilde
x)=x$, and we write
$$
\tilde S_H=\{\tilde x\ ;\ x\in S_H\}.
$$
Let $r=x_1^{\e_1}\dots x_l^{\e_l}$ in $R_H$. Write $\tilde
r=\tilde x_1^{\e_1}\dots \tilde x_l^{\e_l}\in G$. Since
$r$ is a relator of $H$, we have $\rho(\tilde r)=1$. Thus,
$S_K$ being a generating set of the kernel of $\rho$, one
may choose a word $w_r$ over $S_K$ such that both $\tilde
r$ and $w_r$ represent the same element of $G$. Set
$$
R_1=\{\tilde rw_r^{-1}\ ;\ r\in R_H\}.
$$
Let $\tilde x\in\tilde S_H$ and $y\in S_K$. Since $K$ is a
normal subgroup of $G$, $\tilde xy\tilde x^{-1}$ is also
an element of $K$, thus one may choose a word $v(x,y)$
over $S_K$ such that both $\tilde xy\tilde x^{-1}$ and
$v(x,y)$ represent the same element of $G$. Set
$$
R_2=\{\tilde xy\tilde x^{-1}v(x,y)^{-1}\ ;\ \tilde x\in
\tilde S_H\ {\rm and}\ y\in S_K\}.
$$
The proof of the following lemma is left to the reader.

\medskip\noindent
{\bf Lemma 2.5}\qua {\sl G admits the presentation
$$
G=\langle \tilde S_H\cup S_K\ |\ R_1\cup R_2\cup R_K
\rangle.\eqno{\qed}
$$}

The first exact sequence on which we will apply Lemma 2.5
is the one given in the introduction:
$$
1\to\PM(F,\PP_n)\to\MM(F,\PP_n)\to\Sigma_n\to 1,\leqno (2.1)
$$
where $\Sigma_n$ denotes the symmetric group of $\{1,\dots,n\}$.

The inclusion $\PP_{n-1}\subset\PP_n$ gives rise to a
homomorphism $\f_n:\PM(F,\PP_n)\to\PM(F,\PP_{n-1})$. By
\cite{Bi2}, if $(g,r,n)\neq(1,0,1)$, then we have the
following exact sequence:
$$
1\to\pi_1(F\setminus \PP_{n-1},P_n)\arrow{\iota_n}{}
\PM(F,\PP_n)\arrow{\f_n}{} \PM(F,\PP_{n-1}) \to 1.
\leqno (2.2)
$$

We will need later a more precise description of
the images by $\iota_n$ of certain elements of
$\pi_1(F\setminus\PP_{n-1},P_n)$. Consider an essential
circle $\alpha:S^1\to F\setminus\PP_{n-1}$ such that
$\alpha(1)=P_n$. Here, we assume that $\alpha$ is oriented. 
Let $\xi$ be the element of
$\pi_1(F\setminus\PP_{n-1},P_n)$ represented by
$\alpha$.
We
choose an embedding $A:[0,1]\times S^1\to
F\setminus\PP_{n-1}$ of the annulus such that
$A({1\over 2},z)=\alpha(z)$ for all $z\in S^1$ (see Figure \ref{f4}). Let
$s_0,s_1:S^1\to F\setminus\PP_n$ be the essential circles
defined by
$$
s_0(z)=A(0,z),\quad s_1(z)=A(1,z),\quad z\in S^1,
$$
and let $\sigma_0,\sigma_1$ be the Dehn twists along $s_0$ and
$s_1$, respectively. Then the following holds.

\begin{figure}[ht]
\centerline{\input{pmc_4}}
\caption{\label{f4} Image of a simple circle by $\iota_n$}
\end{figure}

\medskip\noindent
{\bf Lemma 2.6}\qua {\sl We have $\iota_n(\xi)=\sigma_0^{
-1}\sigma_1$. \qed}

\medskip
Now, consider a surface $F_{g,r+m}$ of genus $g$
with $r+m$ boundary components, and a set
$\PP_n=\{P_1,\dots, P_n\}$ of $n$ punctures in the
interior of $F_{g,r+m}$. Choose $m$ boundary curves
$c_1,\dots, c_m:S^1\to \partial F_{g,r+m}$. Let $F_{g,r}$
be the surface of genus $g$ with $r$ boundary components
obtained from $F_{g,r+m}$ by gluing a disk $D_i^2$ along
$c_i$, for all $i=1,\dots, m$, and let
$\PP_{n+m}=\{P_1,\dots,P_n,Q_1,\dots, Q_m\}$ be a set of
punctures in the interior of $F_{g,r}$, where $Q_i$ is
chosen in the interior of $D_i^2$, for all $i=1,\dots, m$.
The proof of the following exact sequence can be found in
\cite{PaRo}.

\medskip\noindent
{\bf Lemma 2.7}\qua {\sl Assume that  
$(g,r,m)\not\in\{(0,0,1), (0,0,2)\}$.
Then we have the exact sequence:
$$
1\to \BZ^m\to\PM(F_{g,r+m},\PP_n)\to\PM(F_{g,r},\PP_{n+m})
\to 1\ ,\leqno(2.3)
$$
where $\BZ^m$ stands for the free abelian group of rank $m$
generated by the Dehn twists along the $c_i$'s.} \qed


\subsection{Geometric representations of Artin groups}

Define a {\it geometric representation} of an Artin group
$A(\Gamma)$ to be a homomorphism from $A(\Gamma)$ to some
punctured mapping class group. In this subparagraph, we
describe some geometric representations of Artin groups
whose properties will be used later in the paper.

The first family of geometric representations has been
introduced by Perron and Vannier for studying geometric
monodromies of simple singularities \cite{PeVa}. A {\it chord
diagram} in the disk $D^2$ is a family $S_1,\dots,
S_{l}:[0,1]\to D^2$ of segments satisfying:

\smallskip
$\bullet$\qua $S_i:[0,1]\to D^2$ is an embedding for all
$i=1,\dots,l$;

\smallskip
$\bullet$\qua $S_i(0),S_i(1)\in\partial D^2$, and
$S_i((0,1))\cap\partial D^2=\emptyset$, for all
$i=1,\dots, l$;

\smallskip
$\bullet$\qua either $S_i$ and $S_j$ are disjoint, or they
intersect transversely in a unique point in the interior
of $D^2$, for $i\neq j$.

\smallskip\noindent
From this data, one can first define a Coxeter matrix
$M=(m_{i,j})_{i,j=1,\dots, l}$ by seting $m_{i,j}=2$ if
$S_i$ and $S_j$ are disjoint, and $m_{i,j}=3$ if $S_i$ and
$S_j$ intersect transversely in a point. The Coxeter graph
$\Gamma$ associated with $M$ is called {\it intersection
diagram} of the chord diagram. It is an ``ordinary'' graph
in the sence that none of the edges has a
label. From the chord diagram we can also define a surface
$F$ by attaching to $D^2$ a handle $H_i$ which joins both
extremities of $S_i$, for all $i=1,\dots,l$ (see Figure \ref{f5}).
Let $\sigma_i$ be the Dehn twist along the circle made
up with the segment $S_i$ together with the central curve
of $H_i$. By Lemma 2.1, one has a geometric representation
$A(\Gamma)\to\MM(F)$ which sends $x_i$ on $\sigma_i$ for
all $i=1,\dots, l$. This geometric representation will be called
{\it Perron-Vannier representation}.

\begin{figure}[ht]
\centerline{\input{pmc_5}}
\caption{\label{f5} Chord diagram and associated surface
and Dehn twists}
\end{figure}

\medskip
If $\Gamma$ is connected, then the Perron-Vannier
representation is injective if and only if $\Gamma$ is of
type $A_l$ or $D_l$ \cite{La2}, \cite{Waj2}. In the case
where $\Gamma$ is of type $A_l$, $D_l$, $E_6$, or $E_7$,
the vertices of $\Gamma$ will be numbered according to
Figure \ref{f6}, and the Dehn twists $\sigma_1,\dots,\sigma_l$
are those represented in Figures \ref{f7a}, \ref{f7b}, \ref{f7c}.

\begin{figure}[ht]
\centerline{\input{pmc_6}}
\caption{\label{f6} Some finite type Coxeter graphs}
\end{figure}

\begin{figure}[ht]
\centerline{\input{pmc_7a}}
\caption{\label{f7a} Perron-Vannier representations of
type $A_l$}
\end{figure}

\begin{figure}[ht]
\centerline{\input{pmc_7b}}
\caption{\label{f7b} Perron-Vannier representations of
type $D_l$}
\end{figure}

\begin{figure}[ht]
\centerline{\input{pmc_7c}}
\caption{\label{f7c} Perron-Vannier representations of
type $E_6$ and $E_7$}
\end{figure}

\medskip
The Perron-Vannier representation of the Artin group of
type $A_{l-1}$ can be extended to a geometric representation of the
Artin group of type $B_l$ as follows.
First, we number the vertices of $B_l$ according to Figure \ref{f6}.
Then $A_{l-1}$ is the subgraph of $B_{l}$ generated by
the vertices $x_2,\dots,x_{l}$. We start from a chord
diagram $S_2,\dots, S_{l}$ whose intersection diagram is
$A_{l-1}$, and we denote by $F$ the associated surface.
For $i=2,\dots, l$, we denote by $s_i$ the essential
circle of $F$ made up with $S_i$ and the central curve of
the handle $H_i$. We can choose two points $P_1,P_2$ in
the interior of $F$ and an arc $a_1$ from $P_1$ to $P_2$
satisfying:

\smallskip
$\bullet$\qua $\{P_1,P_2\}\cap s_i=\emptyset$ for all
$i=2,\dots, l$;

\smallskip
$\bullet$\qua $a_1\cap s_i=\emptyset$ for all $i=3,\dots, l$,
and $a_1$ and $s_2$ intersect transversely in a unique point (see
Figure \ref{f8}).

\smallskip\noindent
Let $\tau_1$ be the braid twist along $a_1$, and let
$\sigma_i$ be the Dehn twist along $s_i$, for
$i=2,\dots,l$. By Lemma 2.3, there is a well defined
homomorphism $A(B_l)\to\MM(F,\{P_1,P_2\})$ which sends
$x_1$ on $\tau_1$, and $x_i$ on $\sigma_i$ for
$i=2,\dots,l$. It is shown in \cite{La1} that this
geometric representation is injective.

\begin{figure}[ht]
\centerline{\input{pmc_8}}
\caption{\label{f8} Perron-Vannier representation of type
$B_l$}
\end{figure}

\medskip
Now, consider a graph $G$ embedded in a surface $F$. Here,
we assume that $G$ has no loop and no multiple-edge. Let
$\PP=\{P_1,\dots,P_n\}$ be the set of vertices of $G$, and
let $a_1,\dots, a_l$ be the edges. Define the Coxeter
matrix $M=(m_{i,j})_{i,j=1,\dots, l}$ by $m_{i,j}=3$ if
$a_i$ and $a_j$ have a common vertex, and $m_{i,j}=2$
otherwise. Denote by $\Gamma$ the Coxeter graph associated
with $M$. By Lemma 2.2, one has a homomorphism
$A(\Gamma)\to \MM(F,\PP)$ which associates with $x_i$ the
braid twist $\tau_i$ along $a_i$, for all $i=1,\dots, l$.
This homomorphism will be called {\it graph
representation} of $A(\Gamma)$. Its image clearly belongs
to the surface braid group of $F$ based at $\PP$. 
The particular case where $F$ is a disk has been
studied by Sergiescu \cite{Ser} to find new presentations
for the Artin braid groups. Graph representations
have been also used by Humphries \cite{Hum2}
to solve some Tits' conjecture.

Assume now that $G$ is a line in a cylinder $F=S^1\times
I$. Let $a_2,\dots, a_l$ be the edges of $G$, and let
$\PP_l=\{P_1,\dots, P_l\}$ be the set of vertices. Choose
an essential circle $s_1:S^1\to F\setminus\PP$ such that:

\smallskip
$\bullet$\qua $s_1$ does not bound a disk in $F$;

\smallskip
$\bullet$\qua $s_1\cap a_i=\emptyset$ for all $i=3,\dots, l$,
and $s_1$ and $a_2$ intersect transversely in a unique
point (see Figure \ref{f9}).

\smallskip\noindent
Let $\sigma_1$ be the Dehn twist along $s_1$, and let
$\tau_i$ be the braid twist along $a_i$ for
$i=2,\dots,l$. By Lemma 2.3, there is a well defined
homomorphism $A(B_l)\to\MM(S^1\times I,\PP_l)$ which sends
$x_1$ on $\sigma_1$, and $x_i$ on $\tau_i$ for $i=2,\dots,
l$. This homomorphism is clearly an extension of the
graph representation of $A(A_{l-1})$ in $\MM(S^1\times
I,\PP_l)$.

\begin{figure}[ht]
\centerline{\input{pmc_9}}
\caption{\label{f9} Graph representation of type $B_{l}$}
\end{figure}

Let $\Gamma$ be a finite type connected graph. Recall
that the {\it quasi-center} of $A(\Gamma)$ is the subgroup of
elements $\alpha$ in $A(\Gamma)$ satisfying $\alpha X \alpha^{-1}
=X$, where $X$ is the natural generating set of $A(\Gamma)$, and that
this subgroup is an infinite cyclic group generated by some special
element of $A(\Gamma)$, called {\it fundamental element}, and denoted by
$\Delta (\Gamma)$. (see \cite{BrSa} and \cite{Del}). The center of
$A(\Gamma)$ is an infinite cyclic group generated by $\Delta (\Gamma)$
if $\Gamma$ is $B_l$, $D_l$ ($l$ even), $E_7$, $E_8$, $F_4$,
$H_3$, $H_4$, and $I_2(p)$ ($p$ even), and by $\Delta^2(\Gamma)$
if $\Gamma$ is $A_l$, $D_l$ ($l$ odd), $E_6$, and $I_2(p)$ ($p$
odd). Explicit expressions of $\Delta(\Gamma)$ and of $\Delta^2 
(\Gamma)$ can be found in \cite{BrSa}. In the
remainder, we will need the following ones.

\medskip\noindent
{\bf Proposition 2.8}\qua (Brieskorn, Saito \cite{BrSa})\qua {\sl
We number the vertices of $A_{l}$, $B_{l}$, $D_{l}$, $E_6$,
and $E_7$ according to Figure \ref{f6}.
\begin{eqnarray*}
\Delta^2(A_{l})&=&(x_1x_2\dots x_l)^{l+1}\ ,\\
\Delta(B_{l})&=&(x_1x_2\dots x_l)^{l}\ ,\\
\Delta(D_{2p})&=&(x_1x_2\dots x_{2p})^{2p-1}\ ,\\
\Delta^2(D_{2p+1})&=&(x_1x_2\dots x_{2p+1})^{4p}
\ ,\\
\Delta^2(E_6)&=&(x_1x_2\dots x_6)^{12}\ ,\\
\Delta(E_7)&=&(x_1x_2\dots x_7)^{15}\ .\\
\end{eqnarray*}}\vglue -4.5em\noindent\hbox{}\hfill\qed

\medskip

We will also need the following well known
equalities (see \cite{Pa2}).

\medskip\noindent
{\bf Proposition 2.9}\qua {\sl We number the vertices of $A_l$, $B_l$,
and $D_l$ according to Figure 6. Then:}

\centerline{\vbox{\halign{
\hfill$#$&\ $#$\ &$#$\hfill\cr
\Delta (A_l) &=& x_1 \dots x_l \cdot \Delta (A_{l-1}), \cr
\Delta (B_l) &=& x_l \dots x_2 x_1 x_2 \dots x_l \cdot \Delta (B_{l-1}), \cr
\Delta (D_l) &=& x_l \dots x_3 x_1 x_2 x_3 \dots x_l \cdot \Delta
(D_{l-1}). \cr}}}\vglue -2.5em\noindent\hbox{}\hfill\qed

\medskip

Our goal now is to determine the images
under Perron-Vannier representations and under
graph representations of some powers of fundamental
elements (Proposition 2.12). To do so, we first
need to know generating sets for the punctured
mapping class groups. So, we prove the following.

\medskip\noindent
{\bf Proposition 2.10}\qua {\sl Let $g\ge 1$ and $r,n\ge 0$.

\smallskip
{\rm(i)}\qua $\PM(F_{g,r+1},\PP_n)$ is generated by the Dehn twists
$a_0,\dots, a_{n+r}, b_1,\dots, b_{2g-1}$, $c$, $d_1, \dots,
d_r$ represented in Figure \ref{f10}.

\smallskip
{\rm(ii)}\qua
$\MM(F_{g,r+1},\PP_n)$ is generated by the Dehn twists
$a_0, \dots, a_r, a_{r+1}$, $b_1,..,b_{2g-1}$, \linebreak
$c$, $d_1,\dots, d_r$, and the braid twists $\tau_1,\dots, \tau_{n
-1}$ represented in Figure \ref{f10}.}

\begin{figure}[ht]
\centerline{\input{pmc_10}}
\caption{\label{f10} Generators for $\PM(F_{g,r+1},\PP_n)$
and $\MM(F_{g,r+1},\PP_n)$}
\end{figure}

\medskip\noindent
{\bf Corollary 2.11}\qua {\sl Let $g\ge 1$ and $n\ge 0$.

\smallskip
{\rm(i)}\qua $\PM(F_{g,0},\PP_n)$ is generated by the Dehn twists
$a_0,\dots, a_n$, $b_1,\dots, b_{2g-1}$, $c$ represented
in Figure \ref{f11}.

\smallskip
{\rm(ii)}\qua $\MM(F_{g,0},\PP_n)$ is generated by the Dehn twists
$a_0,a_1$, $b_1,\dots, b_{2g-1}$, $c$, and the braid
twists $\tau_1,\dots, \tau_{n-1}$ represented in Figure \ref{f11}.}

\begin{figure}[ht]
\centerline{\input{pmc_11}}
\caption{\label{f11} Generators for $\PM(F_{g,0},\PP_n)$
and $\MM(F_{g,0},\PP_n)$}
\end{figure}

\medskip\noindent
{\bf Proof}\qua The key argument of the proof of Proposition 2.10 
is the following
remark stated as Assertion 1, and which we apply to the
exact sequences (2.1), (2.2), and (2.3) of Subsection 2.2.

\medskip\noindent
{\bf Assertion 1}\qua {\sl Let
$$
1\to K\to G\arrow{\rho}{} H\to 1
$$
be an exact sequence, and let $S_H, S_K$ be generating sets
of $H$ and $K$, respectively. For each $x\in S_H$ we
choose $\tilde x\in G$ such that $\rho(\tilde x)=x$, and
we write $\tilde S_H=\{\tilde x;x\in S_H\}$. Then $S_K\cup
\tilde S_H$ generates $G$.}

\medskip
First, we prove by induction on $n$ that
$\PM(F_{g,1},\PP_n)$ is generated by $a_0,\dots, a_n$,
$b_1,\dots, b_{2g-1}$, $c$. The case $n=0$ is proved in
\cite{Hum}. So, we assume that $n>0$. By the inductive
hypothesis, $\PM(F_{g,1},\PP_{n-1})$ is generated by
$a_0,\dots, a_{n-1}$, $b_1,\dots, b_{2g-1}$, $c$. On the other
hand, $\pi_1(F_{g,1}\setminus\PP_{n-1},P_n)$ is the free
group generated by the loops $\alpha_1,\dots, \alpha_n,
\beta_1,\dots, \beta_{2g-1}$ represented in Figure \ref{f12}.
Applying Assertion 1 to the exact sequence (2.2), one has
that $\PM(F_{g,1},\PP_n)$ is generated by $a_0,\dots,
a_{n-1}$, $b_1,\dots, b_{2g-1}$, $c$, $\alpha_1,\dots, \alpha_n$,
$\beta_1,\dots, \beta_{2g-1}$. 
One can directly verify the
following equalities:
$$\matrix{
\alpha_i=\hfill&(b_1a_na_{i-1}b_1a_{n-1})^{-1}\alpha_n^{-1}
(b_1a_na_{i-1}b_1a_{n-1}),\hfill&\quad i=1,\dots,n-1,\hfill\cr
\beta_1=\hfill&(b_1a_{n-1})^{-1}\alpha_n(b_1a_{n-1}),\hfill
\cr
\beta_j=\hfill&(b_jb_{j-1})^{-1}\beta_{j-1}
(b_j b_{j-1}),\hfill&\quad j=2,\dots, 2g-1.\hfill
\cr}
$$
and, from Proposition 2.6, one has:
$$
\alpha_n = a_{n-1}^{-1} a_n,
$$
thus $\PM(F_{g,1},\PP_n)$ is generated by $a_0,\dots, a_n$,
$b_1,\dots, b_{2g-1}$, $c$.

\begin{figure}[ht]
\centerline{\input{pmc_12}}
\caption{\label{f12} Generators for $\pi_1(F_{g,1}
\setminus \PP_{n-1}, P_n)$}
\end{figure}

\medskip
Now, applying Assertion 1 to (2.3), one
has that $\PM(F_{g,r+1},\PP_n)$ is generated by
$a_0,\dots, a_{n+r}$, $b_1, \dots, b_{2g-1}$, $c$, $d_1, \dots,
d_r$.

\medskip\noindent
{\bf Assertion 2}\qua {\sl Let $a_0,a_1,a_2$ be the Dehn
twists and $\tau$ the braid twist in $\MM(S^1\times I,
\{P_1,P_2\})$ represented in Figure \ref{f13}. Then
$$
\tau a_1 \tau a_1=a_0a_2\ .
$$}

\begin{figure}[ht]
\centerline{\input{pmc_13}}
\caption{\label{f13} A relation in $\MM(S^1\times I,
\{P_1,P_2\})$}
\end{figure}

\noindent
{\bf Proof of Assertion 2}\qua We consider the Dehn twist
$a_3$ along a circle which bounds a small disk in
$S^1\times I$ which contains $P_1$, and the Dehn twist
$a_4$ along a circle which bounds a small disk in
$S^1\times I$ which contains $P_2$. As pointed out in
Subsection 2.1, we have $a_3=a_4=1$. The lantern relation
of Lemma 2.4 says:
$$
\tau^2\cdot a_1\cdot \tau a_1\tau^{-1}=a_0a_2a_3a_4\ .
$$
Thus, since $\tau$ commutes with $a_0$ and $a_2$, we have:
$$
\tau a_1\tau a_1=a_0a_2\ .
$$

\medskip
Now, we prove (ii). Applying Assertion 1 to (2.1), one has that $\MM(F_{g,r+1},\PP_n)$ is
generated by $a_0,\dots, a_{n+r}, b_1,\dots, b_{2g
-1},c,d_1,\dots, d_r, \tau_1, \dots, \tau_{n-1}$. But,
Assertion 2 implies
$$
a_{r+i}=\tau_{i-1}a_{r+i-1}\tau_{i-1}a_{r+i-1}
a_{r+i-2}^{-1}
$$
for $i=2,\dots, r$, thus $\MM(F_{g,r+1},\PP_n)$ is
generated by $a_0, \dots, a_{r+1}$, $b_1, \dots, b_{2g-1}$,
$c$, $d_1, \dots, d_r$, $\tau_1, \dots, \tau_{n-1}$. \qed

\medskip\noindent
{\bf Proposition 2.12}\qua {\sl {\rm(i)}\qua For $\Gamma$ equal to
$A_{l}$, $D_{l}$, $E_6$, or $E_7$, we denote by
$\rho_{PV}:A(\Gamma)\to \MM(F)$ the Perron-Vannier
representation of $A(\Gamma)$. In each case, $b_i$ denotes
the Dehn twist represented in the corresponding figure
(Figure 7, 8, or 9),
for $i=1,2,3$. Then:
\begin{eqnarray*}
\rho_{PV}(\Delta^2(A_{2p+1}))&=&b_1b_2,\\
\rho_{PV}(\Delta^4(A_{2p}))&=&b_1,\\
\rho_{PV}(\Delta^2(D_{2p+1}))&=&b_1b_2^{2p-1},\\
\rho_{PV}(\Delta(D_{2p}))&=&b_1b_2b_3^{p-1},\\
\rho_{PV}(\Delta^2(E_6))&=&b_1,\\
\rho_{PV}(\Delta(E_7))&=&b_1b_2^2.\\
\end{eqnarray*}

{\rm(ii)}\qua We denote by $\rho_{PV}: A(B_{l})\to
\MM(F,\{P_1,P_2\})$ the Perron-Vannier representation of
$A(B_{l})$. In each case, $b_i$ denotes the Dehn twist
represented in Figure \ref{f8}, for $i=1,2$. Then:
\begin{eqnarray*}
\rho_{PV}(\Delta(B_{2p}))&=&b_1b_2,\\
\rho_{PV}(\Delta^2(B_{2p+1}))&=&b_1.\\
\end{eqnarray*}

{\rm(iii)}\qua We denote by $\rho_G:A(B_{l})\to\MM( S^1\times
I,\PP_l)$ the graph representation of $A(B_{l})$ in the
punctured mapping class group of the cylinder. Let $b_1,
b_2$ denote the Dehn twists represented in Figure \ref{f9}. Then:
$$
\rho_G(\Delta(B_{l}))=b_1^{l-1}b_2\ .
$$}

Part (i) of Proposition 2.12 is proved in \cite{Mat}
with different techniques from the ones used in this
paper. Matsumoto's proof is based on the study of
geometric monodromies of simple singularities. Our proof
consists first on showing that the image of the considered
element lies in the center of the punctured mapping class
group, and, afterwards, on identifying this image using
the action of the center on some curves.

\medskip\noindent
{\bf Proof}\qua We only prove the equality
$$
\rho(\Delta(B_{2p}))=b_1b_2
$$
of Part (ii): the other equalities can be proved in the
same way.

By Proposition 2.10, $\MM(F,\{P_1,P_2\})$ is generated by
the Dehn twists $a_1,a_2,a_3$, $b_1$, $\sigma_2, \dots,
\sigma_{2p-1}$ and the braid twist $\tau_1$ represented in
Figure \ref{f8}. Since $\Delta(B_{2p})$ is in the center of
$A(B_{2p})$, $\rho_{PV}(\Delta(B_{2p}))$ commutes with
$\tau_1, \sigma_2, \dots, \sigma_{2p-1}$. The Dehn twist
$b_1$ belongs to the center of $\MM(F,\{P_1,P_2\})$, thus
$\rho_{PV}(\Delta(B_{2p}))$ also commutes with $b_1$. Let
$s_i$ be the defining circle of $a_i$, for $i=1,2,3$.
Using the expression of $\Delta(B_{2p})$ given in
Proposition 2.8, we verify that
$\rho_{PV}(\Delta(B_{2p}))(s_i)$ is isotopic to $s_i$,
thus $\rho_{PV}(\Delta(B_{2p}))$ commutes with $a_i$.

So, $\rho_{PV}(\Delta(B_{2p}))$ is an element of the center
of $\MM(F,\{P_1,P_2\})$. By \cite{PaRo}, this center is a
free abelian group of rank 2 generated by $b_1$ and $b_2$.
Thus $\rho_{PV}(\Delta(B_{2p}))=b_1^{q_1}b_2^{q_2}$ for
some $q_1,q_2\in\BZ$.

Now, consider the curve $\gamma$ of Figure \ref{f8}. Clearly, the
only element of the center of $\MM(F,\{P_1,P_2\})$ which
fixes $\gamma$ up to isotopy is the identity.
Using the expression of $\Delta(B_{2p})$ given in
Proposition 2.8, we verify that $\rho_{PV}(\Delta(B_{2p}))
b_1^{-1} b_2^{-1}$ fixes $\gamma$ up to isotopy, thus
$q_1=q_2=1$ and $\rho_{PV}(\Delta(B_{2p})) =b_1b_2$.
\qed


\subsection{Matsumoto's presentation for $\MM(F_{g,1})$ and
$\MM(F_{g,0})$}

This subparagraph is dedicated to the statement of
Matsumoto's presentations for $\MM(F_{g,1})$ and
$\MM(F_{g,0})$.

We first introduce some notation. Let $\Gamma$ be a Coxeter
graph, and let $X$ be a subset of the set $\{x_1,\dots,
x_{l}\}$ of vertices of $\Gamma$. Recall
that $\Gamma_X$ denotes the Coxeter subgraph
generated by $X$, and $A_X$ denotes the parabolic subgroup
of $A(\Gamma)$ generated by $X$. If $\Gamma_X$ is a finite
type connected Coxeter graph, then we denote by $\Delta(X)$
the fundamental element of $A_X$, viewed as an element of
$A(\Gamma)$.

\medskip\noindent
{\bf Theorem 2.13} (Matsumoto \cite{Mat}). {\sl Let $g\ge
1$, and let $\Gamma_g$ be the Coxeter graph drawn in Figure \ref{f14}. 

\smallskip
{\rm(i)}\qua $\MM(F_{g,1})$ is isomorphic with the quotient of
$A(\Gamma_g)$ by the following relations:
$$\matrix{
(1)\quad\hfill&\hfill\Delta^4(y_1,y_2,y_3,z)&=&\Delta^2(
x_0,y_1,y_2,y_3,z)\hfill&\quad{\hbox{\sl if}}\ g\ge 2,\hfill\cr
(2)\quad\hfill&\hfill\Delta^2(y_1,y_2,y_3,y_4,y_5,z)&=&
\Delta(x_0,y_1,y_2,y_3,y_4,y_5,z)\hfill&\quad{\hbox{\sl if}}\ 
g\ge 3.\hfill\cr}
$$

{\rm(ii)}\qua $\MM(F_{g,0})$ is isomorphic with the quotient of
$A(\Gamma_g)$ by the relations (1) and (2) above plus the
following relation:
$$
(3)\quad\matrix{
\hfill (x_0y_1)^6&=&1\hfill&\quad {\hbox{\sl if}}\ g=1,\hfill\cr
\hfill x_0^{2g-2}&=&\Delta^2(y_2,y_3,z,y_4,\dots, y_{2g-1})
\hfill&\quad{\hbox{\sl if}}\ g\ge 2.\hfill\cr}
$$}

\begin{figure}[ht]
\centerline{\input{pmc_14}}
\caption{\label{f14} Coxeter graph associated with $\MM(F_{g,1})$
and with $\MM(F_{g,0})$}
\end{figure}

\medskip
Set $r=n=0$, and consider the Dehn twists $a_0$, $b_1,\dots,
b_{2g-1}$, $c$ of Figure \ref{f10}. By Lemma 2.1, there is a well
defined homomorphism $\rho: A(\Gamma_g)\to \MM(F_{g,1})$
which sends $x_0$ on $a_0$, $y_i$ on $b_i$ for $i=1,\dots,
{2g-1}$, and $z$ on $c$. By \cite{Hum} (see Proposition
2.10), this homomorphism is surjective. By Proposition 2.12,
both $\rho(\Delta^4(y_1,y_2,y_3,z))$ and
$\rho(\Delta^2(x_0,y_1,y_2,y_3,z))$ are equal to the Dehn
twist $\sigma_1$ of Figure \ref{f15}. Similarly, both
$\rho(\Delta^2(y_1,\dots,y_5,z))$ and
$\rho(\Delta(x_0,y_1,\dots,y_5,z))$ are equal to the Dehn
twist $\sigma_2$ of Figure \ref{f15}. 
Let $G_g$ denote the quotient of $A(\Gamma_g)$ by
the relations (1) and (2). So, the homomorphism
$\rho:A(\Gamma_g)\to\MM(F_{g,1})$ induces a surjective
homomorphism $\bar\rho:G_g\to\MM(F_{g,1})$. In order to
prove that this homomorphism is in fact an
isomorphism, Matsumoto \cite{Mat} showed that
the presentation of $G_g$ as a quotient of $A(\Gamma_g)$
is equivalent to Wajnryb's presentation of 
$\MM(F_{g,1})$ \cite{Waj}.

Similar remarks can be made for the presentation of
$\MM(F_{g,0})$.

\begin{figure}[ht]
\centerline{\input{pmc_15}}
\caption{\label{f15} Relations in $\MM(F_{g,1})$}
\end{figure}


\section{The presentation}

Recall that, if $\Gamma$ is a finite type connected Coxeter
graph, then $\Delta(\Gamma)$ denotes 
the fundamental element of
$A(\Gamma)$. If
$\Gamma$ is any Coxeter graph and $X$ is a subset of the
set $\{x_1,\dots, x_l\}$ of vertices of $\Gamma$ such that
$\Gamma_X$ is finite type and connected, then we denote by
$\Delta(X)$ the fundamental element of $A_X=A(\Gamma_X)$
viewed as an element of $A(\Gamma)$.

\medskip\noindent
{\bf Theorem 3.1}\qua {\sl Let $g\ge 1$, let $r,n\ge 0$, and let
$\Gamma_{g,r,n}$ be the Coxeter graph drawn in Figure \ref{f16}.
Then $\MM(F_{g,r+1},\PP_n)$ is isomorphic with the quotient
of $A(\Gamma_{g,r,n})$ by the following relations.

\smallskip
$\bullet$\qua Relations from $\MM(F_{g,1})$:

\medskip
\centerline{\vbox{\halign{
$#$\quad\hfill&\hfill$#$&\ $#$\ &$#$\quad\hfill&$#$\hfill\cr
{\rm(R1)}&\Delta^4(y_1,y_2,y_3,z)&=&\Delta^2
(x_0,y_1,y_2,y_3,z)&{\hbox{\sl if}}\ g\ge 2,\cr
\noalign{\smallskip}
{\rm(R2)}&\Delta^2(y_1,y_2,y_3,y_4,y_5,z)&=&
\Delta(x_0,y_1,y_2,y_3,y_4,y_5,z)&{\hbox{\sl if}}\
g\ge 3.\cr}}}

\medskip
$\bullet$\qua Relations of commutation:

\medskip
\centerline{\vbox{\halign{
$#$\quad\hfill&$#$\ &$#$\quad\hfill&$#$\hfill\cr
{\rm(R3)}&&x_k\Delta^{-1}(x_{i+1},x_j,y_1)x_i
\Delta(x_{i+1},x_j,y_1)\cr
&=&\Delta^{-1}(x_{i+1},x_j,y_1)x_i
\Delta(x_{i+1},x_j,y_1) x_k&{\hbox{\sl if}}\ 0\le
k<j<i\le r,\cr
\noalign{\smallskip}
{\rm(R4)}&& y_2\Delta^{-1}(x_{i+1},x_j,y_1)x_i
\Delta(x_{i+1},x_j,y_1)\cr
&=&\Delta^{-1}(x_{i+1},x_j,y_1)x_i
\Delta(x_{i+1},x_j,y_1) y_2&{\rm if}\ 0\le
j<i\le r\ {\hbox{\sl and}}\ g\ge 2,\cr}}}

\eject
$\bullet$\qua Expressions of the $u_i$'s:

\medskip
\centerline{\vbox{\halign{
$#$\hfill&\hfill$#$&\ $#$\ &$#$\hfill&$#$\hfill\cr
{\rm(R5)}&u_1&=&\Delta(x_0,x_1,y_1,y_2,y_3,z)
\Delta^{-2}(x_1,y_1,y_2,y_3,z)&{\hbox{\sl if}}\ g\ge
2,\cr
\noalign{\smallskip}
{\rm(R6)}& u_{i+1}&=&\Delta(x_i,x_{i+1},y_1,
y_2,y_3,z)\Delta^{-2}(x_{i+1},y_1,y_2,y_3,z)\cr
&&&\quad\Delta^2
(x_0,x_{i+1},y_1)\Delta^{-1}(x_0,x_i,x_{i+1},y_1)\ \
{\hbox{\sl if}}\ 1\le i\le r-1 &, \ g\ge 2.
\cr}}}

\medskip
$\bullet$\qua Other relations:}

\medskip
\centerline{\vbox{\halign{
$#$\quad\hfill&\hfill$#$&\ $#$\ &$#$\quad\hfill&$#$\hfill\cr
{\rm(R7)}&\Delta(x_r,x_{r+1},y_1,v_1)&=&
\Delta^2(x_{r+1},y_1,v_1)&{\hbox{\sl if}}\ n\ge 2,\cr
\noalign{\smallskip}
{\rm(R8a)}&\Delta(x_0,x_1,y_1,y_2,y_3,z)&=&
\Delta^2(x_1,y_1,y_2,y_3,z)&{\hbox{\sl if}}\ n\ge
1, g\ge 2, \ r=0,\cr
\noalign{\smallskip}
{\rm(R8b)}&\multispan 4 $\phantom{=\ }
\Delta(x_r,x_{r+1},y_1,y_2,y_3,z)
\Delta^{-2}(x_{r+1},y_1,y_2,y_3,z)$\hfill\cr
&\multispan 3 $=\ 
\Delta(x_0,x_r,x_{r+1},
y_1)\Delta^{-2}(x_0,x_{r+1},y_1)$\hfill\quad
&{\hbox{\sl if}}\ n\ge 1, g\ge 2,\ r\ge 1.\cr
}}}

\begin{figure}[ht]
\centerline{\input{pmc_16}}
\caption{\label{f16} Coxeter graph associated with
$\MM(F_{g,r+1},\PP_n)$}
\end{figure}

\medskip
Notice that only the relations (R1),
(R2), (R7), and (R8a) remain in the presentation of
$\MM(F_{g,1},\PP_n)$,
and (R8a) has to be replaced by (R8b) if $r\ge 1$.

Assume that $g\ge 2$. From the relations (R5) and (R6) we
see that we can remove $u_1,\dots, u_r$ from the generating set.
However, to do so, one has to add relations comming
from the ones in the Artin group $A(\Gamma_{g,r,n})$. 
For example, one has that
$\Delta(x_0,x_1,y_1,y_2,y_3,z) 
\Delta^{-2}(x_1,y_1,y_2,y_3,z)$ commutes with $y_4$ in the
quotient, since $u_1$ commutes with $y_4$ in
$A(\Gamma_{g,r,n})$.

Consider the Dehn twists $a_0,\dots, a_{r+1}$, $b_1,\dots,
b_{2g-1}$, $c$, $d_1,\dots, d_r$ and the braid twists
$\tau_1,\dots, \tau_{n-1}$ represented in Figure \ref{f10}. From
Subsection 2.1 follows that there is a well defined
homomorphism $\rho: A(\Gamma_{g,r,n})\to
\MM(F_{g,r+1},\PP_n)$ which sends $x_i$ on $a_i$ for
$i=0,\dots, r+1$, $y_i$ on $b_i$ for $i=1,\dots, 2g-1$, $z$
on $c$, $u_i$ on $d_i$ for $i=1,\dots, r$, and $v_i$ on
$\tau_i$ for $i=1,\dots, n-1$. This homomorphism is
surjective by Proposition 2.10. 
If $w_1=w_2$ is one of the relations 
(R1),\dots,(R7), (R8a), (R8b), then we have $\rho(w_1) = 
\rho(w_2)$. This fact can be easily proved using 
Proposition 2.12 in the case of the relations (R1), (R2), 
(R5), (R6), (R7), (R8a), and (R8b), and comes from the 
following reason in the case of the relations (R3) and 
(R4). We have the equality 
$$
\Delta^{-1} (x_{i+1}, x_j, y_1) x_i \Delta (x_{i+1}, x_j, 
y_1) = y_1^{-1} x_{i+1}^{-1} x_j^{-1} y_1^{-1} x_i y_1 x_j 
x_{i+1} y_1,
$$
and the image by $b_1^{-1} a_{i+1}^{-1} a_j^{-1} b_1^{-1}$ of 
the defining circle of $a_i$ is disjoint from the defining 
circle of $a_k$, up to isotopy, if $k<j$, and is disjoint 
from the defining circle of $b_2$, up to isotopy.

Let $G(g,r,n)$ denote the quotient of $A(\Gamma_{g,r,n})$
by the relations (R1),..,(R7), (R8a), (R8b). By the above considerations,
the homomorphism : $$\rho:A(\Gamma_{g,r,n})
\to  \MM(F_{g,r+1},\PP_n)$$ induces a surjective homomorphism
$\bar\rho: G(g,r,n)\to \MM(F_{g,r+1},\PP_n)$. In order to 
prove Theorem 3.1, it remains to show that this
homomorphism is in fact an isomorphism. This will be the
object of Subsection 3.1.

\medskip\noindent
{\bf Theorem 3.2}\qua {\sl Let $g\ge 1$, let $n\ge 1$,
and let $\Gamma_{g,0,n}$ be the Coxeter graph drawn in
Figure \ref{f16}. Then
$\MM(F_{g,0},\PP_n)$ is isomorphic with the quotient of
$A(\Gamma_{g,0,n})$ by the following relations.

\smallskip
$\bullet$\qua Relations from $\MM(F_{g,1},\PP_n)$:

\medskip
\centerline{\vbox{\halign{
$#$\quad\hfill&\hfill$#$&\ $#$\ &$#$\quad\hfill&$#$\hfill\cr
{\rm(R1)}&\Delta^4(y_1,y_2,y_3,z)&=&\Delta^2
(x_0,y_1,y_2,y_3,z)&{\hbox{\sl if}}\ g\ge 2,\cr
\noalign{\smallskip}
{\rm(R2)}&\Delta^2(y_1,y_2,y_3,y_4,y_5,z)&=&
\Delta(x_0,y_1,y_2,y_3,y_4,y_5,z)&{\hbox{\sl if}}\
g\ge 3,\cr
\noalign{\smallskip}
{\rm(R7)}&\Delta(x_0,x_1,y_1,v_1)&=&
\Delta^2(x_1,y_1,v_1)&{\hbox{\sl if}}\ n\ge 2,
\cr
\noalign{\smallskip}
{\rm(R8a)}&\Delta(x_0,x_1,y_1,y_2,y_3,z)&=&
\Delta^2(x_1,y_1,y_2,y_3,z)&{\hbox{\sl if}}\ n\ge
1\ {\hbox{\sl and}}\ g\ge 2.\cr}}}

\medskip
$\bullet$\qua Other relations:}

\medskip
\centerline{\vbox{\halign{
$#$\quad\hfill&\hfill$#$&\ $#$\ &$#$\quad\hfill&$#$\hfill\cr
{\rm(R9a)}& x_0^{2g-n-2}\Delta(x_1,v_1,\dots,
v_{n-1})&=&\Delta^2(z,y_2,\dots,y_{2g-1})&{\hbox{\sl
if}}\ g\ge 2,\cr
\noalign{\smallskip}
{\rm(R9b)}& x_0^n&=&\Delta(x_1,v_1,\dots, 
v_{n-1})&{\hbox{\sl if}}\ g=1,\cr
\noalign{\smallskip}
{\rm(R9c)}&\Delta^4 (x_0,y_1)&=&\Delta^2(v_1,\dots, 
v_{n-1})&{\hbox{\sl if}}\ g=1.\cr}}}

\medskip
Note that, in the above presentation, the relation (R9a),
which holds if $g\ge 2$, has to be replaced by the
relations (R9b) and (R9c) when $g=1$.

Consider the Dehn twists $a_0,a_1$, $b_1,\dots, b_{2g-1}$,
$c$ and the braid twists $\tau_1,.., \tau_{n-1}$
represented in Figure \ref{f11}. From Subsection 2.1 follows that
there is a well defined homomorphism
$\rho_0:A(\Gamma_{g,0,n}) \to \MM(F_{g,0},\PP_n)$ which
sends $x_i$ on $a_i$ for $i=0,1$, $y_i$ on $b_i$ for
$i=1,\dots, 2g-1$, $z$ on $c$, and $v_i$ on $\tau_i$ for
$i=1,\dots, n-1$. This homomorphism is surjective by
Corollary 2.11. 
Let $G_0(g,n)$ denote the quotient of $A(\Gamma_{g,0,n})$ by
the relations
(R1), (R2), (R7), (R8), (R9a), (R9b), and (R9c).
As before, using Proposition 2.12, one can easily
prove that the homomorphism $\rho_0:
A(\Gamma_{g,0,n})\to \MM(F_{g,0},\PP_n)$ induces a surjective
homomorphism $\bar\rho_0: G_0(g,n)\to \MM(F_{g,0},\PP_n)$.
In order to prove Theorem 3.2, it remains to show that this
homomorphism is in fact an isomorphism. This will be the
object of Subsection 3.2.


\subsection{Proof of Theorem 3.1}

The proof of Theorem 3.1 is organized as follows. In the
first step,
starting from Matsumoto's presentation of $\MM(F_{g,1})$
\cite{Mat} (see Theorem 2.13), we determine by induction on
$n$ a presentation of $\PM(F_{g,1},\PP_n)$ (Proposition
3.3), applying Lemma 2.5 to the exact sequence (2.2) of
Subsection 2.2. In the second step, we determine a presentation of
$\PM(F_{g,r+1},\PP_n)$ (Proposition 3.7), applying Lemma 2.5
to the exact sequence (2.3). Finally, we prove Theorem 3.1
applying Lemma 2.5 to the exact sequence (2.1).

\medskip\noindent
{\bf Proposition 3.3}\qua {\sl Let $g\ge 1$, let $n\ge 0$, and
let $P\Gamma_{g,0,n}$ be the Coxeter graph drawn in Figure \ref{f17}.
Then $\PM(F_{g,1},\PP_n)$ is isomorphic with the
quotient of $A(P\Gamma_{g,0,n})$ by the following
relations.

\smallskip
$\bullet$\qua Relations from $\MM(F_{g,1})$:

\medskip
\centerline{\vbox{\halign{
$#$\quad\hfill&\hfill$#$&\ $#$\ &$#$\quad\hfill&$#$\hfill
\cr
{\rm(PR1)}&\Delta^4(y_1,y_2,y_3,z)&=&\Delta^2(x_0,y_1,y_2,y_3,z)
&{\hbox{\sl if}}\ g\ge 2,\cr
\noalign{\smallskip}
{\rm(PR2)}&\Delta^2(y_1,y_2,y_3,y_4,y_5,z)&=&\Delta(x_0,y_1,y_2,
y_3,y_4,y_5,z)&{\hbox{\sl if}}\ g\ge 3.\cr}}}

\medskip
$\bullet$\qua Relations of commutation:

\medskip
\centerline{\vbox{\halign{
$#$\quad\hfill&\ $#$\ &$#$\quad\hfill&$#$\hfill\cr
{\rm(PR3)}&&x_k\Delta^{-1}(x_{i+1},x_j,y_1)x_i\Delta(x_{i+1},
x_j,y_1)\cr
&=&\Delta^{-1}(x_{i+1},x_j,y_1)x_i\Delta(x_{i+1},
x_j,y_1)x_k&{\hbox{\sl if}}\ 0\le k<j<i\le n-1,\cr
\noalign{\smallskip}
{\rm(PR4)}&&y_2\Delta^{-1}(x_{i+1},x_j,y_1)x_i\Delta(x_{i+1},
x_j,y_1)\cr
&=&\Delta^{-1}(x_{i+1},x_j,y_1)x_i\Delta(x_{i+1},
x_j,y_1)y_2 & {\hbox{\sl if}}\ 0\le j<i\le n-1, \ 
g\ge 2.\cr}}}

\medskip
$\bullet$\qua Relations between fundamental elements:}

\medskip
\centerline{\vbox{\halign{
$#$\hfill&$#$\hfill&$#$\hfill\cr
{\rm(PR5)}&\quad \Delta(x_0,x_1,y_1,y_2,y_3,z)\ =\ \Delta^2(x_1,y_1,
y_2,y_3,z)&{\hbox{\sl if}}\ g\ge 2,\cr
\noalign{\smallskip}
{\rm(PR6)}&\phantom{=\ } \Delta(x_i,x_{i+1},y_1,y_2,y_3,z)
\Delta^{-2}(x_{i+1},y_1,y_2,y_3,z)\cr
&=\ \Delta(x_0,x_i,x_{i+1},y_1)\Delta^{-2}(x_0,x_{i+1},y_1)
& {\hbox{\sl if}}\ 1\le i\le n-1, \ g\ge 2.\cr}}}

\begin{figure}[ht]
\centerline{\input{pmc_17}}
\caption{\label{f17} Coxeter graph associated with 
$\PM(F_{g,1},\PP_n)$}
\end{figure}

\medskip
The following lemmas 3.4, 3.5, and 3.6 are preliminary
results to the proof of Proposition 3.3.

\medskip\noindent
{\bf Lemma 3.4}\qua {\sl Let $\Gamma$ be the Coxeter graph
drawn in Figure \ref{f18}, and let $G$ be the quotient of
$A(\Gamma)$ by the following relation:
$$
x_4\Delta^{-1}(x_1,x_3,y)x_2\Delta(x_1,x_3,y)=\Delta^{-1}
(x_1,x_3,y)x_2\Delta(x_1,x_3,y)x_4\ .
$$
Then the following equalities hold in $G$.}

\medskip
\centerline{\vbox{\halign{
\hfill$#$&\ $#$\ &$#$\hfill\cr
x_3\Delta^{-1}(x_2,x_4,y)x_1\Delta(x_2,x_4,y)&=&
\Delta^{-1}(x_2,x_4,y)x_1\Delta(x_2,x_4,y)x_3,\cr
x_2\Delta^{-1}(x_1,x_3,y)x_4\Delta(x_1,x_3,y)&=&
\Delta^{-1}(x_1,x_3,y)x_4\Delta(x_1,x_3,y)x_2,\cr
x_1\Delta^{-1}(x_2,x_4,y)x_3\Delta(x_2,x_4,y)&=&
\Delta^{-1}(x_2,x_4,y)x_3\Delta(x_2,x_4,y)x_1.\cr}}}

\begin{figure}[ht]
\centerline{\input{pmc_18}}
\nocolon
\caption{\label{f18}}
\end{figure}

\medskip\noindent
{\bf Proof}\qua It clearly suffices to prove the first
equality.

\medskip
\halign{
$#$\ &$#$\hfill\cr
&x_3\Delta^{-1}(x_2,x_4,y)x_1\Delta(x_2,x_4,y)x_3^{-1}
\Delta^{-1}(x_2,x_4,y)x_1^{-1}\Delta(x_2,x_4,y)\cr
=&x_3y^{-1}x_2^{-1}x_4^{-1}y^{-1}x_1yx_2x_4yx_3^{-1}
y^{-1}x_2^{-1}x_4^{-1}y^{-1}x_1^{-1}yx_2x_4y\cr
=&y^{-1}\cdot x_3^{-1}yx_3x_2^{-1}x_4^{-1} x_1yx_1^{-1}
x_2x_4x_3^{-1}y^{-1}x_3x_2^{-1}x_4^{-1}x_1y^{-1}x_1^{-1}
x_2x_4\cdot y\cr
=&y^{-1} x_2^{-1} x_3^{-1} \cdot x_2 y x_2^{-1} x_1 x_3
x_4^{-1} y x_4 x_1^{-1} x_3^{-1} x_2 y^{-1} x_2^{-1} x_1
x_3 x_4^{-1} y^{-1} x_4 x_1^{-1} x_3^{-1} \cr
&\cdot x_3 x_2y\cr
=& y^{-1} x_2^{-1} x_3^{-1} \cdot y^{-1} x_2 y x_1 x_3 y
x_4 y^{-1} x_1^{-1} x_3^{-1} y^{-1} x_2^{-1} y x_1 x_3 y
x_4^{-1} y^{-1} x_1^{-1} x_3^{-1} \cdot x_3 x_2 y\cr
=& y^{-1} x_2^{-1} x_3^{-1} y^{-1} \cdot x_2 \Delta (x_1,
x_3, y) x_4 \Delta^{-1} (x_1, x_3, y) x_2^{-1} \Delta
(x_1, x_3, y) x_4^{-1} \cr
&\Delta^{-1} (x_1, x_3, y) \cdot y
x_3 x_2 y\cr
=& 1.\cr}\vglue -2em\noindent\hbox{}\hfill\qed

\medskip\noindent
{\bf Lemma 3.5}\qua {\sl We number the vertices of the Coxeter
graph $D_l$ according to Figure \ref{f6}. Then the following equalities
hold in $A(D_l)$.}

\medskip
\centerline{\vbox{\halign{
$#$\ &$#$\hfill\cr
&\Delta^{-1} (x_2, \dots, x_{l-1}) x_1^{-1} x_2 \Delta (x_2,
\dots, x_{l-1}) \Delta^{-1} (x_2, \dots, x_l) x_2^{-1} x_1
\Delta (x_2, \dots, x_l)\cr
=&x_l \Delta^{-1} (x_2, \dots, x_{l-1}) x_1^{-1} x_2 \Delta
(x_2, \dots, x_{l-1}) x_l^{-1},\cr
\cr
&\Delta^{-1} (x_2, \dots, x_l) x_2^{-1} x_1 \Delta (x_2,
\dots, x_l) \Delta^{-1} (x_2, \dots, x_{l-1}) x_2^{-1} x_1
\Delta (x_2, \dots, x_{l-1})\cr
=&x_{l-1} \Delta^{-1} (x_2, \dots, x_l) x_2^{-1} x_1 \Delta
(x_2, \dots, x_l) x_{l-1}^{-1}.\cr}}}

\medskip\noindent
{\bf Proof}\qua

\medskip
\halign{
$#$\ &$#$\hfill\cr
&x_l^{-1} \Delta^{-1} (x_2, \dots, x_{l-1}) x_1^{-1} x_2
\Delta (x_2, \dots, x_{l-1}) \Delta^{-1} (x_2, \dots, x_l)
x_2^{-1} x_1 \cr 
&\Delta (x_2, \dots, x_l) x_l
\quad \Delta^{-1} (x_2, \dots, x_{l-1}) x_2^{-1} x_1
\Delta (x_2, \dots, x_{l-1}) \cr
=&x_l^{-1} \Delta^{-1} (x_2, \dots, x_{l-2}) (x_{l-1}^{-1}
\dots x_2^{-1})x_2 x_1^{-1} (x_l^{-1} \dots x_2^{-1})  
x_2^{-1} x_1 x_2 \Delta (x_2, \dots, x_l)\cr 
&\quad \Delta^{-1} (x_2, \dots, x_{l-1}) x_1 x_2^{-1} (x_2
\dots x_{l-1}) \Delta (x_2, \dots, x_{l-2}) \cr
=&\Delta^{-1} (x_2, \dots, x_{l-2}) x_l^{-1} (x_{l-1}^{-1}
\dots x_3^{-1}) x_1^{-1} (x_l^{-1} \dots x_2^{-1}) x_1
(x_2 \dots x_l) x_1 \cr
&(x_3 \dots x_{l-1}) \Delta (x_2,
\dots, x_{l-2}) \cr
=&\Delta^{-1} (x_2, \dots, x_{l-2}) (x_l^{-1} \dots 
x_3^{-1}) x_1^{-1} (x_l^{-1} \dots x_3^{-1}) (x_3 \dots
x_l) x_1 (x_3 \dots x_l) \cr
&\Delta (x_2, \dots, x_{l-2}) \cr
=&1.\cr}

\medskip
\halign{
$#$\ &$#$\hfill\cr
&\Delta^{-1} (x_2, \dots, x_l) x_2^{-1} x_1 \Delta (x_2,
\dots, x_l) \Delta^{-1} (x_2, \dots, x_{l-1}) x_2^{-1} x_1
\Delta (x_2, \dots, x_{l-1})  \cr 
&x_{l-1}\Delta^{-1} (x_2, \dots, x_l) x_1^{-1} x_2 \Delta
(x_2, \dots, x_l) x_{l-1}^{-1}\cr
=&\Delta^{-1} (x_2, \dots, x_l) x_2^{-1} x_1 (x_2 \dots
x_l) x_2^{-1} x_1 x_2 \Delta (x_2, \dots, x_{l-1})
\Delta^{-1} (x_2, \dots, x_l) x_1^{-1} \cr 
&x_2 x_3^{-1} \quad \Delta (x_2, \dots, x_l) \cr
=&\Delta^{-1} (x_2, \dots, x_l) x_1 (x_3 \dots x_l) x_1
(x_l^{-1} \dots x_2^{-1}) x_1^{-1} x_2 x_3^{-1} \Delta
(x_2, \dots, x_l) \cr
=&\Delta^{-1} (x_2, \dots, x_l) x_3 x_1 (x_3 \dots x_l)
(x_l^{-1} \dots x_3^{-1}) x_1^{-1} x_3^{-1} \Delta (x_2,
\dots, x_l) \cr
=&1. \cr}\vglue -2em\noindent\hbox{}\hfill\qed

\medskip
Several algorithms to solve the word problem in finite type
Artin groups are known (see \cite{BrSa}, \cite{Del},
\cite{Cha}, \cite{DePa}). We use the one of \cite{DePa}
implemented in a Maple program to prove the following.

\medskip\noindent
{\bf Lemma 3.6}\qua {\sl {\rm(i)}\qua We number the vertices of $D_6$
 according to Figure \ref{f6}. Let

\medskip
\centerline{\vbox{\halign{
\hfill$#$&\ $#$\ &$#$\hfill\cr
w_1 &=& \Delta^{-1} (x_1, x_3) x_1^{-1} x_2 \Delta (x_1,
x_3) \cr
w_2 &=& \Delta^{-1} (x_1, x_3, x_4) x_1^{-1} x_2 \Delta
(x_1, x_3, x_4) \cr
w_3 &=& \Delta^{-1} (x_1, x_3, x_4, x_5) x_1^{-1} x_2
\Delta (x_1, x_3, x_4, x_5) \cr}}}

\medskip\noindent
Then the following equality holds in $A(D_6)$.
$$
x_2^{-1} x_1 w_1^{-1} w_2^{-1} w_3^{-1} x_6 w_3 x_6^{-1} 
w_1 = \Delta^{-2} (x_2, x_3, \dots, x_6) \Delta (x_1, x_2, 
x_3, \dots, x_6).
$$

{\rm(ii)}\qua We number the vertices of $D_4$ according to Figure \ref{f6}.
Let
$$
w=x_2^{-1} \Delta^{-1} (x_1, x_3, x_4) x_1^{-1} x_2 \Delta
(x_1, x_3, x_4) x_2.
$$
Then the following equality holds in $A(D_4)$.
$$
x_1^{-1} x_2 w = \Delta^{-2} (x_1, x_3, x_4) \Delta (x_1,
x_2, x_3, x_4). \eqno{\qed}
$$}

\noindent
{\bf Proof of Proposition 3.3}\qua We set $r=0$ and we
consider the Dehn twists $a_0, \dots, a_n$ $b_1, \dots,
b_{2g-1}$, $c$ represented in Figure \ref{f10}. From Subsection
2.1 follows that there is a well defined homomorphism
$\rho: A(P\Gamma_{g,0,n}) \to \PM( F_{g,1}, \PP_n)$ which
sends $x_i$ on $a_i$ for $i=0,\dots, n$, $y_i$ on $b_i$
for $i=1,\dots, 2g-1$, and $z$ on $c$. This homomorphism
is surjective by Proposition 2.10.
Let $PG(g,0,n)$ denote the quotient of $A (P\Gamma_{g,0,n})$
by the relations (PR1),\dots,(PR6). 
One can easily prove using Proposition 2.12 that: if 
$w_1=w_2$ is one of the relations (PR1),\dots,(PR6), then 
$\rho(w_1) = \rho(w_2)$. So,
the homomorphism $\rho: A( P\Gamma_{g,0,n}) \to \PM
(F_{g,1}, \PP_n)$ induces a surjective homomorphism :
$$\bar\rho: PG(g,0,n) \to \PM (F_{g,1}, \PP_n).$$
Now, we prove by induction on $n$ that $\bar\rho$ is an
isomorphism. The case $n=0$ is proved in \cite{Mat} (see
Theorem 2.13). So, we assume that $n>0$. By the inductive
hypothesis, $\PM (F_{g,1}, \PP_{n-1})$ is isomorphic with $PG
(g,0,n-1)$. On the other hand, $\pi_1 (F_{g,1} \setminus
\PP_{n-1}, P_n)$ is the free group $F(\alpha_1, \dots,
\alpha_n, \beta_1, \dots, \beta_{2g-1})$ freely generated
by the loops $\alpha_1, \dots, \alpha_n$, $\beta_1, \dots,
\beta_{2g-1}$ represented in Figure \ref{f12}. Applying Lemma 2.5
to the exact sequence (2.2) of Subsection 2.2, one has
that $\PM (F_{g,1}, \PP_n)$ is isomorphic with the
quotient of the free product $PG (g,0,n-1) \ast F(
\alpha_1, \dots, \alpha_n, \beta_1, \dots, \beta_{2g-1})$
by the following relations.

\smallskip
$\bullet$\qua Relations involving the $\alpha_i$'s:

\medskip
\centerline{\vbox{\halign{
$#$\quad\hfill&\hfill$#$&\ $#$\ &$#$\quad\hfill&$#$\hfill
\cr
{\rm(PT1)} &x_j \alpha_i x_j^{-1} &=& \alpha_i &{\rm for}\ 0\le
j<i \le n,\cr
{\rm(PT2)} &x_j \alpha_i x_j^{-1} &=& \alpha_{j+1}^{-1} \alpha_i
\alpha_{j+1} &{\rm for}\ 1\le i\le j\le n-1,\cr
{\rm(PT3)} &y_1 \alpha_i y_1^{-1} &=& \beta_1^{-1} \alpha_i &{\rm
for}\ 1\le i\le n,\cr
{\rm(PT4)} &y_j \alpha_i y_j^{-1} &=& \alpha_i &{\rm for}\ 1\le
i\le n\ {\rm and}\ 2\le j\le 2g-1,\cr
{\rm(PT5)} &z \alpha_i z^{-1} &=& \alpha_i &{\rm for}\ 1\le i\le
n.\cr}}}

\medskip
$\bullet$\qua Relations involving the $\beta_i$'s:

\medskip
\centerline{\vbox{\halign{
$#$\quad\hfill&\hfill$#$&\ $#$\ &$#$\quad\hfill&$#$\hfill
\cr
{\rm(PT6)} &x_j \beta_1 x_j^{-1} &=& \beta_1 \alpha_{j+1} &{\rm
for}\ 0\le j\le n-1,\cr
{\rm(PT7)} &x_j \beta_i x_j^{-1} &=& \beta_i &{\rm for}\ 0\le
j\le n-1\ {\rm and}\ 2\le i\le 2g-1,\cr
{\rm(PT8)} &y_j \beta_i y_j^{-1} &=& \beta_i &{\rm for}\ j\neq 
i-1\ {\rm and}\ j\neq i+1,\cr
{\rm(PT9)} &y_{i-1} \beta_i y_{i-1}^{-1} &=& \beta_i \beta_{i-1}
&{\rm for}\ 2\le i\le 2g-1,\cr
{\rm(PT10)} &y_{i+1} \beta_i y_{i+1}^{-1} &=& \beta_{i+1}^{-1}
\beta_i &{\rm for}\ 1\le i\le 2g-2,\cr
{\rm(PT11)} &z \beta_3 z^{-1} &=& \beta_3 \beta_2 \beta_1
\alpha_1 \beta_1^{-1},\cr
{\rm(PT12)} &z \beta_i z^{-1} &=& \beta_i &{\rm for}\ i\neq
3.\cr}}}

\medskip
Consider the homomorphism $f: PG (g,0,n-1) \ast F
(\alpha_1, \dots, \alpha_n, \beta_1, \dots, \beta_{2g-1})
\to PG(g,0,n)$ defined by:

\medskip
\centerline{\vbox{\halign{
\hfill$#$&\,$#$\,&$#$\quad\hfill&$#$\hfill&\ $#$\hfill\cr
f(x_i) &=& x_i &{\rm for}\ 0\le i\le n-1,\cr
f(y_i) &=& y_i &{\rm for}\ 1\le i\le 2g-1,\cr
f(z) &=& z,\cr
f(\alpha_i) &=&\multispan 2$\!x_{n-1}^{-1} \Delta^{-1} (x_n,
x_{i-1}, y_1) x_n^{-1} x_{n-1} \Delta (x_n, x_{i-1}, y_1)
x_{n-1}$ \hfill &{\rm for}\ 1\le i\le n-1,\cr
f(\alpha_n) &=& x_{n-1}^{-1} x_n,\cr
f(\beta_i) &=&\multispan 2$\!\Delta^{-1} (x_{n-1}, y_1, 
\dots, y_i) x_{n-1}^{-1} x_n \Delta (x_{n-1}, y_1, \dots,
y_i)$ \hfill &{\rm for}\ 1\le i\le 2g-1.\cr}}}

\medskip\noindent
{\bf Assertion 1}\qua {\sl $f$ induces a homomorphism $\bar f:
\PM (F_{g,1}, \PP_n) \to PG (g,0,n)$.}

\medskip
One can easily verify on the generators of $PG (g,0,n)$
that $\bar f \circ \bar\rho$ is the identity of $PG
(g,0,n)$. So, Assertion 1 shows that $\bar \rho$ is
injective and, therefore, finishes the proof of
Proposition 3.3.

\medskip\noindent
{\bf Proof of Assertion 1}\qua We have to show that: if $w_1
=w_2$ is one of the relations (PT1),\dots,(PT12), then
$f(w_1) = f(w_2)$.

By an {\it easy case} we mean a relation $w_1=w_2$ such
that the equality $f(w_1)=f(w_2)$ in $PG(g,0,n)$ is a
direct consequence of the braid relations in $A
(P\Gamma_{g,0,n})$. For instance, (PT5), (PT6), and (PT8)
are easy cases.

$\bullet$\qua Relation (PT1): (PT1) is an easy case if either
$j=i-1$ or $i=n$. So, we assume that $0\le j< i-1< n-1$.
Then:

\medskip
\halign{
$#$\ &$#$\hfill\cr
&f(x_j \alpha_i x_j^{-1}) f( \alpha_i)^{-1}\cr
=& x_j x_{n-1}^{-1} \Delta^{-1} (x_n, x_{i-1}, y_1) x_n^{-1}
x_{n-1} \Delta (x_n, x_{i-1}, y_1) x_{n-1} x_j^{-1} 
x_{n-1}^{-1}  \cr
&\Delta^{-1} (x_n, x_{i-1}, y_1)
x_{n-1}^{-1}x_n \Delta (x_n, x_{i-1} ,y_1) x_{n-1}\cr
=& x_{n-1}^{-1} x_{i-1}^{-1} \cdot x_j \Delta^{-1} (x_n,
x_{i-1}, y_1) x_{n-1} \Delta (x_n, x_{i-1}, y_1) x_j^{-1}
\Delta^{-1} (x_n, x_{i-1}, y_1) x_{n-1}^{-1}\cr
& \Delta (x_n,x_{i-1} ,y_1)  \cdot x_{i-1} x_{n-1}\cr
=&1 \quad({\rm by\ (PR3)}).\cr}

\medskip
$\bullet$\qua Relation (PT2): (PT2) is an easy case if $j=n-1$.
So, we assume that $j<n-1$. Then:

\medskip
\halign{
$#$\ &$#$\hfill\cr
&f (x_j \alpha_i x_j^{-1}) f (\alpha_{j+1}^{-1} \alpha_i
\alpha_{j+1})^{-1}\cr
=& x_j x_{n-1}^{-1} \Delta^{-1} (x_n, x_{i-1}, y_1) x_n^{-1}
x_{n-1} \Delta (x_n, x_{i-1}, y_1) x_{n-1} x_j^{-1} 
x_{n-1}^{-1} \Delta^{-1} (x_n, x_j, y_1) \cr
&x_{n-1}^{-1} x_n \Delta (x_n, x_j, y_1) x_{n-1} x_{n-1}^{-1} 
\Delta^{-1} (x_n, x_{i-1}, y_1) x_{n-1}^{-1} x_n \Delta
(x_n, x_{i-1}, y_1) x_{n-1}\cr
&  x_{n-1}^{-1}\Delta^{-1} (x_n,x_j, y_1)x_n^{-1} x_{n-1} 
\Delta (x_n, x_j, y_1) x_{n-1}\cr
=& x_j x_{n-1}^{-1} x_{i-1}^{-1} \Delta^{-1} (x_n, x_{i-1},
y_1) x_{n-1} \Delta (x_n, x_{i-1}, y_1) \Delta^{-1} (x_n,
x_j, y_1) x_{n-1}^{-1} \cr 
&\Delta (x_n, x_j, y_1) \Delta^{-1} (x_n, x_{i-1}, y_1) x_{n-1}^{-1} \Delta
(x_n, x_{i-1}, y_1) x_{i-1} \Delta^{-1} (x_n, x_j, y_1)
x_{n-1} \cr
&\Delta (x_n, x_j, y_1) x_{n-1} x_j^{-1}\cr
=& x_j x_{n-1}^{-1} x_{i-1}^{-1} \Delta^{-1} (x_n, x_{i-1},
y_1) x_{n-1} \Delta (x_n, x_{i-1}, y_1) \Delta^{-1} (x_n,
x_j, y_1) x_{n-1}^{-1}\cr
& \Delta (x_n, x_j, y_1) \Delta^{-1} (x_n, x_{i-1}, y_1) x_{n-1}^{-1} \Delta
(x_n, x_{i-1}, y_1) \Delta^{-1} (x_n, x_j, y_1) x_{n-1}
\cr
&\Delta (x_n, x_j, y_1) x_{i-1} x_{n-1} x_j^{-1} \ ({\rm by\ (PR3)})\cr
=& x_j x_{n-1}^{-1} x_{i-1}^{-1} y_1^{-1} x_n^{-1} 
x_{i-1}^{-1} y_1^{-1} x_{n-1} y_1 x_n x_{i-1} y_1 y_1^{-1}
x_n^{-1} x_j^{-1} y_1^{-1} x_{n-1}^{-1} y_1 x_n x_j \cr
&y_1 y_1^{-1} x_n^{-1} x_{i-1}^{-1} y_1^{-1} x_{n-1}^{-1} y_1 
x_n x_{i-1} y_1 y_1^{-1} x_n^{-1} x_j^{-1} 
y_1^{-1} x_{n-1} y_1 x_n x_j y_1 x_{i-1} x_{n-1} x_j^{-1}\cr
=& x_j x_{n-1}^{-1} x_{i-1}^{-1} y_1^{-1} x_n^{-1} 
x_{i-1}^{-1} x_{n-1} y_1 x_{n-1}^{-1} x_{i-1} x_j^{-1} 
x_{n-1} y_1^{-1} x_{n-1}^{-1} x_j x_{i-1}^{-1} x_{n-1}\cr
&y_1^{-1} x_{n-1}^{-1} 
x_{i-1} x_j^{-1} x_{n-1} y_1 x_{n-1}^{-1} x_j x_n y_1 x_{i-1} x_{n-1} 
x_j^{-1}\cr}

\eject

\halign{
$#$\ &$#$\hfill\cr
=& x_j x_{n-1}^{-1} x_{i-1}^{-1} y_1^{-1} x_n^{-1} x_{n-1}
y_1 x_{i-1} y_1^{-1} y_1 x_j^{-1} y_1^{-1} y_1 x_{i-1}^{-1}
y_1^{-1} y_1 x_j y_1^{-1} x_{n-1}^{-1}  x_n y_1 x_{i-1}\cr
&x_{n-1} x_j^{-1}\cr
=& 1.\cr}

\medskip
$\bullet$\qua Relation (PT3): (PT3) is an easy case if $i=n$. 
So, we assume that $i<n$. Then:

\medskip
\halign{
$#$\ &$#$\hfill\cr
&f (y_1 \alpha_i y_1^{-1}) f (\beta_1^{-1} \alpha_i)^{-1} \cr
=& y_1 x_{n-1}^{-1} \Delta^{-1} (x_n, x_{i-1}, y_1) x_n^{-1}
x_{n-1} \Delta (x_n, x_{i-1}, y_1) x_{n-1} y_1^{-1} 
x_{n-1}^{-1} \cr
&\Delta^{-1} (x_n, x_{i-1}, y_1) x_{n-1}^{-1}x_n \Delta (x_n, x_{i-1}, y_1) x_{n-1} \Delta^{-1}
(x_{n-1}, y_1) x_{n-1}^{-1} x_n \Delta (x_{n-1}, y_1)\cr
=& y_1 x_{n-1}^{-1} y_1^{-1} x_n^{-1} x_{i-1}^{-1} y_1^{-1}
x_n^{-1} x_{n-1} y_1 x_n x_{i-1} y_1 x_{n-1} y_1^{-1}
x_{n-1}^{-1} y_1^{-1} x_n^{-1} x_{i-1}^{-1} y_1^{-1} 
x_{n-1}^{-1} \cr
&x_n y_1 x_n x_{i-1} y_1 x_{n-1} x_{n-1}^{-1} y_1^{-1} x_{n-1}^{-1} x_n y_1
x_{n-1}\cr
=& x_{n-1}^{-1} y_1^{-1} x_{n-1} x_n^{-1} x_{i-1}^{-1}
y_1^{-1} x_n^{-1} x_{n-1} y_1 x_n x_{i-1} y_1 x_{n-1}
x_{n-1}^{-1} y_1^{-1} x_{n-1}^{-1} x_n^{-1} x_{i-1}^{-1}
y_1^{-1} \cr
&x_{n-1}^{-1} x_n y_1 x_n x_{i-1}x_{n-1}^{-1} x_n y_1 x_{n-1}\cr
=& x_{n-1}^{-1} y_1^{-1} x_{n-1} x_n^{-1} x_{i-1}^{-1}
y_1^{-1} x_n^{-1} x_{n-1} y_1 y_1^{-1} x_{n-1}^{-1} 
y_1^{-1} y_1 x_n y_1 x_{i-1} x_n x_{n-1}^{-1} y_1 
x_{n-1}\cr
=&1.\cr} 

\medskip
$\bullet$\qua Relation (PT4): (PT4) is an easy case if either 
$i=n$ or $j\ge 3$. So, we assume that $j=2$ and $i\le n-1$. 
Then: 

\medskip
\halign{
$#$\ &$#$\hfill\cr
&y_2 f(\alpha_i) y_2^{-1}\cr
=& y_2 x_{n-1}^{-1} \Delta^{-1} (x_n, x_{i-1} ,y_1) x_n^{-1} 
x_{n-1} \Delta (x_n, x_{i-1}, y_1) x_{n-1} y_2^{-1} \cr
=& x_{n-1}^{-1} x_{i-1}^{-1} y_2 \Delta^{-1} (x_n, x_{i-1} 
,y_1) x_{n-1} \Delta (x_n, x_{i-1}, y_1) y_2^{-1} x_{n-1}  
\cr
=& x_{n-1}^{-1} x_{i-1}^{-1} \Delta^{-1} (x_n, x_{i-1} ,y_1) 
x_{n-1} \Delta (x_n, x_{i-1}, y_1) x_{n-1} \quad ({\rm by\ 
(PR4)})  \cr
=& f(\alpha_i).\cr}

\medskip
$\bullet$\qua Relation (PT7): (PT7) is an easy case if $j=n-1$. 
So, we assume that $j\le n-2$. We prove by induction on 
$i\ge 2$ that $x_j$ and $f(\beta_i)$ commute. Assume first 
that $i=2$. (PR4) and Lemma 3.4 imply:
$$
x_j \Delta^{-1} (x_{n-1}, y_1, y_2) x_n \Delta (x_{n-1},
y_1, y_2) = \Delta^{-1} (x_{n-1}, y_1, y_2) x_n \Delta
(x_{n-1}, y_1, y_2) x_j,
$$
and this last equality implies:
$$
x_j f(\beta_2) x_j^{-1} = f(\beta_2).
$$
Now, we assume that $i>2$. The first equality of Lemma 3.5 
implies:
$$
f(\beta_i) = f(\beta_{i-1}) y_i f(\beta_{i-1})^{-1} 
y_i^{-1}. 
$$
Thus, since $x_j$ commutes 
with $y_i$ and
with $f(\beta_{i-1})$ (inductive 
hypothesis), $x_j$ also commutes 
with $f(\beta_i)$.

$\bullet$\qua Relation (PT9): The equality
$$
y_{i-1} f(\beta_i) y_{i-1}^{-1} = f(\beta_i) f(\beta_{i-1})
$$
is a straightforward consequence of the second equality of 
Lemma 3.5.

$\bullet$\qua Relation (PT10): The equality
$$
y_{i+1} f(\beta_i) y_{i+1}^{-1} = f(\beta_{i+1})^{-1} 
f(\beta_i)
$$
is a straightforward consequence of the first equality of 
Lemma 3.5.

$\bullet$\qua Relation (PT11): Assume first that $n=1$. Then:

\medskip
\centerline{\vbox{\halign{
$#$\ &$#$\hfill\cr
&f(\alpha_1)^{-1} f(\beta_1)^{-1} f(\beta_2)^{-1}
f(\beta_3)^{-1} z f(\beta_3) z^{-1} f(\beta_1)\cr
=& \Delta^{-2} (x_1, y_1, y_2, y_3, z) \Delta (x_0, x_1,
y_1, y_2, y_3, z) \quad ({\rm by\ Lemma\ 3.6.(i)})\cr
=&1 \quad ({\rm by\ (PR5)}).\cr}}}

\medskip\noindent
Now, assume that $n\ge 2$. Lemma 3.6.(i) implies:

\medskip
\centerline{\vbox{\halign{
$#$\ &$#$\hfill\cr
&x_n^{-1} x_{n-1} f(\beta_1)^{-1} f(\beta_2)^{-1}
f(\beta_3)^{-1} z f(\beta_3) z^{-1} f(\beta_1)\cr 
=& \Delta^{-2} (x_n, y_1, y_2, y_3, z) \Delta (x_{n-1},
x_n, y_1, y_2, y_3, z),\cr}}}

\medskip\noindent
and Lemma 3.6.(ii) implies:
$$
x_n^{-1} x_{n-1} f(\alpha_1) = \Delta^{-2} (x_0, x_n, y_1)
\Delta (x_0, x_{n-1}, x_n, y_1).
$$
Thus:

\medskip
\centerline{\vbox{\halign{
$#$\ &$#$\hfill\cr
&f(\alpha_1)^{-1} f(\beta_1)^{-1} f(\beta_2)^{-1}
f(\beta_3)^{-1} z f(\beta_3) z^{-1} f(\beta_1)\cr
=&\Delta^{-1} (x_0, x_{n-1}, x_n, y_1) \Delta^2 (x_0, x_n,
y_1) \Delta^{-2} (x_n, y_1, y_2, y_3, z) \Delta (x_{n-1},
x_n, y_1, y_2, y_3, z)\cr
=& 1 \quad ({\rm by\ (PR6)}).\cr}}}

\medskip
$\bullet$\qua Relation (PT12): (PT12) is an easy case if 
$i=1,2$. We prove by induction on $i\ge 4$ that $z$ and 
$f(\beta_i)$ commute. Recall first that the first equality 
of Lemma 3.5 implies:
$$
f(\beta_i) = f(\beta_{i-1}) y_i f(\beta_{i-1})^{-1} 
y_i^{-1}.
$$
Assume that $i=4$. Then:

\medskip
\centerline{\vbox{\halign{
$#$\ &$#$\hfil\cr
& z f(\beta_4) z^{-1}\cr
=& z f(\beta_3) y_4 f(\beta_3)^{-1} y_4^{-1} z^{-1}\cr
=& f(\beta_3) f(\beta_2) f(\beta_1) f(\alpha_1)
f(\beta_1)^{-1} y_4 f(\beta_1) f(\alpha_1)^{-1}
f(\beta_1)^{-1} f(\beta_2)^{-1} f(\beta_3)^{-1} y_4^{-1}\cr
&\hfill{\rm by\ (PT11)}\cr
=& f(\beta_3) y_4 f(\beta_3)^{-1} y_4^{-1} \quad ({\rm by\
(PT4)\ and\ (PT8)})\cr
=& f(\beta_4).\cr}}}

\medskip\noindent
Now, we assume that $i>4$. Then $z$ commutes with 
$f(\beta_i)$, since it commutes 
with $y_i$ and with $f(\beta_{i-1})$ 
(inductive hypothesis). \qed

\medskip
Now, in view of Proposition 3.3, and applying Lemma 2.5 to 
the exact sequences (2.3) of Subsection 2.2, one has 
immediately the following presentation for $\PM (F_{g,r+1}, 
\PP_n)$. 

\medskip\noindent
{\bf Proposition 3.7}\qua {\sl Let $g,r\ge 1$, let $n\ge 0$, 
and let $P\Gamma_{g,r,n}$ be the Coxeter graph drawn in 
Figure \ref{f19}. Then $\PM (F_{g,r+1}, \PP_n)$ is isomorphic with 
the quotient of $A (P\Gamma_{g,r,n})$ by the following 
relations.

\smallskip
$\bullet$\qua Relations from $\MM (F_{g,1})$:

\medskip
\centerline{\vbox{\halign{
$#$  \hfill &\hfill $#$ &\ $#$\ &$#$  \hfill &$#$ 
\hfill \cr
{\rm(PR1)} \ &\Delta^4 (y_1, y_2, y_3, z) &=& \Delta^2 (x_0, y_1, 
y_2, y_3, z) &{\hbox{\sl if}}\ g\ge 2,\cr
\noalign{\smallskip}
{\rm(PR2)}\ &\Delta^2 (y_1, y_2, y_3, y_4, y_5, z) &=& \Delta 
(x_0, y_1, y_2, y_3, y_4, y_5, z) &{\hbox{\sl if}}\ g\ge 3. 
\cr}}}

\medskip
$\bullet$\qua Relations of commutation:

\medskip
\centerline{\vbox{\halign{
$#$  \hfill &$#$\ &$#$  \hfill &$#$ \hfill \cr
{\rm(PR3)} &&x_k \Delta^{-1} (x_{i+1}, x_j, y_1) x_i \Delta 
(x_{i+1}, x_j, y_1)\cr
&=&\Delta^{-1} (x_{i+1}, x_j, y_1) x_i \Delta (x_{i+1}, x_j, 
y_1) x_k &{\hbox{\sl if}}\ 0\le k< j< i\le r+n-1, \cr
\noalign{\smallskip}
{\rm(PR4)} &&y_2 \Delta^{-1} (x_{i+1}, x_j, y_1) x_i \Delta 
(x_{i+1}, x_j, y_1)\cr
&=&\Delta^{-1} (x_{i+1}, x_j, y_1) x_i \Delta (x_{i+1}, x_j, 
y_1) y_2 &{\hbox{\sl if}}\ 0\le j< i\le r+n-1, \cr}}}

\medskip
$\bullet$\qua Relations between fundamental elements:}

\medskip
\centerline{\vbox{\halign{
$#$  \hfill &\hfill $#$ &\ $#$\ &$#$  \hfill &$#$ 
\hfill \cr
{\rm(PR5a)} &u_1 &=& \Delta (x_0, x_1, y_1, y_2, y_3,z) 
\Delta^{-2} (x_1, y_1, y_2, y_3, z), \cr
\noalign{\smallskip}
{\rm(PR6a)} &u_{i+1} &=&\Delta (x_i, x_{i+1}, y_1, y_2, y_3, z) 
\Delta^{-2} (x_{i+1}, y_1, y_2, y_3, z) \cr 
&&& \Delta^2 (x_0, x_{i+1}, y_1) \Delta^{-1} (x_0, 
x_i, x_{i+1}, y_1) \quad {\hbox{\sl if}} \quad 1\le i\le r-1, &\cr
\noalign{\smallskip}
{\rm(PR6b)} &\multispan 3 $\Delta (x_i, x_{i+1}, y_1, y_2, 
y_3, z) \Delta^{-
2} (x_{i+1}, y_1, y_2, y_3, z)$\hfill \cr 
&&=& \Delta (x_0, x_i, x_{i+1}, y_1) \Delta^{-2} (x_0, 
x_{i+1}, y_1) \quad {\hbox{\sl if}} \quad r\le i\le n+r-1.&\cr  
\cr}}}

\begin{figure}[ht]
\centerline{\input{pmc_19}}
\caption{\label{f19} Coxeter graph associated with $\PM 
(F_{g,r+1}, \PP_n)$}
\end{figure}

\medskip
Let $PG(g,r,n)$ denote the quotient of $A( 
P\Gamma_{g,r,n})$ by the relations (PR1),(PR2),
(PR3),(PR4),(PR5a), (PR6a), (PR6b). Consider the Dehn twists $a_0, 
\dots, a_{n+r}$, $b_1, \dots, b_{2g-1}$, $c$, $d_1, \dots, 
d_r$ represented in Figure \ref{f10}. Then an isomorphism $\bar 
\rho : PG (g,r,n) \to \PM (F_{g,r+1}, \PP_n)$ between $PG 
(g,r,n)$ and $\PM (F_{g,r+1}, \PP_n)$ is given by $\bar 
\rho (x_i) =a_i$ for $i=0, \dots, n+r$, $\bar \rho (y_i) 
=b_i$ for $i=1, \dots, 2g-1$, $\bar \rho (z) =c$, and $\bar 
\rho (u_i) =d_i$ for $i=1, \dots, r$.

As in Lemma 3.6, we use the algorithm of \cite{DePa} 
to prove the following.

\medskip\noindent
{\bf Lemma 3.8}\qua {\sl {\rm(i)}\qua We number the vertices of 
the Coxeter graph $D_6$ 
according to Figure \ref{f6}. Then the following equality holds in 
$A(D_6)$.

\medskip
\halign{
$#$\ &$#$\hfill\cr
&\Delta^2 (x_1, x_3, \dots, x_6) \Delta^{-1} (x_1, x_2, 
x_3, \dots, x_6) 
=x_6 x_5 x_4 x_3 x_1 x_2^{-1} x_3^{-1} x_4^{-1} x_5^{-1} 
x_6^{-1} x_5 x_4\cr
& x_3 x_2x_1^{-1} x_3^{-1} x_4^{-1} x_5^{-
1} x_4 x_3 x_1 x_2^{-1} x_3^{-1} x_4^{-1} x_2 x_3 x_2 
x_1^{-1} x_3^{-1} x_2^{-1}.\cr}

\medskip
{\rm(ii)}\qua We number the vertices of the Coxeter graph $D_4$ 
according to Figure \ref{f6}. Then the following equality holds in 
$A(D_4)$.}
$$
\Delta (x_1, x_2, x_3, x_4) \Delta^{-2} (x_1, x_3, x_4) = 
x_2 x_3 x_2^{-1} x_1 x_3^{-1} x_2^{-1} x_4 x_3 x_2 x_1^{-1} 
x_3^{-1} x_4^{-1}. \eqno{\qed}
$$

\noindent
{\bf Proof of Theorem 3.1}\qua 
Recall that $\Gamma_{g,r,n}$ denotes the Coxeter graph 
drawn in Figure \ref{f16}, and that $G(g,r,n)$ denotes the 
quotient of $A(\Gamma_{g,r,n})$ by the relations 
(R1),\dots,(R7), (R8a), (R8b).
Recall also that there is a well defined 
epimorphism $\bar\rho : G(g,r,n) \to \MM (F_{g,r+1}, 
\PP_n)$ which sends $x_i$ on $a_i$ for $i=0, \dots, r+1$, 
$y_i$ on $b_i$ for $i=1, \dots, 2g-1$, $z$ on $c$, $u_i$ on 
$d_i$ for $i=1, \dots, r$, and $v_i$ on $\tau_i$ for $i=1, 
\dots, n-1$. Our aim now is to construct a homomorphism 
$\bar f: \MM (F_{g,r+1}, \PP_n) \to G(g,r,n)$ such that 
$\bar f \circ \bar \rho$ is the identity of $G(g,r,n)$. The 
existence of such a homomorphism clearly proves that $\bar 
\rho$ is an isomorphism.

We set $A_0=x_r$, $A_1= x_{r+1}$, and
$$
A_i = x_r^{1-i} \Delta (x_{r+1}, v_1, \dots, v_{i-1}) \quad 
{\rm for}\ i=2, \dots, n.
$$
These expressions are viewed as elements of $G(g,r,n)$. 
Note that, by Proposition 2.12, we have $\bar \rho (A_i) = 
a_{r+i}$ for all $i=0, 1, \dots, n$.

\medskip\noindent
{\bf Assertion 1}\qua {\sl {\rm(i)}\qua The following relations hold in 
$G(g,r,n)$:

\medskip
\centerline{\vbox{\halign{
$#$ \quad \hfill& \hfill $#$ &\ $#$\ &$#$ \quad \hfill &$#$ 
\hfill \cr
{\rm(T1)} &A_{i-1} A_{i+1} &=& v_i A_i v_i A_i \cr
&&=& A_i v_i A_i v_i &{\hbox{\sl for}}\ 1\le i\le n-1,\cr
{\rm(T2)} &A_i A_j &=& A_j A_i &{\hbox{\sl for}}\ 0\le i< j\le n,\cr
{\rm(T3)} &A_i v_j &=& v_j A_i &{\hbox{\sl for}}\ i \neq j,\cr
{\rm(T4)} &y_1 A_i y_1 &=& A_i y_1 A_i &{\hbox{\sl for}}\ 0\le i\le n.
\cr}}}

\medskip
{\rm(ii)}\qua The relations (T1),\dots,(T4) imply that there is a 
well defined homomorphism $h_i: A(B_4) \to G(g,r,n)$ which 
sends $x_1$ on $v_i$, $x_2$ on $A_i$, $x_3$ on $y_1$, and 
$x_4$ on $A_{i-1}$. Then the following relation holds in 
$G(g,r,n)$:}
$$
{\rm(T5)} \quad h_i (\Delta (x_1, x_2, x_3, x_4)) = h_i 
(\Delta^2 (x_1, x_2, x_3)) \quad {\hbox{\sl for}}\ 1\le i
\le n.
$$

\noindent
{\bf Proof of Assertion 1}\qua $\bullet$\qua Relation (T1): 

\medskip
\halign{
\hfill $#$ &\ $#$\ &$#$ \hfill \cr
A_{i+1} &=& x_r^{-i} \Delta (x_{r+1}, v_1, \dots, v_i)\cr
&=& x_r^{-i} v_i v_{i-1} \dots v_1 x_{r+1} v_1 \dots v_{i-
1} v_i \Delta (x_{r+1}, v_1, \dots, v_{i-1})
\quad ({\rm by\  2.9})\cr
&=& x_r^{-i} v_i \Delta (x_{r+1}, v_1, \dots, v_{i-1}) 
\Delta^{-1} (x_{r+1}, v_1, \dots, v_{i-2}) v_i \cr
&&\Delta (x_{r+1}, v_1, \dots, v_{i-1})\cr
&=& x_r^{i-2} \Delta^{-1} (x_{r+1}, v_1, \dots, v_{i-2}) 
v_i x_r^{1-i} \Delta (x_{r+1}, v_1, \dots, v_{i-1}) v_i 
x_r^{1-i} \cr
&&\Delta (x_{r+1}, v_1, \dots, v_{i-1})\cr
&=& A_{i-1}^{-1} v_i A_i v_i A_i. \cr}

\medskip\noindent
Similarly:
$$
A_{i+1} = A_{i-1}^{-1} A_i v_i A_i v_i.
$$

$\bullet$\qua The relations (T2) and (T3) are direct 
consequences of the ``braid'' 
relations in $A(\Gamma_{g,r,n})$.

$\bullet$\qua Now, we prove (T4) and (T5) by induction on $i$. 
First, assume $i=1$. Then (T4) follows from the ``braid'' 
relation $y_1 x_{r+1} y_1 = x_{r+1} y_1 x_{r+1}$ in 
$A(\Gamma_{g,r,n})$, and (T5) follows from the relation 
(R7) in the definition of $G(g,r,n)$. 

Now, assume $i>1$. Then the relation (T4) follows from the 
following sequence of equalities.

\medskip
\halign{
$#$\ &$#$ \hfill\cr
& A_i y_1 A_i y_1^{-1} A_i^{-1} y_1^{-1}\cr
=& A_{i-2}^{-1} v_{i-1} A_{i-1} v_{i-1} A_{i-1} y_1 A_{i-1}
v_{i-1} A_{i-1} v_{i-1} A_{i-2}^{-1} y_1^{-1} A_{i-2} v_{i-
1}^{-1} A_{i-1}^{-1} v_{i-1}^{-1} A_{i-1}^{-1} y_1^{-1}\cr 
&\quad ({\rm by\ (T1)}) \cr
=& A_{i-2}^{-1} \cdot v_{i-1} A_{i-1} v_{i-1} A_{i-1} y_1 
A_{i-1} v_{i-1} A_{i-1} y_1 A_{i-2}^{-1} y_1^{-1} A_{i-
1}^{-1} v_{i-1}^{-1} A_{i-1}^{-1} y_1^{-1} A_{i-2}^{-1}\cr 
&\cdot A_{i-2} \quad ({\rm by\ (T2),(T3),\ induction}) \cr
=& A_{i-2}^{-1} \cdot h_{i-1} (\Delta^2 (x_1, x_2, x_3) 
\Delta^{-1} (x_1, x_2, x_3, x_4)) \cdot A_{i-2} \quad ({\rm 
by\ Proposition\ 2.9}) \cr
=& 1 \quad ({\rm by\ induction}). \cr} 

\medskip\noindent
The Relation (T5) follows from the following sequence of 
equalities.

\medskip
\halign{
$#$\ &$#$ \hfill\cr
& h_i (\Delta^{-1} (x_1, x_2, x_3, x_4) \Delta^2 (x_1, x_2, 
x_3)) \cr
=& A_{i-1}^{-1} y_1^{-1} A_i^{-1} v_i^{-1} A_i^{-1} y_1^{-
1} A_{i-1}^{-1} y_1 A_i v_i y_1 A_i v_i y_1 A_i v_i \quad 
({\rm by\ Propositions\ 2.8\ , \ 2.9}) \cr
=& A_{i-1}^{-1} y_1^{-1} A_{i-2} v_{i-1}^{-1} A_{i-1}^{-1} 
v_{i-1}^{-1} A_{i-1}^{-1} v_i^{-1} v_{i-1}^{-1} A_{i-1}^{-
1} v_{i-1}^{-1} A_{i-1}^{-1} A_{i-2} y_1^{-1} A_{i-1}^{-1} 
y_1 A_{i-2}^{-1} A_{i-1} \cr
&v_{i-1} A_{i-1} v_{i-1}v_i y_1 A_{i-2}^{-1} A_{i-1} v_{i-1} A_{i-1} v_{i-1} 
v_i y_1 v_{i-1} A_{i-1} v_{i-1} A_{i-1} A_{i-2}^{-1}v_i 
\quad (\ {\rm T1}\ ) \cr
=&A_{i-2} \cdot A_{i-1}^{-1} A_{i-2}^{-1} y_1^{-1} A_{i-2} 
v_{i-1}^{-1} A_{i-1}^{-1} v_{i-1}^{-1} A_{i-1}^{-1} v_i^{-
1} v_{i-1}^{-1} A_{i-1}^{-1} v_{i-1}^{-1} A_{i-1}^{-1} 
A_{i-2} A_{i-1} y_1^{-1}\cr
&A_{i-1}^{-1} A_{i-2}^{-1} A_{i-1}v_{i-1} A_{i-1} v_{i-1} v_i y_1 A_{i-2}^{-1} A_{i-1} 
v_{i-1} A_{i-1} v_{i-1} v_i y_1 v_{i-1} A_{i-1} v_{i-1} 
A_{i-1} v_i \cr
&\cdot A_{i-2}^{-1} ({\rm by\ (T2),(T3),\ induction}) \cr
=&A_{i-2} A_{i-1}^{-1} v_{i-1}^{-1} \cdot y_1 A_{i-2}^{-1} 
y_1^{-1} A_{i-1}^{-1} v_{i-1}^{-1} v_i^{-1} A_{i-1}^{-1} 
v_{i-1}^{-1} A_{i-1}^{-1} y_1^{-1} A_{i-2}^{-
1} y_1 A_{i-1} v_{i-1} v_i y_1\cr
&A_{i-2}^{-1} A_{i-1}v_{i-1} A_{i-1} v_{i-1} v_i y_1 v_{i-1} A_{i-1} 
v_{i-1} v_i v_{i-1}^{-1} \cdot v_{i-1} A_{i-1} A_{i-2}^{-1} \cr
&({\rm by\ (T2),(T3),\ induction}) \cr
=&A_{i-2} A_{i-1}^{-1} v_{i-1}^{-1}  y_1 \cdot A_{i-2}^{-1} 
y_1^{-1} A_{i-1}^{-1} v_{i-1}^{-1} v_i^{-1} y_1  A_{i-2} 
\cdot h_{i-1} (\Delta^{-1} (x_1, x_2, x_3 ,x_4) \cr
&\Delta(x_1, x_2, x_3)) \cdot y_1A_{i-1} v_{i-1} v_i y_1 A_{i-2}^{-1} A_{i-1} v_{i-1} 
A_{i-1} v_{i-1} v_i v_{i-1} y_1 A_{i-1} v_i^{-1} v_{i-1} 
v_i y_1\cr
& \cdot  y_1^{-1} v_{i-1} A_{i-1} A_{i-2}^{-1}({\rm by\ Proposition\ 2.9}) \cr
=&A_{i-2} A_{i-1}^{-1} v_{i-1}^{-1} y_1 \cdot A_{i-2}^{-1} 
y_1^{-1} A_{i-1}^{-1} v_{i-1}^{-1} v_i^{-1} y_1  A_{i-2} 
v_{i-1}^{-1} A_{i-1}^{-1} y_1^{-1} v_{i-1}^{-1} A_{i-1}^{-
1} y_1^{-1} v_{i-1}^{-1}\cr
& A_{i-1}^{-1} y_1^{-1} y_1 A_{i-
1}v_{i-1} v_i y_1 A_{i-2}^{-1} A_{i-1} v_{i-1} A_{i-1} 
v_{i-1} v_i v_{i-1} y_1 A_{i-1} v_i^{-1} v_{i-1} v_i y_1 \cr
&y_1^{-1} v_{i-1} A_{i-1} A_{i-2}^{-1} \quad ({\rm 
by\ induction}) \cr
=&A_{i-2} A_{i-1}^{-1} v_{i-1}^{-1} y_1 \cdot A_{i-2}^{-1} 
y_1^{-1} A_{i-1}^{-1} y_1 v_{i-1}^{-1} v_i^{-1} v_{i-1}^{-
1} A_{i-1}^{-1} A_{i-2} y_1^{-1} A_{i-2}^{-1} v_{i-1}^{-1} 
v_i v_{i-1} \cr
&A_{i-1} v_i v_{i-1} v_i y_1 A_{i-1} v_i^{-1} y_1 v_{i-1} v_i \cdot y_1^{-1} 
v_{i-1} A_{i-1} A_{i-2}^{-1} ({\rm \ (T2),(T3),\ induction}) \cr
=&A_{i-2} A_{i-1}^{-1} v_{i-1}^{-1} y_1 \cdot A_{i-2}^{-1} 
A_{i-1} y_1^{-1} A_{i-1}^{-1} v_{i-1}^{-1} v_i^{-1} v_{i-
1}^{-1} A_{i-1}^{-1} y_1^{-1} A_{i-2}^{-1} y_1 v_i v_{i-1}   
v_i^{-1} A_{i-1} \cr
&v_i v_{i-1} y_1 A_{i-1}y_1 v_{i-1} v_i \cdot y_1^{-1} v_{i-1} A_{i-1} A_{i-
2}^{-1} \quad ({\rm by\ (T2),(T3),\ induction}) \cr
=&A_{i-2} A_{i-1}^{-1} v_{i-1}^{-1} y_1 A_{i-1} \cdot A_{i-
2}^{-1} y_1^{-1} A_{i-1}^{-1} v_{i-1}^{-1} v_i^{-1} v_{i-
1}^{-1} A_{i-1}^{-1} y_1^{-1} A_{i-2}^{-1} v_i \cdot 
y_1 v_{i-1} A_{i-1} y_1 \cr
&v_{i-1} A_{i-1} y_1 v_{i-1} 
A_{i-1}\cdot v_i\cdot A_{i-1}^{-1} y_1^{-1} v_{i-1} A_{i-1} A_{i-
2}^{-1} \quad ({\rm \ (T2),(T3),\ induction}) \cr
=&A_{i-2} A_{i-1}^{-1} v_{i-1}^{-1} y_1 A_{i-1} \cdot A_{i-
2}^{-1} y_1^{-1} A_{i-1}^{-1} v_{i-1}^{-1} v_i^{-1} v_{i-
1}^{-1} A_{i-1}^{-1} y_1^{-1} A_{i-2}^{-1} v_i \cr
&\cdot h_{i-1} (\Delta (x_1, x_2, x_3)) \cdot v_i\cdot 
A_{i-1}^{-1} y_1^{-1} v_{i-1} A_{i-1} A_{i-2}^{-1} 
\quad ({\rm by\ Proposition\ 2.8}) \cr
=&A_{i-2} A_{i-1}^{-1} v_{i-1}^{-1} y_1 A_{i-1} \cdot A_{i-
2}^{-1} y_1^{-1} A_{i-1}^{-1} v_{i-1}^{-1} v_i^{-1} v_{i-
1}^{-1} A_{i-1}^{-1} y_1^{-1} A_{i-2}^{-1} v_i  A_{i-2} y_1 
A_{i-1} \cr
&v_{i-1} A_{i-1} y_1 A_{i-2} v_i \cdot A_{i-1}^{-1} y_1^{-1} 
v_{i-1} A_{i-1} A_{i-2}^{-1}
 \quad ({\rm by\ induction}) \cr
=&A_{i-2} A_{i-1}^{-1} v_{i-1}^{-1} y_1 A_{i-1} \cdot A_{i-
2}^{-1} y_1^{-1} A_{i-1}^{-1} v_{i-1}^{-1} v_i^{-1} v_i 
v_{i-1} v_i^{-1} A_{i-1} y_1 A_{i-2} v_i \cdot A_{i-1}^{-1} 
y_1^{-1} \cr
&v_{i-1} A_{i-1} A_{i-2}^{-1} 
({\rm by\ (T2),(T3),\ induction}) \cr
=&1 \quad ({\rm by\ (T2),(T3),\ induction}) \cr}

\medskip\noindent
{\bf Assertion 2}\qua {\sl 
Recall that $P\Gamma_{g,r,n}$ denotes the Coxeter graph
drawn in Figure \ref{f19}.
There is a well defined 
homomorphism $g: A (P\Gamma_{g,r,n}) \to G(g,r,n)$ which 
sends $x_i$ on $x_i$ for $i=0, \dots, r+1$, $x_{r+i}$ on 
$A_i$ for $i=2, \dots, n$, $y_i$ on $y_i$ for $i=1, \dots, 
2g-1$, $z$ on $z$, and $u_i$ on $u_i$ for $i=1, \dots, r$.}

\medskip\noindent
{\bf Proof of Assertion 2}\qua We have to verify that the 
following relations hold in $G(g,r,n)$.

\medskip
\centerline{\vbox{\halign{
$#$\quad\hfill&\hfill$#$& $#$ &$#$\quad\hfill&$#$\hfill\cr
{\rm(T6)} &A_i A_j &=& A_j A_i &{\rm for}\ 1\le i\le j\le n,\cr
{\rm(T7)} &x_i A_j &=& A_j x_i &{\rm for}\ 0\le i\le r\ {\rm 
and}\ 1\le j\le n, \cr
{\rm(T8)} &y_1 A_i y_1 &=& A_i y_1 A_i &{\rm for}\ 1\le i\le 
n,\cr
{\rm(T9)} &A_i y_j &=& y_j A_i &{\rm for}\ 1\le i\le n\ {\rm 
and}\ 2\le j\le 2g-1, \cr
{\rm(T10)} &A_i z &=& z A_i &{\rm for}\ 1\le i\le n, \cr
{\rm(T11)} &A_i u_j &=& u_j A_i &{\rm for}\ 1\le i\le n\ {\rm 
and}\ 1\le j\le r. \cr}}}

\medskip\noindent
The relations (T6) and (T8) hold by Assertion 1, and the 
other relations are direct consequences of the ``braid'' 
relations in $A (\Gamma_{g,r,n})$.

\medskip
Recall that $PG (g,r,n)$ denotes the quotient of $A( 
P\Gamma_{g,r,n})$ 
by the relations (PR1),\dots,(PR4), (PR5a), (PR6a), (PR6b),
and that this quotient is isomorphic with $\PM (F_{g,r+1}, 
\PP_n)$ (see Proposition 3.7).

\medskip\noindent
{\bf Assertion 3}\qua {\sl 
The homomorphism $g: A 
(P\Gamma_{g,r,n}) \to G(g,r,n)$ induces a homomorphism 
$\bar g: PG (g,r,n) \to G(g,r,n)$.}

\medskip\noindent
{\bf Proof of Assertion 3}\qua It suffices to show that the 
following relations hold in $G(g,r,n)$.

\medskip
\centerline{\vbox{\halign{
$#$\hfill& $#$ &$#$\hfill&$#$\hfill\cr
{\rm(T12)} &&g (x_k \Delta^{-1} (x_{i+1}, x_j, y_1) x_i \Delta 
(x_{i+1}, x_j, y_1)) \cr
&=&g (\Delta^{-1} (x_{i+1}, x_j, y_1) x_i \Delta (x_{i+1}, 
x_j, y_1) x_k) \ {\rm for}\ 0\le k< j< i\le r+n-1, &\cr
\noalign{\smallskip}
{\rm(T13)} &&g (y_2 \Delta^{-1} (x_{i+1}, x_j, y_1) x_i \Delta 
(x_{i+1}, x_j, y_1)) \cr
&=&g (\Delta^{-1} (x_{i+1}, x_j, y_1) x_i \Delta (x_{i+1}, 
x_j, y_1) y_2) \ {\rm for}\ 0\le j< i\le r+n-1,& \cr
\noalign{\smallskip}
{\rm(T14)} &&g (\Delta (x_i, x_{i+1}, y_1, y_2, y_3, z) 
\Delta^{-2} (x_{i+1}, y_1, y_2, y_3, z))\cr
&=& g (\Delta (x_0, x_i, x_{i+1}, y_1) \Delta^{-2} (x_0, 
x_{i+1}, y_1)) \ {\rm for}\ \ r+1\le i\le r+n-1.& \cr}}}

\medskip
$\bullet$\qua Relation (T12): for $i\ge r+1$ and $j< i-1$, we 
have:

\medskip
\halign{
$#$\hfill& $#$ &$#$\hfill\cr
{\rm(E1)} &&g (\Delta^{-1} (x_{i+1}, x_j, y_1) x_i \Delta 
(x_{i+1}, x_j, y_1)) \cr
&=& y_1^{-1} g(x_j)^{-1} A_{i-r+1}^{-1} y_1^{-1} A_{i-r} 
y_1 A_{i-r+1} g(x_j) y_1\cr
&=& y_1^{-1} g(x_j)^{-1} A_{i-r-1} v_{i-r}^{-1} A_{i-r}^{-
1} v_{i-r}^{-1} A_{i-r}^{-1} y_1^{-1} A_{i-r} y_1 A_{i-r} 
v_{i-r} A_{i-r} v_{i-r} \cr
&& A_{i-r-1}^{-1} g(x_j) y_1 \quad ({\rm by\ (T1)}) \cr
&=&v_{i-r}^{-1} y_1^{-1} g(x_j)^{-1} A_{i-r}^{-1} A_{i-r-1}   
v_{i-r}^{-1} A_{i-r}^{-1} A_{i-r} y_1 A_{i-r}^{-1} A_{i-r} 
v_{i-r}A_{i-r-1}^{-1} A_{i-r}\cr
&& g(x_j) y_1 v_{i-r} \quad ({\rm by\ (T2),(T3),(T4)}) \cr
&=&v_{i-r}^{-1} y_1^{-1} g(x_j)^{-1} A_{i-r}^{-1} y_1^{-1}  
A_{i-r-1} y_1 A_{i-r} g(x_j) y_1 v_{i-r} \quad ({\rm by\ 
(T2),(T3),(T4)}) \cr
&=& v_{i-r}^{-1} g (\Delta^{-1} (x_i, x_j, y_1) x_{i-1} 
\Delta (x_i, x_j, y_1)) v_{i-r}. \cr}

\medskip\noindent
For $i\ge r+1$ and $j=i-1$ we have:

\medskip
\halign{
$#$\quad\hfill& $#$ &$#$\hfill\cr
{\rm(E2)} &&g (\Delta^{-1} (x_{i+1}, x_{i-1}, y_1) x_i \Delta 
(x_{i+1}, x_{i-1}, y_1))\cr
&=& y_1^{-1} A_{i-r-1}^{-1} A_{i-r+1}^{-1} y_1^{-1} A_{i-r} 
y_1 A_{i-r+1} A_{i-r-1} y_1\cr
&=& y_1^{-1} A_{i-r-1}^{-1} A_{i-r-1} v_{i-r}^{-1} A_{i-
r}^{-1} v_{i-r}^{-1} A_{i-r}^{-1} y_1^{-1} A_{i-r} y_1 
A_{i-r} v_{i-r} A_{i-r} v_{i-r}\cr
&& A_{i-r-1}^{-1} A_{i-r-1}y_1 \quad ({\rm by\ (T1)}) \cr
&=& v_{i-r}^{-1} y_1^{-1} A_{i-r}^{-1} v_{i-r}^{-1} A_{i-
r}^{-1} A_{i-r} y_1 A_{i-r}^{-1} A_{i-r} v_{i-r} A_{i-r}  
y_1 v_{i-r} \cr
&&({\rm by\ (T2),(T3),(T4)}) \cr
&=& v_{i-r}^{-1} y_1^{-1} y_1 A_{i-r} y_1^{-1} y_1 v_{i-r} 
\quad ({\rm by\ (T2),(T3),(T4)}) \cr
&=& v_{i-r}^{-1} A_{i-r} v_{i-r}. \cr}

\medskip\noindent
First, assume that $i\le r$. Then the relation (T12) follows 
from the relation (R3) in the definition of $G(g,r,n)$. 
Now, we assume that $j<r\le i\le r+n-1$, and we prove by 
induction on $i$ that the relation (T12) holds. The case 
$i=r$ follows from the relation (R3) in the definition of 
$G(g,r,n)$, and the case $i>r$ follows from the inductive 
hypothesis and from the equality (E1) above. Now, we assume 
that $r\le j< i\le r+n-1$, and we prove, again by induction 
on $i$, that the relation (T12) holds. The case $i=j+1$ 
follows from the equality (E2) above, and the case $i>j+1$ 
follows from the inductive hypothesis and from the equality 
(E1).

$\bullet$\qua The relation (T13) can be shown in the same 
manner as the relation (T12).

$\bullet$\qua Relation (T14): We prove by induction on $i\ge 
\sup\{ r,1 \}$ that the relation (T14) holds in $G(g,r,n)$. 
If $i=r \ge 1$, then the relation (T14) follows from the 
relation (R8b) in the definition of $G(g,r,n)$. Assume 
$r=0$ and $i=1$. Then:

\medskip
\halign{
$#$\ &$#$\hfill\cr
&g (\Delta^2 (x_2, y_1, y_2, y_3, z) \Delta^{-1} (x_1, x_2, 
y_1, y_2, y_3, z) \Delta (x_0, x_1, x_2, y_1) \Delta^{-2} 
(x_0, x_2, y_1)) \cr
=& z y_3 y_2 y_1 A_2 A_1^{-1} y_1^{-1} y_2^{-1} y_3^{-1} 
z^{-1} y_3 y_2 y_1 A_1 A_2^{-1} y_1^{-1} y_2^{-1} y_3^{-1} 
y_2 y_1 A_2 A_1^{-1} y_1^{-1} y_2^{-1} A_1 y_1\cr
& A_1 A_2^{-1}y_1^{-1} A_1^{-1} \cdot A_1 y_1 A_1^{-1} A_2 y_1^{-1} A_1^{-1} A_0 y_1 A_1 
A_2^{-1} y_1^{-1} A_0^{-1} \quad ({\rm by\ Lemma\ 3.8}) \cr
=& z y_3 y_2 y_1 v_1 A_1 v_1 A_1 A_0^{-1} A_1^{-1} y_1^{-1} 
y_2^{-1} y_3^{-1} z^{-1} y_3 y_2 y_1 A_1 A_0 A_1^{-1} 
v_1^{-1} A_1^{-1} v_1^{-1} y_1^{-1} y_2^{-1}\cr
& y_3^{-1} y_2 y_1 v_1 A_1 v_1A_1 A_0^{-1} A_1^{-1} 
y_1^{-1} y_2^{-1} A_0 y_1 A_1 
A_0 A_1^{-1} v_1^{-1} A_1^{-1} v_1^{-1} y_1^{-1} A_0^{-1} 
\quad ({\rm  T1}) \cr
=& v_1 \cdot z y_3 y_2 y_1  A_1 A_0^{-1} y_1^{-1} y_2^{-1} 
y_3^{-1} z^{-1} y_3 y_2 y_1 A_0 A_1^{-1} y_1^{-1} y_2^{-1} 
y_3^{-1} y_2 y_1 A_1 A_0^{-1} y_1^{-1} y_2^{-1}\cr
& A_0 y_1 A_0 A_1^{-1}y_1^{-1} A_0^{-1} \cdot v_1^{-1} \quad ({\rm by\ 
(T2),(T3),(T4)}) \cr
=& v_1 \cdot \Delta^2 (x_1, y_1, y_2 ,y_3, z) \Delta^{-1} 
(x_0, x_1, y_1, y_2, y_3, z) \cdot v_1^{-1} \quad ({\rm by\ 
Lemma\ 3.8}) \cr
=& 1 \quad ({\rm by\ (R8a)}). \cr}

\medskip\noindent
Now, we assume that $i> \sup\{r,1\}$. Then:

\medskip
\halign{
$#$\ &$#$\hfill\cr
&g( \Delta^2 (x_{i+1}, y_1, y_2, y_3, z) \Delta^{-1} (x_i, 
x_{i+1}, y_1, y_2, y_3, z) \Delta (x_0, x_i, x_{i+1}, y_1) 
\cr
&\Delta^{-2} (x_0, x_{i+1}, y_1)) \cr
= &z y_3 y_2 y_1 A_{i-r+1} A_{i-r}^{-1} y_1^{-1} y_2^{-1} 
y_3^{-1} z^{-1} y_3 y_2 y_1 A_{i-r} A_{i-r+1}^{-1} y_1^{-1} 
y_2^{-1} y_3^{-1} y_2 y_1 A_{i-r+1}\cr 
&A_{i-r}^{-1}y_1^{-1} 
y_2^{-1} A_{i-r}y_1 A_{i-r} A_{i-r+1}^{-1} y_1^{-1} A_{i-r}^{-1} 
\cdot A_{i-r} y_1 A_{i-r}^{-1} A_{i-r+1} y_1^{-1} A_{i-
r}^{-1} x_0 y_1\cr
& A_{i-r} A_{i-r+1}^{-1} y_1^{-1} x_0^{-1} 
\quad ({\rm by\ Lemma\ 3.8}) \cr
=& z y_3 y_2 y_1 v_{i-r} A_{i-r} v_{i-r} A_{i-r} A_{i-r-
1}^{-1} A_{i-r}^{-1} y_1^{-1} y_2^{-1} y_3^{-1} z^{-1} y_3 
y_2 y_1 A_{i-r} A_{i-r-1} A_{i-r}^{-1}\cr
& v_{i-r}^{-1} A_{i-
r}^{-1} v_{i-r}^{-1} y_1^{-1}y_2^{-1} y_3^{-1} y_2 y_1 v_{i-r} A_{i-r} v_{i-r} 
A_{i-r} A_{i-r-1}^{-1} A_{i-r}^{-1} y_1^{-1} y_2^{-1} x_0 
y_1 A_{i-r} \cr
&A_{i-r-1} A_{i-r}^{-1} v_{i-r}^{-1} A_{i-r}^{-
1} v_{i-r}^{-1} y_1^{-1} x_0^{-1} 
({\rm by\ (T1)}) \cr
=& v_{i-r} \cdot z y_3 y_2 y_1 A_{i-r} A_{i-r-1}^{-1} 
y_1^{-1} y_2^{-1} y_3^{-1} z^{-1} y_3 y_2 y_1 A_{i-r-1} A_{i-r}^{-1}
y_1^{-1} y_2^{-1} y_3^{-1} y_2 y_1 \cr
&A_{i-r} A_{i-r-1}^{-1} 
y_1^{-1} y_2^{-1}x_0 y_1 A_{i-r-1} A_{i-
r}^{-1} y_1^{-1} x_0^{-1} \cdot v_{i-r}^{-1} \quad ({\rm 
by\ (T2),(T3),(T4)}) \cr
=& v_{i-r} \cdot g (\Delta^2 (x_i, y_1, y_2, y_3, z) 
\Delta^{-1} (x_{i-1}, x_i, y_1, y_2, y_3, z) \Delta (x_0, 
x_{i-1}, x_i, y_1)\cr
& \Delta^{-2} (x_0, x_i, y_1))\cdot v_{i-r}^{-1} \quad ({\rm by\ Lemma\ 3.8}) \cr
=& 1 \quad ({\rm by\ induction}). \cr} 

\medskip
Let $V_1, \dots, V_{n-1}$ denote the natural generators of 
the Artin group $A(A_{n-1})$, numbered according to Figure \ref{f6}. 
Applying Lemma 2.5 to the exact sequence (2.1) of 
Subsection 2.2, one has that $\MM (F_{g,r+1}, \PP_n)$ is 
isomorphic with the quotient of the free product $PG 
(g,r,n) \ast A(A_{n-1})$ by the following relations.

$\bullet$\qua Relations from $\Sigma_n$:

$
{\rm(T15)} \quad V_i^2 = \Delta^2 (x_{r+i-1}, x_{r+i+1}, y_1) 
\Delta^{-1} (x_{r+i-1}, x_{r+i}, x_{r+i+1}, y_1) $

\hskip2cm ${\rm for}\ 1\le i\le n-1.$

$\bullet$\qua Relations from conjugation by the $V_i$'s:

\medskip
\centerline{\vbox{\halign{
$#$\ &$#$\hfill\cr
{\rm(T16)} &V_i w V_i^{-1} = w\ {\rm for}\ 1\le i\le n-1\ {\rm 
and}\cr
& w \in \{x_0, \dots, x_{r+i-1}, x_{r+i+1}, \dots, x_{r+n},y_1,
\dots, y_{2g-1}, z, u_1, \dots, u_r\}, \cr
\noalign{\smallskip}
{\rm(T17)} &V_i x_{r+i} V_i^{-1} = y_1 x_{r+i-1} x_{r+i}^{-1} 
y_1^{-1} x_{r+i+1} y_1 x_{r+i} x_{r+i-1}^{-1} y_1^{-1} 
\quad {\rm for}\ 1\le i\le n-1. \cr}}}

\medskip
We can easily prove using Proposition 2.12 that the
relation (T15) ``holds'' in 
\break
$\MM (F_{g,r+1}, \PP_n)$.
The relation (T16) is 
obvious, while the relation (T17) has to be verified by 
hand.

Now, the homomorphism $\bar g: PG (g,r,n) \to G(g,r,n)$ 
extends to a homomorphism $f: PG (g,r,n) \ast A(A_{n-1}) 
\to G(g,r,n)$ which sends $V_i$ on $v_i$ for all $i=1, \dots, 
n-1$.

\medskip\noindent
{\bf Assertion 4}\qua {\sl The homomorphism $f: PG(g,r,n) \ast 
A(A_{n-1}) \to G(g,r,n)$ induces a homomorphism $\bar f: \MM 
(F_{g,r+1}, \PP_n) \to G(g,r,n)$.}

\medskip
One can easily verify on the generators of $G(g,r,n)$ that 
$\bar f \circ \bar \rho$ is the identity of $G(g,r,n)$. So, 
Assertion 4 finishes the construction of $\bar f$ and the 
proof of Theorem 3.1.

\medskip\noindent
{\bf Proof of Assertion 4}\qua We have to show that: if 
$w_1=w_2$ is one of the relations (T15), (T16), (T17), then 
$f(w_1) = f(w_2)$.

$\bullet$\qua Relation (T15): 

\medskip
\halign{
$#$\ &$#$\hfill\cr
& f( \Delta^{-1} (x_{r+i-1}, x_{r+i}, x_{r+i+1}, y_1) 
\Delta^2 (x_{r+i-1}, x_{r+i+1}, y_1)) \cdot v_i^{-2} \cr
=& A_i^{-1} y_1^{-1} A_{i-1}^{-1} A_{i+1}^{-1} y_1^{-1} 
A_i^{-1} y_1 A_{i-1} A_{i+1} y_1 A_{i-1} A_{i+1} v_i^{-2} \cr
&\quad ({\rm by\ Propositions\ 2.8\ and\ 2.9})\cr
=& A_i^{-1} y_1^{-1} A_{i-1}^{-1} A_{i-1} v_i^{-1} A_i^{-1} 
v_i^{-1} A_i^{-1} y_1^{-1} A_i^{-1} y_1 A_{i-1} A_{i-1}^{-
1} A_i v_i A_i v_i y_1 A_{i-1} A_{i-1}^{-1} A_i v_i\cr
& A_i v_i v_i^{-2} ({\rm by\ (T1)})\cr
=& A_i^{-1} y_1^{-1} v_i^{-1} A_i^{-1} v_i^{-1} A_i^{-1} A_i y_1^{-
1} A_i^{-1}A_i v_i A_i v_i y_1 A_i v_i A_i v_i^{-1} \quad 
({\rm by\ (T2),(T3),(T4)})\cr
=& A_i^{-1} y_1^{-1} v_i^{-1} y_1 A_i^{-1} y_1^{-1} v_i y_1 
v_i^{-1} A_i v_i A_i  \quad ({\rm by\ (T1),\dots,(T4)})\cr
=& 1 \quad ({\rm by\ (T2),(T3),(T4)}). \cr}

\medskip
$\bullet$\qua The relation (T16) is a direct consequence of the 
braid relations in $A( \Gamma_{g,r,n})$.

$\bullet$\qua Relation (T17): 

\medskip
\halign{
$#$\ &$#$\hfill\cr
&f (y_1 x_{r+i-1} x_{r+i}^{-1} y_1^{-1} x_{r+i+1} y_1 
x_{r+i} x_{r+i-1}^{-1} y_1^{-1}) v_i f(x_{r+i}^{-1}) v_i^{-
1} \cr
=& y_1 A_{i-1} A_i^{-1} y_1^{-1} A_{i+1} y_1 A_i A_{i-1}^{-
1} y_1^{-1} v_i A_i^{-1} v_i^{-1}\cr
=& y_1 A_i^{-1} A_{i-1} y_1^{-1} A_{i-1}^{-1} A_i v_i A_i 
v_i y_1 A_i v_i A_{i-1}^{-1} y_1^{-1} A_i^{-1} v_i^{-1} 
\quad ({\rm by\ (T1),(T2),(T3)}) \cr
=& y_1 A_i^{-1} y_1^{-1} A_{i-1}^{-1} y_1 A_i v_i A_i 
v_i y_1 A_i v_i A_{i-1}^{-1} y_1^{-1} A_i^{-1} v_i^{-1} 
\quad ({\rm by\ (T4)}) \cr
=& A_i^{-1} y_1^{-1} A_i A_{i-1}^{-1} y_1 A_i v_i A_i 
v_i y_1 A_i v_i A_{i-1}^{-1} y_1^{-1} A_i^{-1} v_i^{-1} 
\quad ({\rm by\ (T4)}) \cr
=& A_i^{-1} y_1^{-1} A_{i-1}^{-1} y_1 A_i v_i y_1 A_i 
v_i y_1 A_i v_i A_{i-1}^{-1} y_1^{-1} A_i^{-1} v_i^{-1} 
\quad ({\rm by\ (T2),(T3),(T4)}) \cr
=& A_i^{-1} y_1^{-1} A_{i-1}^{-1} \cdot h_i (\Delta (x_1, 
x_2, x_3)) \cdot A_{i-1}^{-1} y_1^{-1} A_i^{-1} v_i^{-1} 
\quad ({\rm by\ Proposition\ 2.8}) \cr
=& A_i^{-1} y_1^{-1} A_{i-1}^{-1} A_{i-1} y_1 A_i v_i A_i 
y_1 A_{i-1} A_{i-1}^{-1} y_1^{-1} A_i^{-1} v_i^{-1}  
({\rm by\ (T5)\,\ Proposition\ 2.9}) \cr
=& 1. \cr}\vglue -2em\noindent\hbox{}\hfill\qed


\subsection{Proof of Theorem 3.2}

Let $c_1: S^1 \to \partial F_{g,1}$ be the boundary curve of 
$F_{g,1}$. We regard $F_{g,0}$ as obtained from $F_{g,1}$ by 
gluing a disk $D^2$ along $c_1$,
and we denote by $\varphi: \MM (F_{g,1}, 
\PP_n) \to \MM (F_{g,0}, \PP_n)$ the homomorphism induced 
by the inclusion of $F_{g,1}$ in $F_{g,0}$. The next 
proposition is the key of the proof of Theorem 3.2.

\medskip\noindent
{\bf Proposition 3.9}\qua {\sl {\rm(i)}\qua Let $g\ge 2$, and let $a_n, 
a_n'$ be the Dehn twists represented in Figure \ref{f20}. Then 
$\varphi$ is surjective and its kernel is the normal 
subgroup of $\MM (F_{g,1}, \PP_n)$ normaly generated by $\{ 
a_n^{-1} a_n' \}$.

\smallskip
{\rm(ii)}\qua Let $g=1$, and let $e,e'$ be the Dehn twists 
represented in Figure \ref{f20}. Then $\varphi$ is surjective and 
its kernel is the normal subgroup of $\MM (F_{1,1}, \PP_n)$ 
normaly generated by $\{ a_n^{-1} a_0, e^{-1} e' \}$.}

\begin{figure}[ht]
\centerline{\input{pmc_20}}
\caption{\label{f20} Relations in $\MM (F_{g,0}, \PP_n)$}
\end{figure}

\medskip\noindent
{\bf Proof}\qua We choose a point $Q$ in the interior of the 
disk $D^2$, and we denote by $\MM_Q (F_{g,0}, \PP_n \cup \{ 
Q \})$ the subgroup of $\MM (F_{g,0}, \PP_n \cup \{ Q \})$ 
of isotopy classes of elements of $\HH (F_{g,0}, \PP_n \cup 
\{Q\})$ that fix $Q$. An easy algebraic argument on the 
exact sequences (2.1), (2.2), and (2.3) of Subsection 2.2 
shows that we have the following exact sequences.
$$
\displaylines{
(2.2.a) \hfill 1 \to \pi_1 (F_{g,0} \setminus \PP_n, Q) \to 
\MM_Q (F_{g,0}, \PP_n \cup \{Q\}) \arrow{ \varphi_1}{} \MM 
(F_{g,0}, \PP_n) \to 1, \hfill \cr
(2.3.a) \hfill 1 \to \BZ \to \MM (F_{g,1}, \PP_n) \arrow{ 
\varphi_2}{} \MM_Q (F_{g,0}, \PP_n \cup \{Q\}) \to 1. 
\hfill \cr}
$$
Moreover, we have $\varphi = \varphi_1 \circ \varphi_2$.

A first consequence of these exact sequences is that 
$\varphi$ is surjective. Now, we use them for finding a 
normal generating set of $\ker \varphi$. 

The group $\pi_1 (F_{g,0} \setminus \PP_n, Q)$ is the free 
group freely generated by the loops $\bar \alpha_1, \dots, 
\bar \alpha_n$, $\bar \beta_1, \dots, \bar \beta_{2g-1}$ 
represented in Figure \ref{f21}. One can easily verify by hand that 
the following equalities hold in $\MM_Q (F_{g,0}, \PP_n \cup 
\{Q\})$:

\medskip
\centerline{\vbox{\halign{
\hfill$#$&\ $#$\ &$#$\quad\hfill&$#$\hfill\cr
\bar \alpha_i &=& \varphi_2 (b_1 a_n' a_i b_1 a_n)^{-1} 
\cdot \bar \alpha_n^{-1} \cdot \varphi_2 (b_1 a_n' a_i b_1 
a_n) &{\rm for}\ i=1, \dots, n-1,\cr
\bar \beta_1 &=& \varphi_2 (b_1 a_n)^{-1} \cdot \bar 
\alpha_n \cdot \varphi_2 (b_1 a_n), \cr
\bar \beta_j &=& \varphi_2 (b_j b_{j-1})^{-1} \cdot 
\bar \beta_{j-1} \cdot \varphi_2 (b_j b_{j-1}) &{\rm 
for}\ j= 2, \dots, 2g-1. \cr}}} 

\medskip\noindent
Moreover, by Lemma 2.6, we have:
$$
\bar \alpha_n = \varphi_2 (a_n^{-1} a_n').
$$
On the other hand, by Lemma 2.7, the Dehn twists $\sigma_1$ along 
the boundary curve of $F_{g,1}$ generates the kernel 
of $\varphi_2$.  So, the kernel of $\varphi$ is the normal 
subgroup normaly generated by $\{ a_n^{-1} a_n', \sigma_1\}$.

\begin{figure}[ht]
\centerline{\input{pmc_21}}
\caption{\label{f21} Generators of $\pi_1 (F_{g,0} \setminus
\PP_n, Q)$}
\end{figure}
 
\medskip 
Now, assume $g\ge 2$. Let $G'$ denote the quotient of $\MM 
(F_{g,1}, \PP_n)$ by the relation $a_n=a_n'$. Define a {\it 
spinning pair} of Dehn twists to be a pair $(\sigma, \sigma')$
of Dehn twists 
conjugated to $(a_n,a_n')$, namely, a pair $(\sigma, 
\sigma')$ of Dehn twists satisfying: there exists $\xi \in 
\MM (F_{g,1}, \PP_n)$ such that $\sigma = \xi a_n \xi^{-1}$ 
and $\sigma' = \xi a_n' \xi^{-1}$. Note that we have the 
equality $\sigma = \sigma'$ in $G'$ if $(\sigma, \sigma')$ 
is a spinning pair. Consider the Dehn twists $e_1, e_2, e_3, 
e_1', e_2', e_3'$ represented in Figure \ref{f22}. The pairs $(e_1, 
e_1')$, $(e_2, e_2')$, $(e_3, e_3')$ are spinning pairs, thus 
we have the equalities $e_1=e_1'$, $e_2 = e_2'$, $e_3 = 
e_3'$ in $G'$. Moreover, the lantern relation of Lemma 2.4 
implies:
$$
e_1 e_2 e_3 \sigma_1 = e_1' e_2' e_3'.
$$
Thus, the equality $\sigma_1=1$ holds in $G'$. This shows that 
the kernel of $\varphi$ is the normal subgroup of $\MM 
(F_{g,1}, \PP_n)$ normaly generated by $\{a_n^{-1} a_n'\}$.

\begin{figure}[ht]
\centerline{\input{pmc_22}}
\caption{\label{f22} Lantern relation in $\MM (F_{g,1}, 
\PP_n)$}
\end{figure}

\medskip
Now, we assume $g=1$. Then $a_n'=a_0$. Let $G'$ be the 
quotient of $\MM (F_{1,1}, \PP_n)$ by the relation 
$a_n=a_0$. By Proposition 2.12, we have the following 
equalities in $G'$.

\medskip
\centerline{\vbox{\halign{
\hfill$#$&\ $#$\ &$#$\hfill\cr
\sigma_1 e &=& (a_0 b_1 a_n a_0 b_1 a_0)^2 = (a_0 b_1 a_0 a_0 
b_1 a_0)^2,\cr
e' &=& (a_0 b_1 a_0)^4. \cr}}}

\medskip\noindent
Thus, we have the equality $\sigma_1=e^{-1} e'$ in $G'$. So, the 
kernel of $\varphi$ is the normal subgroup of $\MM (F_{1,1}, 
\PP_n)$ normaly generated by $\{ a_n^{-1} a_0, e^{-1} e'\}$. 
\qed

\medskip\noindent
{\bf Proof of Theorem 3.2}\qua Recall that $\Gamma_{g,0,n}$ 
denotes the Coxeter graph drawn in Figure \ref{f16}, and that 
$G(g,0,n)$ denotes the quotient of $A(\Gamma_{g,0,n})$ by 
the relations (R1), (R2), (R7), (R8a). By Theorem 3.1, there 
is an isomorphism $\bar \rho: G(g,0,n) \to \MM (F_{g,1}, 
\PP_n)$ which sends $x_i$ on $a_i$ for $i=0,1$, $y_i$ on 
$b_i$ for $i=1, \dots, 2g-1$, $z$ on $c$, and $v_i$ on 
$\tau_i$ for $i=1, \dots, n-1$.

First, assume $g\ge 2$. Let $G_0(g,n)$ denote the quotient 
of $G(g,0,n)$ by the relation (R9a). Proposition 2.12 
implies:

\medskip
\centerline{\vbox{\halign{
\hfill$#$&\ $#$\ &$#$\hfill\cr
a_n &=& \bar \rho (x_0^{1-n} \Delta (x_1, v_1, \dots, v_{n-
1})),\cr
\noalign{\smallskip}
a_n' &=& \bar \rho (x_0^{3-2g} \Delta (z, y_2, \dots, 
y_{2g-1})).\cr}}}

\medskip\noindent
Thus, by Proposition 3.9, $\bar \rho$ induces an isomorphism :
$$\bar \rho_0: G_0(g,n) \to \MM (F_{g,0}, \PP_n).$$

Now, assume $g=1$. Let $G_0(1,n)$ denote the quotient of 
$G(1,0,n)$ by the relations (R9b), (R9c). Proposition 2.12 
implies:

\medskip
\centerline{\vbox{\halign{
\hfill$#$&\ $#$\ &$#$\hfill\cr
a_n &=& \bar \rho (x_0^{1-n} \Delta (x_1, v_1, \dots, v_{n-
1})),\cr
\noalign{\smallskip}
e &=& \bar \rho (\Delta^2 (v_1, \dots, v_{n-1})),\cr
\noalign{\smallskip}
e' &=& \bar \rho (\Delta^4 (x_0, y_1)).\cr}}}

\medskip\noindent
Thus, by Proposition 3.9, $\bar\rho$ induces an isomorphism : 
$$\bar \rho_0:G_0(1,n) \to \MM (F_{1,0}, \PP_n).\eqno{\qed}$$


\Addresses\recd

\end{document}

%% file: agtout.tex
\def\ifplaintex{\expandafter\ifx\csname documentclass\endcsname\relax}

\def\gtp{{\mathsurround=0pt\it $\cal G\mskip-2mu$eometry \&\ 
$\cal T\!\!$opology $\cal P\!$ublications}}  

\def\recd{{\small Received:\qua\receiveddate\ifx\reviseddate\relax
\else\qquad Revised:\qua\reviseddate\fi\par}} 

\def\lognumber#1{\def\thelognumber{#1}}
\def\volumenumber#1{\def\thevolumenumber{#1}}
\def\volumeyear#1{\def\thevolumeyear{#1}}
\def\papernumber#1{\def\thepapernumber{#1}}
\def\pagenumbers#1#2{\def\startpage{#1}\def\finishpage{#2}}
\def\published#1{\def\publishdate{#1}}

\def\received#1{\def\receiveddate{#1}}

\def\accepted#1{\def\accepteddate{#1}}

\def\asciiauthors#1{\def\theasciiauthors{#1}}
\def\asciiaddress#1{\def\theasciiaddress{#1}}

\def\coverauthors#1{\def\thecoverauthors{#1}}
\long\def\asciiabstract#1{\long\def\theasciiabstract{#1}}

\let\\\par\let\thelognumber\relax\let\thevolumenumber\relax
\let\thepapernumber\relax\let\thevolumeyear\relax\let\startpage\relax
\let\finishpage\relax\let\publishdate\relax\let\receiveddate\relax
\let\reviseddate\relax\let\accepteddate\relax\let\theasciititle\relax
\let\theasciiauthors\relax\let\theasciiaddress\relax
\let\theasciiabstract\relax

\let\thecoverauthors\relax\let\theasciiemail\relax

\ifplaintex
\font\logobig=cmssbx10 scaled 3836
\font\logomed=cmssbx10 scaled 2557
\else
\font\logobig=cmssbx10 scaled 4200
\font\logomed=cmssbx10 scaled 2800
\fi

\long\def\makeagttitle{   
\count0=\startpage
\agt\hfill      
\hbox to 45truept{\vbox to 0pt{\vglue -13truept{\logomed A\kern -.37em{\logobig 
T}\kern -.38em G}\vss}\hss}
\break
{\small Volume \thevolumenumber\ (\thevolumeyear)
\startpage--\finishpage\nl
Published: \publishdate}

\vglue .25truein

{\parskip=0pt\leftskip 0pt plus
1fil\def\\{\par\smallskip}{\Large\bf\thetitle}\par\medskip} \vglue
0.05truein

{\parskip=0pt\leftskip 0pt plus 1fil\def\\{\par}{\sc\theauthors}
\par\medskip}%
 
\vglue 0.03truein 

{\small\leftskip 25truept\rightskip 25truept{\bf Abstract}\stdspace\theabstract

{\bf AMS Classification}\stdspace\theprimaryclass
\ifx\thesecondaryclass\relax\else; \thesecondaryclass\fi\par
{\bf Keywords}\stdspace \thekeywords\par}\vglue 7truept

}   

\ifplaintex
\font\phead=cmsl9 scaled 950
\font\pnum=cmbx10 scaled 913
\font\pfoot=cmsl9 scaled 950
\headline{\vbox to 0pt{\vskip -4.5mm\line{\small\phead\ifnum
\count0=\startpage ISSN 1472-2739 (on-line) 1472-2747 (printed)
\hfill {\pnum\folio}\else\ifodd\count0\def\\{ }%
\ifx\theshorttitle\relax\thetitle\else\theshorttitle\fi\hfill{\pnum\folio}
\else\def\\{ and }{\pnum\folio}\hfill\ifx\theshortauthors\relax\theauthors
\else\theshortauthors\fi\fi\fi}\vss}}
\footline{\vbox to 0pt{\vglue 0mm\line{\small\pfoot\ifnum\count0=\startpage
\copyright\ \gtp\hfill\else
\agt, Volume \thevolumenumber\ (\thevolumeyear)\hfill\fi}\vss}}
\else
\font=cmsl9 scaled 1050
\font\lnum=cmbx10 
\font=cmsl9 scaled 1050
\makeatletter
\def\@oddhead{{\small\lhead\ifnum\count0=\startpage ISSN 1472-2739 
(on-line) 1472-2747 (printed)\hfill {\lnum\number\count0}\else\ifodd\count0
\def\\{ }\ifx\theshorttitle\relax \thetitle \else\theshorttitle\fi\hfill
{\lnum\number\count0}\else\def\\{ and }{\lnum\number\count0}
\hfill\ifx\theshortauthors\relax 
\theauthors\else\theshortauthors\fi\fi\fi}}\def\@evenhead{\@oddhead}
\def\@oddfoot{\small\lfoot\ifnum\count0=\startpage\copyright\ \gtp\hfill\else
\agt, Volume \thevolumenumber\ (\thevolumeyear)\hfill\fi}
\def\@evenfoot{\@oddfoot}
\makeatother
\fi
\let\maketitlepage\makeagttitle
\let\makeshorttitle\maketitlepage
\let\maketitle\maketitlepage

\newwrite\gtoutfile
\long\gdef\makeheadfile{  
{\def\\{, }\def\s{ }
\immediate\openout\gtoutfile head.xxx
\immediate\write\gtoutfile{To: math@arxiv.org}
\immediate\write\gtoutfile{Subject: put OR rep NNNNN:ppppp}
\immediate\write\gtoutfile{--text follows this line--}
\immediate\write\gtoutfile{Proxy-for: \ifx\theasciiauthors\relax
\theauthors\else\theasciiauthors\fi\s<\ifx\theasciiemail\relax\theemail\else\theasciiemail\fi>}
\immediate\write\gtoutfile{\noexpand\\}
\immediate\write\gtoutfile{Authors: \ifx\theasciiauthors\relax
\theauthors\else\theasciiauthors\fi}
{\def\\{ }\immediate\write\gtoutfile{Title: \ifx\theasciititle\relax
\thetitle\else\theasciititle\fi}}
\immediate\write\gtoutfile{Subj-class: GT or SG, GR etc}
\immediate\write\gtoutfile{MSC-class: \theprimaryclass\ifx\thesecondaryclass\relax\else, \thesecondaryclass\fi}
\immediate\write\gtoutfile{Journal-ref: Algebraic and Geometric Topology \thevolumenumber\s
(\thevolumeyear) \startpage-\finishpage}
\immediate\write\gtoutfile{Comments: Published by Algebraic and
Geometric Topology at}
\immediate\write\gtoutfile{\s\s\s  http://www.maths.warwick.ac.uk/agt/AGTVol\thevolumenumber/agt-\thevolumenumber-\thepapernumber.abs.html}
\immediate\write\gtoutfile{\noexpand\\}
\immediate\write\gtoutfile{}
\ifx\theasciiabstract\relax
\immediate\write\gtoutfile{\theabstract}\else
\immediate\write\gtoutfile{\theasciiabstract}\fi
\immediate\write\gtoutfile{}
\immediate\write\gtoutfile{\noexpand\\}
\immediate\write\gtoutfile{}
\immediate\closeout\gtoutfile}}  

\def\maketitlepage{\makeagttitle\makeheadfile}
\let\makeshorttitle\maketitlepage
\let\maketitle\maketitlepage

\def\ifplaintex{\expandafter\ifx\csname documentclass\endcsname\relax}

\def\gtp{{\mathsurround=0pt\it $\cal G\mskip-2mu$eometry \&\ 
$\cal T\!\!$opology $\cal P\!$ublications}}  

\def\recd{{\small Received:\qua\receiveddate\ifx\reviseddate\relax
\else\qquad Revised:\qua\reviseddate\fi\par}} 

\def\lognumber#1{\def\thelognumber{#1}}
\def\volumenumber#1{\def\thevolumenumber{#1}}
\def\volumeyear#1{\def\thevolumeyear{#1}}
\def\papernumber#1{\def\thepapernumber{#1}}
\def\pagenumbers#1#2{\def\startpage{#1}\def\finishpage{#2}}
\def\published#1{\def\publishdate{#1}}

\def\received#1{\def\receiveddate{#1}}

\def\accepted#1{\def\accepteddate{#1}}

\def\asciiauthors#1{\def\theasciiauthors{#1}}
\def\asciiaddress#1{\def\theasciiaddress{#1}}

\def\coverauthors#1{\def\thecoverauthors{#1}}
\long\def\asciiabstract#1{\long\def\theasciiabstract{#1}}

\let\\\par\let\thelognumber\relax\let\thevolumenumber\relax
\let\thepapernumber\relax\let\thevolumeyear\relax\let\startpage\relax
\let\finishpage\relax\let\publishdate\relax\let\receiveddate\relax
\let\reviseddate\relax\let\accepteddate\relax\let\theasciititle\relax
\let\theasciiauthors\relax\let\theasciiaddress\relax
\let\theasciiabstract\relax

\let\thecoverauthors\relax\let\theasciiemail\relax

\ifplaintex
\font\logobig=cmssbx10 scaled 3836
\font\logomed=cmssbx10 scaled 2557
\else
\font\logobig=cmssbx10 scaled 4200
\font\logomed=cmssbx10 scaled 2800
\fi

\long\def\makeagttitle{   
\count0=\startpage
\agt\hfill      
\hbox to 45truept{\vbox to 0pt{\vglue -13truept{\logomed A\kern -.37em{\logobig 
T}\kern -.38em G}\vss}\hss}
\break
{\small Volume \thevolumenumber\ (\thevolumeyear)
\startpage--\finishpage\nl
Published: \publishdate}

\vglue .25truein

{\parskip=0pt\leftskip 0pt plus
1fil\def\\{\par\smallskip}{\Large\bf\thetitle}\par\medskip} \vglue
0.05truein

{\parskip=0pt\leftskip 0pt plus 1fil\def\\{\par}{\sc\theauthors}
\par\medskip}%
 
\vglue 0.03truein 

{\small\leftskip 25truept\rightskip 25truept{\bf Abstract}\stdspace\theabstract

{\bf AMS Classification}\stdspace\theprimaryclass
\ifx\thesecondaryclass\relax\else; \thesecondaryclass\fi\par
{\bf Keywords}\stdspace \thekeywords\par}\vglue 7truept

}   

\ifplaintex
\font\phead=cmsl9 scaled 950
\font\pnum=cmbx10 scaled 913
\font\pfoot=cmsl9 scaled 950
\headline{\vbox to 0pt{\vskip -4.5mm\line{\small\phead\ifnum
\count0=\startpage ISSN 1472-2739 (on-line) 1472-2747 (printed)
\hfill {\pnum\folio}\else\ifodd\count0\def\\{ }%
\ifx\theshorttitle\relax\thetitle\else\theshorttitle\fi\hfill{\pnum\folio}
\else\def\\{ and }{\pnum\folio}\hfill\ifx\theshortauthors\relax\theauthors
\else\theshortauthors\fi\fi\fi}\vss}}
\footline{\vbox to 0pt{\vglue 0mm\line{\small\pfoot\ifnum\count0=\startpage
\copyright\ \gtp\hfill\else
\agt, Volume \thevolumenumber\ (\thevolumeyear)\hfill\fi}\vss}}
\else
\font=cmsl9 scaled 1050
\font\lnum=cmbx10 
\font=cmsl9 scaled 1050
\makeatletter
\def\@oddhead{{\small\lhead\ifnum\count0=\startpage ISSN 1472-2739 
(on-line) 1472-2747 (printed)\hfill {\lnum\number\count0}\else\ifodd\count0
\def\\{ }\ifx\theshorttitle\relax \thetitle \else\theshorttitle\fi\hfill
{\lnum\number\count0}\else\def\\{ and }{\lnum\number\count0}
\hfill\ifx\theshortauthors\relax 
\theauthors\else\theshortauthors\fi\fi\fi}}\def\@evenhead{\@oddhead}
\def\@oddfoot{\small\lfoot\ifnum\count0=\startpage\copyright\ \gtp\hfill\else
\agt, Volume \thevolumenumber\ (\thevolumeyear)\hfill\fi}
\def\@evenfoot{\@oddfoot}
\makeatother
\fi
\let\maketitlepage\makeagttitle
\let\makeshorttitle\maketitlepage
\let\maketitle\maketitlepage

\newwrite\gtoutfile
\long\gdef\makeheadfile{  
{\def\\{, }\def\s{ }
\immediate\openout\gtoutfile head.xxx
\immediate\write\gtoutfile{To: math@arxiv.org}
\immediate\write\gtoutfile{Subject: put OR rep NNNNN:ppppp}
\immediate\write\gtoutfile{--text follows this line--}
\immediate\write\gtoutfile{Proxy-for: \ifx\theasciiauthors\relax
\theauthors\else\theasciiauthors\fi\s<\ifx\theasciiemail\relax\theemail\else\theasciiemail\fi>}
\immediate\write\gtoutfile{\noexpand\\}
\immediate\write\gtoutfile{Authors: \ifx\theasciiauthors\relax
\theauthors\else\theasciiauthors\fi}
{\def\\{ }\immediate\write\gtoutfile{Title: \ifx\theasciititle\relax
\thetitle\else\theasciititle\fi}}
\immediate\write\gtoutfile{Subj-class: GT or SG, GR etc}
\immediate\write\gtoutfile{MSC-class: \theprimaryclass\ifx\thesecondaryclass\relax\else, \thesecondaryclass\fi}
\immediate\write\gtoutfile{Journal-ref: Algebraic and Geometric Topology \thevolumenumber\s
(\thevolumeyear) \startpage-\finishpage}
\immediate\write\gtoutfile{Comments: Published by Algebraic and
Geometric Topology at}
\immediate\write\gtoutfile{\s\s\s  http://www.maths.warwick.ac.uk/agt/AGTVol\thevolumenumber/agt-\thevolumenumber-\thepapernumber.abs.html}
\immediate\write\gtoutfile{\noexpand\\}
\immediate\write\gtoutfile{}
\ifx\theasciiabstract\relax
\immediate\write\gtoutfile{\theabstract}\else
\immediate\write\gtoutfile{\theasciiabstract}\fi
\immediate\write\gtoutfile{}
\immediate\write\gtoutfile{\noexpand\\}
\immediate\write\gtoutfile{}
\immediate\closeout\gtoutfile}}  

\def\maketitlepage{\makeagttitle\makeheadfile}
\let\makeshorttitle\maketitlepage
\let\maketitle\maketitlepage

%% file: pmc_1.tex
\begin{picture}(0,0)%
\special{psfile=pmc_1.ps}%
\end{picture}%
\setlength{\unitlength}{1579sp}%
\begingroup\makeatletter\ifx\SetFigFont\undefined
\def\x#1#2#3#4#5#6#7\relax{\def\x{#1#2#3#4#5#6}}%
\expandafter\x\fmtname xxxxxx\relax \def\y{splain}%
\ifx\x\y   
\gdef\SetFigFont#1#2#3{%
  \ifnum #1<17\tiny\else \ifnum #1<20\small\else
  \ifnum #1<24\normalsize\else \ifnum #1<29\large\else
  \ifnum #1<34\Large\else \ifnum #1<41\LARGE\else
     \huge\fi\fi\fi\fi\fi\fi
  \csname #3\endcsname}%
\else
\gdef\SetFigFont#1#2#3{\begingroup
  \count@#1\relax \ifnum 25<\count@\count@25\fi
  \def\x{\endgroup\@setsize\SetFigFont{#2pt}}%
  \expandafter\x
    \csname \romannumeral\the\count@ pt\expandafter\endcsname
    \csname @\romannumeral\the\count@ pt\endcsname
  \csname #3\endcsname}%
\fi
\fi\endgroup
\begin{picture}(8160,2758)(3571,-3690)
\put(4636,-1276){\makebox(0,0)[lb]{\smash{\SetFigFont{10}{12.0}{rm}$s$}}}
\end{picture}

%% file: pmc_2.tex
\begin{picture}(0,0)%
\special{psfile=pmc_2.ps}%
\end{picture}%
\setlength{\unitlength}{1973sp}%
\begingroup\makeatletter\ifx\SetFigFont\undefined
\def\x#1#2#3#4#5#6#7\relax{\def\x{#1#2#3#4#5#6}}%
\expandafter\x\fmtname xxxxxx\relax \def\y{splain}%
\ifx\x\y   
\gdef\SetFigFont#1#2#3{%
  \ifnum #1<17\tiny\else \ifnum #1<20\small\else
  \ifnum #1<24\normalsize\else \ifnum #1<29\large\else
  \ifnum #1<34\Large\else \ifnum #1<41\LARGE\else
     \huge\fi\fi\fi\fi\fi\fi
  \csname #3\endcsname}%
\else
\gdef\SetFigFont#1#2#3{\begingroup
  \count@#1\relax \ifnum 25<\count@\count@25\fi
  \def\x{\endgroup\@setsize\SetFigFont{#2pt}}%
  \expandafter\x
    \csname \romannumeral\the\count@ pt\expandafter\endcsname
    \csname @\romannumeral\the\count@ pt\endcsname
  \csname #3\endcsname}%
\fi
\fi\endgroup
\begin{picture}(6360,2414)(3571,-3240)
\put(4321,-1141){\makebox(0,0)[lb]{\smash{\SetFigFont{10}{14.4}{rm}$s$}}}
\put(4366,-2131){\makebox(0,0)[lb]{\smash{\SetFigFont{10}{14.4}{rm}$a$}}}
\end{picture}

%% file: pmc_3.tex
\begin{picture}(0,0)%
\special{psfile=pmc_3.ps}%
\end{picture}%
\setlength{\unitlength}{1973sp}%
\begingroup\makeatletter\ifx\SetFigFont\undefined
\def\x#1#2#3#4#5#6#7\relax{\def\x{#1#2#3#4#5#6}}%
\expandafter\x\fmtname xxxxxx\relax \def\y{splain}%
\ifx\x\y   
\gdef\SetFigFont#1#2#3{%
  \ifnum #1<17\tiny\else \ifnum #1<20\small\else
  \ifnum #1<24\normalsize\else \ifnum #1<29\large\else
  \ifnum #1<34\Large\else \ifnum #1<41\LARGE\else
     \huge\fi\fi\fi\fi\fi\fi
  \csname #3\endcsname}%
\else
\gdef\SetFigFont#1#2#3{\begingroup
  \count@#1\relax \ifnum 25<\count@\count@25\fi
  \def\x{\endgroup\@setsize\SetFigFont{#2pt}}%
  \expandafter\x
    \csname \romannumeral\the\count@ pt\expandafter\endcsname
    \csname @\romannumeral\the\count@ pt\endcsname
  \csname #3\endcsname}%
\fi
\fi\endgroup
\begin{picture}(3960,4185)(4726,-4606)
\put(6526,-781){\makebox(0,0)[lb]{\smash{\SetFigFont{10}{14.4}{rm}$e_4$}}}
\put(4726,-2941){\makebox(0,0)[lb]{\smash{\SetFigFont{10}{14.4}{rm}$e_1$}}}
\put(5941,-2896){\makebox(0,0)[lb]{\smash{\SetFigFont{10}{14.4}{rm}$b$}}}
\put(6526,-4606){\makebox(0,0)[lb]{\smash{\SetFigFont{10}{14.4}{rm}$e_2$}}}
\put(8686,-2896){\makebox(0,0)[lb]{\smash{\SetFigFont{10}{14.4}{rm}$e_3$}}}
\put(6706,-3751){\makebox(0,0)[lb]{\smash{\SetFigFont{10}{14.4}{rm}$c$}}}
\put(7561,-2851){\makebox(0,0)[lb]{\smash{\SetFigFont{10}{14.4}{rm}$a$}}}
\end{picture}

%% file: pmc_4.tex
\begin{picture}(0,0)%
\special{psfile=pmc_4.ps}%
\end{picture}%
\setlength{\unitlength}{1973sp}%
\begingroup\makeatletter\ifx\SetFigFont\undefined
\def\x#1#2#3#4#5#6#7\relax{\def\x{#1#2#3#4#5#6}}%
\expandafter\x\fmtname xxxxxx\relax \def\y{splain}%
\ifx\x\y   
\gdef\SetFigFont#1#2#3{%
  \ifnum #1<17\tiny\else \ifnum #1<20\small\else
  \ifnum #1<24\normalsize\else \ifnum #1<29\large\else
  \ifnum #1<34\Large\else \ifnum #1<41\LARGE\else
     \huge\fi\fi\fi\fi\fi\fi
  \csname #3\endcsname}%
\else
\gdef\SetFigFont#1#2#3{\begingroup
  \count@#1\relax \ifnum 25<\count@\count@25\fi
  \def\x{\endgroup\@setsize\SetFigFont{#2pt}}%
  \expandafter\x
    \csname \romannumeral\the\count@ pt\expandafter\endcsname
    \csname @\romannumeral\the\count@ pt\endcsname
  \csname #3\endcsname}%
\fi
\fi\endgroup
\begin{picture}(2999,2788)(4472,-4586)
\put(6841,-3031){\makebox(0,0)[lb]{\smash{\SetFigFont{10}{14.4}{rm}$P_n$}}}
\put(5581,-2086){\makebox(0,0)[lb]{\smash{\SetFigFont{10}{14.4}{rm}$\alpha$}}}
\put(5896,-3256){\makebox(0,0)[lb]{\smash{\SetFigFont{10}{14.4}{rm}$\sigma_0$}}}
\put(7471,-3346){\makebox(0,0)[lb]{\smash{\SetFigFont{10}{14.4}{rm}$\sigma_1$}}}
\end{picture}

%% file: pmc_5.tex
\begin{picture}(0,0)%
\special{psfile=pmc_5.ps}%
\end{picture}%
\setlength{\unitlength}{1381sp}%
\begingroup\makeatletter\ifx\SetFigFont\undefined
\def\x#1#2#3#4#5#6#7\relax{\def\x{#1#2#3#4#5#6}}%
\expandafter\x\fmtname xxxxxx\relax \def\y{splain}%
\ifx\x\y   
\gdef\SetFigFont#1#2#3{%
  \ifnum #1<17\tiny\else \ifnum #1<20\small\else
  \ifnum #1<24\normalsize\else \ifnum #1<29\large\else
  \ifnum #1<34\Large\else \ifnum #1<41\LARGE\else
     \huge\fi\fi\fi\fi\fi\fi
  \csname #3\endcsname}%
\else
\gdef\SetFigFont#1#2#3{\begingroup
  \count@#1\relax \ifnum 25<\count@\count@25\fi
  \def\x{\endgroup\@setsize\SetFigFont{#2pt}}%
  \expandafter\x
    \csname \romannumeral\the\count@ pt\expandafter\endcsname
    \csname @\romannumeral\the\count@ pt\endcsname
  \csname #3\endcsname}%
\fi
\fi\endgroup
\begin{picture}(8487,6847)(2956,-8207)
\put(4366,-2671){\makebox(0,0)[lb]{\smash{\SetFigFont{8}{9.6}{rm}$\sigma_1$}}}
\put(5986,-5641){\makebox(0,0)[lb]{\smash{\SetFigFont{10}{13.2}{rm}$S_1$}}}
\put(6841,-4876){\makebox(0,0)[lb]{\smash{\SetFigFont{10}{13.2}{rm}$S_2$}}}
\put(8281,-5866){\makebox(0,0)[lb]{\smash{\SetFigFont{10}{13.2}{rm}$S_3$}}}
\put(6931,-1771){\makebox(0,0)[lb]{\smash{\SetFigFont{8}{9.6}{rm}$\sigma_2$}}}
\put(9901,-3166){\makebox(0,0)[lb]{\smash{\SetFigFont{8}{9.6}{rm}$\sigma_3$}}}
\end{picture}

%% file: pmc_6.tex
\begin{picture}(0,0)%
\special{psfile=pmc_6.ps}%
\end{picture}%
\setlength{\unitlength}{1579sp}%
\begingroup\makeatletter\ifx\SetFigFont\undefined
\def\x#1#2#3#4#5#6#7\relax{\def\x{#1#2#3#4#5#6}}%
\expandafter\x\fmtname xxxxxx\relax \def\y{splain}%
\ifx\x\y   
\gdef\SetFigFont#1#2#3{%
  \ifnum #1<17\tiny\else \ifnum #1<20\small\else
  \ifnum #1<24\normalsize\else \ifnum #1<29\large\else
  \ifnum #1<34\Large\else \ifnum #1<41\LARGE\else
     \huge\fi\fi\fi\fi\fi\fi
  \csname #3\endcsname}%
\else
\gdef\SetFigFont#1#2#3{\begingroup
  \count@#1\relax \ifnum 25<\count@\count@25\fi
  \def\x{\endgroup\@setsize\SetFigFont{#2pt}}%
  \expandafter\x
    \csname \romannumeral\the\count@ pt\expandafter\endcsname
    \csname @\romannumeral\the\count@ pt\endcsname
  \csname #3\endcsname}%
\fi
\fi\endgroup
\begin{picture}(12060,5279)(1756,-5565)
\put(12511,-5506){\makebox(0,0)[lb]{\smash{\SetFigFont{10}{12.0}{rm}$x_7$}}}
\put(1756,-781){\makebox(0,0)[lb]{\smash{\SetFigFont{10}{12.0}{rm}$A_l$}}}
\put(3511,-1096){\makebox(0,0)[lb]{\smash{\SetFigFont{10}{12.0}{rm}$x_2$}}}
\put(5761,-1096){\makebox(0,0)[lb]{\smash{\SetFigFont{10}{12.0}{rm}$x_l$}}}
\put(2611,-1096){\makebox(0,0)[lb]{\smash{\SetFigFont{10}{12.0}{rm}$x_1$}}}
\put(7876,-781){\makebox(0,0)[lb]{\smash{\SetFigFont{10}{12.0}{rm}$B_l$}}}
\put(12871,-1096){\makebox(0,0)[lb]{\smash{\SetFigFont{10}{12.0}{rm}$x_l$}}}
\put(10576,-1096){\makebox(0,0)[lb]{\smash{\SetFigFont{10}{12.0}{rm}$x_3$}}}
\put(9676,-1096){\makebox(0,0)[lb]{\smash{\SetFigFont{10}{12.0}{rm}$x_2$}}}
\put(8731,-1096){\makebox(0,0)[lb]{\smash{\SetFigFont{10}{12.0}{rm}$x_1$}}}
\put(1846,-2806){\makebox(0,0)[lb]{\smash{\SetFigFont{10}{12.0}{rm}$D_l$}}}
\put(2296,-2266){\makebox(0,0)[lb]{\smash{\SetFigFont{10}{12.0}{rm}$x_1$}}}
\put(2251,-3301){\makebox(0,0)[lb]{\smash{\SetFigFont{10}{12.0}{rm}$x_2$}}}
\put(1801,-4606){\makebox(0,0)[lb]{\smash{\SetFigFont{10}{12.0}{rm}$E_6$}}}
\put(9181,-601){\makebox(0,0)[lb]{\smash{\SetFigFont{10}{12.0}{rm}$4$}}}
\put(6796,-2986){\makebox(0,0)[lb]{\smash{\SetFigFont{10}{12.0}{rm}$x_l$}}}
\put(4501,-2986){\makebox(0,0)[lb]{\smash{\SetFigFont{10}{12.0}{rm}$x_4$}}}
\put(3646,-2986){\makebox(0,0)[lb]{\smash{\SetFigFont{10}{12.0}{rm}$x_3$}}}
\put(3376,-4291){\makebox(0,0)[lb]{\smash{\SetFigFont{10}{12.0}{rm}$x_2$}}}
\put(2476,-4291){\makebox(0,0)[lb]{\smash{\SetFigFont{10}{12.0}{rm}$x_1$}}}
\put(8101,-4606){\makebox(0,0)[lb]{\smash{\SetFigFont{10}{12.0}{rm}$E_7$}}}
\put(6256,-4291){\makebox(0,0)[lb]{\smash{\SetFigFont{10}{12.0}{rm}$x_5$}}}
\put(5311,-4291){\makebox(0,0)[lb]{\smash{\SetFigFont{10}{12.0}{rm}$x_4$}}}
\put(4366,-4291){\makebox(0,0)[lb]{\smash{\SetFigFont{10}{12.0}{rm}$x_3$}}}
\put(9991,-4246){\makebox(0,0)[lb]{\smash{\SetFigFont{10}{12.0}{rm}$x_2$}}}
\put(9091,-4246){\makebox(0,0)[lb]{\smash{\SetFigFont{10}{12.0}{rm}$x_1$}}}
\put(4906,-5551){\makebox(0,0)[lb]{\smash{\SetFigFont{10}{12.0}{rm}$x_6$}}}
\put(13816,-4246){\makebox(0,0)[lb]{\smash{\SetFigFont{10}{12.0}{rm}$x_6$}}}
\put(12961,-4246){\makebox(0,0)[lb]{\smash{\SetFigFont{10}{12.0}{rm}$x_5$}}}
\put(12016,-4246){\makebox(0,0)[lb]{\smash{\SetFigFont{10}{12.0}{rm}$x_4$}}}
\put(11026,-4246){\makebox(0,0)[lb]{\smash{\SetFigFont{10}{12.0}{rm}$x_3$}}}
\end{picture}

%% file: pmc_7a.tex
\begin{picture}(0,0)%
\special{psfile=pmc_7a.ps}%
\end{picture}%
\setlength{\unitlength}{1973sp}%
\begingroup\makeatletter\ifx\SetFigFont\undefined
\def\x#1#2#3#4#5#6#7\relax{\def\x{#1#2#3#4#5#6}}%
\expandafter\x\fmtname xxxxxx\relax \def\y{splain}%
\ifx\x\y   
\gdef\SetFigFont#1#2#3{%
  \ifnum #1<17\tiny\else \ifnum #1<20\small\else
  \ifnum #1<24\normalsize\else \ifnum #1<29\large\else
  \ifnum #1<34\Large\else \ifnum #1<41\LARGE\else
     \huge\fi\fi\fi\fi\fi\fi
  \csname #3\endcsname}%
\else
\gdef\SetFigFont#1#2#3{\begingroup
  \count@#1\relax \ifnum 25<\count@\count@25\fi
  \def\x{\endgroup\@setsize\SetFigFont{#2pt}}%
  \expandafter\x
    \csname \romannumeral\the\count@ pt\expandafter\endcsname
    \csname @\romannumeral\the\count@ pt\endcsname
  \csname #3\endcsname}%
\fi
\fi\endgroup
\begin{picture}(10380,6093)(1576,-7756)
\put(11746,-4291){\makebox(0,0)[lb]{\smash{\SetFigFont{10}{14.4}{rm}$b_1$}}}
\put(2071,-2356){\makebox(0,0)[lb]{\smash{\SetFigFont{10}{14.4}{rm}Type $A_{2p}$}}}
\put(4456,-2131){\makebox(0,0)[lb]{\smash{\SetFigFont{10}{14.4}{rm}$\sigma_1$}}}
\put(5401,-1951){\makebox(0,0)[lb]{\smash{\SetFigFont{10}{14.4}{rm}$\sigma_2$}}}
\put(7156,-1951){\makebox(0,0)[lb]{\smash{\SetFigFont{10}{14.4}{rm}$\sigma_4$}}}
\put(1576,-5641){\makebox(0,0)[lb]{\smash{\SetFigFont{10}{14.4}{rm}Type $A_{2p+1}$}}}
\put(4456,-5371){\makebox(0,0)[lb]{\smash{\SetFigFont{10}{14.4}{rm}$\sigma_1$}}}
\put(5356,-5236){\makebox(0,0)[lb]{\smash{\SetFigFont{10}{14.4}{rm}$\sigma_2$}}}
\put(7336,-5371){\makebox(0,0)[lb]{\smash{\SetFigFont{10}{14.4}{rm}$\sigma_4$}}}
\put(6301,-2131){\makebox(0,0)[lb]{\smash{\SetFigFont{10}{14.4}{rm}$\sigma_3$}}}
\put(9856,-2041){\makebox(0,0)[lb]{\smash{\SetFigFont{10}{14.4}{rm}$\sigma_{2p}$}}}
\put(6346,-5371){\makebox(0,0)[lb]{\smash{\SetFigFont{10}{14.4}{rm}$\sigma_3$}}}
\put(10756,-5461){\makebox(0,0)[lb]{\smash{\SetFigFont{10}{14.4}{rm}$\sigma_{2p+1}$}}}
\put(11701,-7756){\makebox(0,0)[lb]{\smash{\SetFigFont{10}{14.4}{rm}$b_2$}}}
\put(11386,-3391){\makebox(0,0)[lb]{\smash{\SetFigFont{10}{14.4}{rm}$b_1$}}}
\end{picture}

%% file: pmc_7b.tex
\begin{picture}(0,0)%
\special{psfile=pmc_7b.ps}%
\end{picture}%
\setlength{\unitlength}{1973sp}%
\begingroup\makeatletter\ifx\SetFigFont\undefined
\def\x#1#2#3#4#5#6#7\relax{\def\x{#1#2#3#4#5#6}}%
\expandafter\x\fmtname xxxxxx\relax \def\y{splain}%
\ifx\x\y   
\gdef\SetFigFont#1#2#3{%
  \ifnum #1<17\tiny\else \ifnum #1<20\small\else
  \ifnum #1<24\normalsize\else \ifnum #1<29\large\else
  \ifnum #1<34\Large\else \ifnum #1<41\LARGE\else
     \huge\fi\fi\fi\fi\fi\fi
  \csname #3\endcsname}%
\else
\gdef\SetFigFont#1#2#3{\begingroup
  \count@#1\relax \ifnum 25<\count@\count@25\fi
  \def\x{\endgroup\@setsize\SetFigFont{#2pt}}%
  \expandafter\x
    \csname \romannumeral\the\count@ pt\expandafter\endcsname
    \csname @\romannumeral\the\count@ pt\endcsname
  \csname #3\endcsname}%
\fi
\fi\endgroup
\begin{picture}(10065,6723)(2206,-5884)
\put(5356,-5101){\makebox(0,0)[lb]{\smash{\SetFigFont{10}{14.4}{rm}$\sigma_3$}}}
\put(2251,-961){\makebox(0,0)[lb]{\smash{\SetFigFont{10}{14.4}{rm}Type $D_{2p}$}}}
\put(2206,-4156){\makebox(0,0)[lb]{\smash{\SetFigFont{10}{14.4}{rm}Type $D_{2p+1}$}}}
\put(5671,-2311){\makebox(0,0)[lb]{\smash{\SetFigFont{10}{14.4}{rm}$\sigma_2$}}}
\put(6526,-826){\makebox(0,0)[lb]{\smash{\SetFigFont{10}{14.4}{rm}$\sigma_4$}}}
\put(7606,-691){\makebox(0,0)[lb]{\smash{\SetFigFont{10}{14.4}{rm}$\sigma_5$}}}
\put(11296,-601){\makebox(0,0)[lb]{\smash{\SetFigFont{10}{14.4}{rm}$\sigma_{2p}$}}}
\put(12016,-2986){\makebox(0,0)[lb]{\smash{\SetFigFont{10}{14.4}{rm}$b_2$}}}
\put(4456,-1951){\makebox(0,0)[lb]{\smash{\SetFigFont{10}{14.4}{rm}$b_3$}}}
\put(5581, 29){\makebox(0,0)[lb]{\smash{\SetFigFont{10}{14.4}{rm}$\sigma_1$}}}
\put(5176,-1636){\makebox(0,0)[lb]{\smash{\SetFigFont{10}{14.4}{rm}$\sigma_3$}}}
\put(10081,-511){\makebox(0,0)[lb]{\smash{\SetFigFont{10}{14.4}{rm}$\sigma_{2p-1}$}}}
\put(11971,479){\makebox(0,0)[lb]{\smash{\SetFigFont{10}{14.4}{rm}$b_1$}}}
\put(5896,-3346){\makebox(0,0)[lb]{\smash{\SetFigFont{10}{14.4}{rm}$\sigma_1$}}}
\put(7696,-4066){\makebox(0,0)[lb]{\smash{\SetFigFont{10}{14.4}{rm}$\sigma_5$}}}
\put(6796,-4111){\makebox(0,0)[lb]{\smash{\SetFigFont{10}{14.4}{rm}$\sigma_4$}}}
\put(10351,-4021){\makebox(0,0)[lb]{\smash{\SetFigFont{10}{14.4}{rm}$\sigma_{2p+1}$}}}
\put(11791,-5326){\makebox(0,0)[lb]{\smash{\SetFigFont{10}{14.4}{rm}$b_1$}}}
\put(4636,-5371){\makebox(0,0)[lb]{\smash{\SetFigFont{10}{14.4}{rm}$b_2$}}}
\put(5986,-5776){\makebox(0,0)[lb]{\smash{\SetFigFont{10}{14.4}{rm}$\sigma_2$}}}
\end{picture}

%% file: pmc_7c.tex
\begin{picture}(0,0)%
\special{psfile=pmc_7c.ps}%
\end{picture}%
\setlength{\unitlength}{1973sp}%
\begingroup\makeatletter\ifx\SetFigFont\undefined
\def\x#1#2#3#4#5#6#7\relax{\def\x{#1#2#3#4#5#6}}%
\expandafter\x\fmtname xxxxxx\relax \def\y{splain}%
\ifx\x\y   
\gdef\SetFigFont#1#2#3{%
  \ifnum #1<17\tiny\else \ifnum #1<20\small\else
  \ifnum #1<24\normalsize\else \ifnum #1<29\large\else
  \ifnum #1<34\Large\else \ifnum #1<41\LARGE\else
     \huge\fi\fi\fi\fi\fi\fi
  \csname #3\endcsname}%
\else
\gdef\SetFigFont#1#2#3{\begingroup
  \count@#1\relax \ifnum 25<\count@\count@25\fi
  \def\x{\endgroup\@setsize\SetFigFont{#2pt}}%
  \expandafter\x
    \csname \romannumeral\the\count@ pt\expandafter\endcsname
    \csname @\romannumeral\the\count@ pt\endcsname
  \csname #3\endcsname}%
\fi
\fi\endgroup
\begin{picture}(8670,4701)(451,-4339)
\put(451,-3211){\makebox(0,0)[lb]{\smash{\SetFigFont{10}{14.4}{rm}Type $E_7$}}}
\put(3826,-16){\makebox(0,0)[lb]{\smash{\SetFigFont{10}{14.4}{rm}$\sigma_1$}}}
\put(4816,-61){\makebox(0,0)[lb]{\smash{\SetFigFont{10}{14.4}{rm}$\sigma_2$}}}
\put(6526,-106){\makebox(0,0)[lb]{\smash{\SetFigFont{10}{14.4}{rm}$\sigma_4$}}}
\put(5581,-2761){\makebox(0,0)[lb]{\smash{\SetFigFont{10}{14.4}{rm}$\sigma_4$}}}
\put(5671, 29){\makebox(0,0)[lb]{\smash{\SetFigFont{10}{14.4}{rm}$\sigma_3$}}}
\put(7741, 74){\makebox(0,0)[lb]{\smash{\SetFigFont{10}{14.4}{rm}$\sigma_5$}}}
\put(6211,-1096){\makebox(0,0)[lb]{\smash{\SetFigFont{10}{14.4}{rm}$\sigma_6$}}}
\put(8821,-1321){\makebox(0,0)[lb]{\smash{\SetFigFont{10}{14.4}{rm}$b_1$}}}
\put(3916,-2176){\makebox(0,0)[lb]{\smash{\SetFigFont{10}{14.4}{rm}$\sigma_1$}}}
\put(3376,-3841){\makebox(0,0)[lb]{\smash{\SetFigFont{10}{14.4}{rm}$\sigma_2$}}}
\put(4726,-2896){\makebox(0,0)[lb]{\smash{\SetFigFont{10}{14.4}{rm}$\sigma_3$}}}
\put(2566,-4156){\makebox(0,0)[lb]{\smash{\SetFigFont{10}{14.4}{rm}$b_2$}}}
\put(6616,-2941){\makebox(0,0)[lb]{\smash{\SetFigFont{10}{14.4}{rm}$\sigma_5$}}}
\put(7606,-2761){\makebox(0,0)[lb]{\smash{\SetFigFont{10}{14.4}{rm}$\sigma_6$}}}
\put(6301,-3931){\makebox(0,0)[lb]{\smash{\SetFigFont{10}{14.4}{rm}$\sigma_7$}}}
\put(8821,-4156){\makebox(0,0)[lb]{\smash{\SetFigFont{10}{14.4}{rm}$b_1$}}}
\put(541,-241){\makebox(0,0)[lb]{\smash{\SetFigFont{10}{14.4}{rm}Type $E_6$}}}
\end{picture}

%% file: pmc_8.tex
\begin{picture}(0,0)%
\special{psfile=pmc_8.ps}%
\end{picture}%
\setlength{\unitlength}{1973sp}%
\begingroup\makeatletter\ifx\SetFigFont\undefined
\def\x#1#2#3#4#5#6#7\relax{\def\x{#1#2#3#4#5#6}}%
\expandafter\x\fmtname xxxxxx\relax \def\y{splain}%
\ifx\x\y   
\gdef\SetFigFont#1#2#3{%
  \ifnum #1<17\tiny\else \ifnum #1<20\small\else
  \ifnum #1<24\normalsize\else \ifnum #1<29\large\else
  \ifnum #1<34\Large\else \ifnum #1<41\LARGE\else
     \huge\fi\fi\fi\fi\fi\fi
  \csname #3\endcsname}%
\else
\gdef\SetFigFont#1#2#3{\begingroup
  \count@#1\relax \ifnum 25<\count@\count@25\fi
  \def\x{\endgroup\@setsize\SetFigFont{#2pt}}%
  \expandafter\x
    \csname \romannumeral\the\count@ pt\expandafter\endcsname
    \csname @\romannumeral\the\count@ pt\endcsname
  \csname #3\endcsname}%
\fi
\fi\endgroup
\begin{picture}(10080,6243)(1081,-6889)
\put(10306,-6676){\makebox(0,0)[lb]{\smash{\SetFigFont{10}{14.4}{rm}$b_1$}}}
\put(1351,-2761){\makebox(0,0)[lb]{\smash{\SetFigFont{10}{14.4}{rm}Type $B_{2p}$}}}
\put(3691,-2941){\makebox(0,0)[lb]{\smash{\SetFigFont{10}{14.4}{rm}$\tau_1$}}}
\put(3736,-2311){\makebox(0,0)[lb]{\smash{\SetFigFont{10}{14.4}{rm}$\sigma_2$}}}
\put(4186,-3121){\makebox(0,0)[lb]{\smash{\SetFigFont{10}{14.4}{rm}$\sigma_3$}}}
\put(5491,-2671){\makebox(0,0)[lb]{\smash{\SetFigFont{10}{14.4}{rm}$\sigma_4$}}}
\put(5806,-3076){\makebox(0,0)[lb]{\smash{\SetFigFont{10}{14.4}{rm}$\sigma_5$}}}
\put(8686,-3121){\makebox(0,0)[lb]{\smash{\SetFigFont{10}{14.4}{rm}$\sigma_{2p-1}$}}}
\put(10126,-2716){\makebox(0,0)[lb]{\smash{\SetFigFont{10}{14.4}{rm}$\sigma_{2p}$}}}
\put(11161,-1006){\makebox(0,0)[lb]{\smash{\SetFigFont{10}{14.4}{rm}$b_1$}}}
\put(10936,-4516){\makebox(0,0)[lb]{\smash{\SetFigFont{10}{14.4}{rm}$b_2$}}}
\put(4456,-1411){\makebox(0,0)[lb]{\smash{\SetFigFont{10}{14.4}{rm}$a_1$}}}
\put(6211,-3931){\makebox(0,0)[lb]{\smash{\SetFigFont{10}{14.4}{rm}$a_2$}}}
\put(1081,-5956){\makebox(0,0)[lb]{\smash{\SetFigFont{10}{14.4}{rm}Type $B_{2p+1}$}}}
\put(3646,-5596){\makebox(0,0)[lb]{\smash{\SetFigFont{10}{14.4}{rm}$\sigma_2$}}}
\put(4366,-5281){\makebox(0,0)[lb]{\smash{\SetFigFont{10}{14.4}{rm}$\sigma_3$}}}
\put(7651,-2401){\makebox(0,0)[lb]{\smash{\SetFigFont{10}{14.4}{rm}$a_3$}}}
\put(10486,-3526){\makebox(0,0)[lb]{\smash{\SetFigFont{10}{14.4}{rm}$\gamma$}}}
\put(3691,-6226){\makebox(0,0)[lb]{\smash{\SetFigFont{10}{14.4}{rm}$\tau_1$}}}
\put(5266,-5416){\makebox(0,0)[lb]{\smash{\SetFigFont{10}{14.4}{rm}$\sigma_4$}}}
\put(6166,-5281){\makebox(0,0)[lb]{\smash{\SetFigFont{10}{14.4}{rm}$\sigma_5$}}}
\put(8911,-5371){\makebox(0,0)[lb]{\smash{\SetFigFont{10}{14.4}{rm}$\sigma_{2p+1}$}}}
\end{picture}

%% file: pmc_9.tex
\begin{picture}(0,0)%
\special{psfile=pmc_9.ps}%
\end{picture}%
\setlength{\unitlength}{1973sp}%
\begingroup\makeatletter\ifx\SetFigFont\undefined
\def\x#1#2#3#4#5#6#7\relax{\def\x{#1#2#3#4#5#6}}%
\expandafter\x\fmtname xxxxxx\relax \def\y{splain}%
\ifx\x\y   
\gdef\SetFigFont#1#2#3{%
  \ifnum #1<17\tiny\else \ifnum #1<20\small\else
  \ifnum #1<24\normalsize\else \ifnum #1<29\large\else
  \ifnum #1<34\Large\else \ifnum #1<41\LARGE\else
     \huge\fi\fi\fi\fi\fi\fi
  \csname #3\endcsname}%
\else
\gdef\SetFigFont#1#2#3{\begingroup
  \count@#1\relax \ifnum 25<\count@\count@25\fi
  \def\x{\endgroup\@setsize\SetFigFont{#2pt}}%
  \expandafter\x
    \csname \romannumeral\the\count@ pt\expandafter\endcsname
    \csname @\romannumeral\the\count@ pt\endcsname
  \csname #3\endcsname}%
\fi
\fi\endgroup
\begin{picture}(8161,1799)(3571,-3166)
\put(10936,-3166){\makebox(0,0)[lb]{\smash{\SetFigFont{10}{14.4}{rm}$b_2$}}}
\put(4411,-3166){\makebox(0,0)[lb]{\smash{\SetFigFont{10}{14.4}{rm}$b_1$}}}
\put(6166,-2626){\makebox(0,0)[lb]{\smash{\SetFigFont{10}{14.4}{rm}$\sigma_1$}}}
\put(5581,-2131){\makebox(0,0)[lb]{\smash{\SetFigFont{10}{14.4}{rm}$\tau_2$}}}
\put(6751,-2131){\makebox(0,0)[lb]{\smash{\SetFigFont{10}{14.4}{rm}$\tau_3$}}}
\put(9496,-2131){\makebox(0,0)[lb]{\smash{\SetFigFont{10}{14.4}{rm}$\tau_l$}}}
\end{picture}

%% file: pmc_10.tex
\begin{picture}(0,0)%
\special{psfile=pmc_10.ps}%
\end{picture}%
\setlength{\unitlength}{1579sp}%
\begingroup\makeatletter\ifx\SetFigFont\undefined
\def\x#1#2#3#4#5#6#7\relax{\def\x{#1#2#3#4#5#6}}%
\expandafter\x\fmtname xxxxxx\relax \def\y{splain}%
\ifx\x\y   
\gdef\SetFigFont#1#2#3{%
  \ifnum #1<17\tiny\else \ifnum #1<20\small\else
  \ifnum #1<24\normalsize\else \ifnum #1<29\large\else
  \ifnum #1<34\Large\else \ifnum #1<41\LARGE\else
     \huge\fi\fi\fi\fi\fi\fi
  \csname #3\endcsname}%
\else
\gdef\SetFigFont#1#2#3{\begingroup
  \count@#1\relax \ifnum 25<\count@\count@25\fi
  \def\x{\endgroup\@setsize\SetFigFont{#2pt}}%
  \expandafter\x
    \csname \romannumeral\the\count@ pt\expandafter\endcsname
    \csname @\romannumeral\the\count@ pt\endcsname
  \csname #3\endcsname}%
\fi
\fi\endgroup
\begin{picture}(14460,6527)(1996,-5598)
\put(3376,-1096){\makebox(0,0)[lb]{\smash{\SetFigFont{10}{12.0}{rm}$a_{r-1}$}}}
\put(5536,-1096){\makebox(0,0)[lb]{\smash{\SetFigFont{10}{12.0}{rm}$a_{r+2}$}}}
\put(8281,-3526){\makebox(0,0)[lb]{\smash{\SetFigFont{10}{12.0}{rm}$b_1$}}}
\put(2251,-4381){\makebox(0,0)[lb]{\smash{\SetFigFont{10}{12.0}{rm}$d_1$}}}
\put(6976,-3751){\makebox(0,0)[lb]{\smash{\SetFigFont{10}{12.0}{rm}$a_0$}}}
\put(3151,-2176){\makebox(0,0)[lb]{\smash{\SetFigFont{10}{12.0}{rm}$a_2$}}}
\put(3961,-331){\makebox(0,0)[lb]{\smash{\SetFigFont{10}{12.0}{rm}$a_r$}}}
\put(4771,-1006){\makebox(0,0)[lb]{\smash{\SetFigFont{10}{12.0}{rm}$a_{r+1}$}}}
\put(7831,-961){\makebox(0,0)[lb]{\smash{\SetFigFont{10}{12.0}{rm}$a_{r+n-1}$}}}
\put(9676,-1411){\makebox(0,0)[lb]{\smash{\SetFigFont{10}{12.0}{rm}$a_{r+n}$}}}
\put(9406,-2806){\makebox(0,0)[lb]{\smash{\SetFigFont{10}{12.0}{rm}$b_2$}}}
\put(10216,-2176){\makebox(0,0)[lb]{\smash{\SetFigFont{10}{12.0}{rm}$b_3$}}}
\put(11476,-2221){\makebox(0,0)[lb]{\smash{\SetFigFont{10}{12.0}{rm}$b_4$}}}
\put(12556,-1996){\makebox(0,0)[lb]{\smash{\SetFigFont{10}{12.0}{rm}$b_5$}}}
\put(14986,-2041){\makebox(0,0)[lb]{\smash{\SetFigFont{10}{12.0}{rm}$b_{2g-1}$}}}
\put(10846,-3886){\makebox(0,0)[lb]{\smash{\SetFigFont{10}{12.0}{rm}$c$}}}
\put(2386,569){\makebox(0,0)[lb]{\smash{\SetFigFont{10}{12.0}{rm}$d_r$}}}
\put(4726,-1951){\makebox(0,0)[lb]{\smash{\SetFigFont{10}{12.0}{rm}$\tau_1$}}}
\put(5671,-1996){\makebox(0,0)[lb]{\smash{\SetFigFont{10}{12.0}{rm}$\tau_2$}}}
\put(2206,-2491){\makebox(0,0)[lb]{\smash{\SetFigFont{10}{12.0}{rm}$d_2$}}}
\put(8101,-1951){\makebox(0,0)[lb]{\smash{\SetFigFont{10}{12.0}{rm}$\tau_{n-1}$}}}
\put(3286,-3526){\makebox(0,0)[lb]{\smash{\SetFigFont{10}{12.0}{rm}$a_1$}}}
\end{picture}

%% file: pmc_11.tex
\begin{picture}(0,0)%
\special{psfile=pmc_11.ps}%
\end{picture}%
\setlength{\unitlength}{1579sp}%
\begingroup\makeatletter\ifx\SetFigFont\undefined
\def\x#1#2#3#4#5#6#7\relax{\def\x{#1#2#3#4#5#6}}%
\expandafter\x\fmtname xxxxxx\relax \def\y{splain}%
\ifx\x\y   
\gdef\SetFigFont#1#2#3{%
  \ifnum #1<17\tiny\else \ifnum #1<20\small\else
  \ifnum #1<24\normalsize\else \ifnum #1<29\large\else
  \ifnum #1<34\Large\else \ifnum #1<41\LARGE\else
     \huge\fi\fi\fi\fi\fi\fi
  \csname #3\endcsname}%
\else
\gdef\SetFigFont#1#2#3{\begingroup
  \count@#1\relax \ifnum 25<\count@\count@25\fi
  \def\x{\endgroup\@setsize\SetFigFont{#2pt}}%
  \expandafter\x
    \csname \romannumeral\the\count@ pt\expandafter\endcsname
    \csname @\romannumeral\the\count@ pt\endcsname
  \csname #3\endcsname}%
\fi
\fi\endgroup
\begin{picture}(13393,3870)(236,-3571)
\put(2026,-1141){\makebox(0,0)[lb]{\smash{\SetFigFont{10}{12.0}{rm}$\tau_1$}}}
\put(3691,-3571){\makebox(0,0)[lb]{\smash{\SetFigFont{10}{12.0}{rm}$a_0$}}}
\put(5491,-2581){\makebox(0,0)[lb]{\smash{\SetFigFont{10}{12.0}{rm}$b_1$}}}
\put(7561,-1366){\makebox(0,0)[lb]{\smash{\SetFigFont{10}{12.0}{rm}$b_3$}}}
\put(2071,-241){\makebox(0,0)[lb]{\smash{\SetFigFont{10}{12.0}{rm}$a_1$}}}
\put(2836,-61){\makebox(0,0)[lb]{\smash{\SetFigFont{10}{12.0}{rm}$a_2$}}}
\put(5176,-196){\makebox(0,0)[lb]{\smash{\SetFigFont{10}{12.0}{rm}$a_{n-1}$}}}
\put(7201,-646){\makebox(0,0)[lb]{\smash{\SetFigFont{10}{12.0}{rm}$a_n$}}}
\put(6796,-1321){\makebox(0,0)[lb]{\smash{\SetFigFont{10}{12.0}{rm}$b_2$}}}
\put(9676,-1276){\makebox(0,0)[lb]{\smash{\SetFigFont{10}{12.0}{rm}$b_5$}}}
\put(12466,-1231){\makebox(0,0)[lb]{\smash{\SetFigFont{10}{12.0}{rm}$b_{2g-1}$}}}
\put(7966,-2986){\makebox(0,0)[lb]{\smash{\SetFigFont{10}{12.0}{rm}$c$}}}
\put(2926,-1186){\makebox(0,0)[lb]{\smash{\SetFigFont{10}{12.0}{rm}$\tau_2$}}}
\put(5311,-1051){\makebox(0,0)[lb]{\smash{\SetFigFont{10}{12.0}{rm}$\tau_{n-1}$}}}
\put(8551,-1411){\makebox(0,0)[lb]{\smash{\SetFigFont{10}{12.0}{rm}$b_4$}}}
\end{picture}

%% file: pmc_12.tex
\begin{picture}(0,0)%
\special{psfile=pmc_12.ps}%
\end{picture}%
\setlength{\unitlength}{1579sp}%
\begingroup\makeatletter\ifx\SetFigFont\undefined
\def\x#1#2#3#4#5#6#7\relax{\def\x{#1#2#3#4#5#6}}%
\expandafter\x\fmtname xxxxxx\relax \def\y{splain}%
\ifx\x\y   
\gdef\SetFigFont#1#2#3{%
  \ifnum #1<17\tiny\else \ifnum #1<20\small\else
  \ifnum #1<24\normalsize\else \ifnum #1<29\large\else
  \ifnum #1<34\Large\else \ifnum #1<41\LARGE\else
     \huge\fi\fi\fi\fi\fi\fi
  \csname #3\endcsname}%
\else
\gdef\SetFigFont#1#2#3{\begingroup
  \count@#1\relax \ifnum 25<\count@\count@25\fi
  \def\x{\endgroup\@setsize\SetFigFont{#2pt}}%
  \expandafter\x
    \csname \romannumeral\the\count@ pt\expandafter\endcsname
    \csname @\romannumeral\the\count@ pt\endcsname
  \csname #3\endcsname}%
\fi
\fi\endgroup
\begin{picture}(12576,10553)(2080,-10974)
\put(10531,-9511){\makebox(0,0)[lb]{\smash{\SetFigFont{7}{8.4}{rm}$p+1$}}}
\put(4771,-736){\makebox(0,0)[lb]{\smash{\SetFigFont{10}{12.0}{rm}$i-1$}}}
\put(7336,-781){\makebox(0,0)[lb]{\smash{\SetFigFont{10}{12.0}{rm}$n$}}}
\put(5941,-736){\makebox(0,0)[lb]{\smash{\SetFigFont{10}{12.0}{rm}$i$}}}
\put(7201,-1321){\makebox(0,0)[lb]{\smash{\SetFigFont{10}{12.0}{rm}$\alpha_i$}}}
\put(8731,-9511){\makebox(0,0)[lb]{\smash{\SetFigFont{10}{12.0}{rm}$p$}}}
\put(11656,-5326){\makebox(0,0)[lb]{\smash{\SetFigFont{10}{12.0}{rm}$\beta_{2p-1}$}}}
\put(10711,-5776){\makebox(0,0)[lb]{\smash{\SetFigFont{10}{12.0}{rm}$p$}}}
\put(10171,-9016){\makebox(0,0)[lb]{\smash{\SetFigFont{10}{12.0}{rm}$\beta_{2p}$}}}
\end{picture}

%% file: pmc_13.tex
\begin{picture}(0,0)%
\special{psfile=pmc_13.ps}%
\end{picture}%
\setlength{\unitlength}{1973sp}%
\begingroup\makeatletter\ifx\SetFigFont\undefined
\def\x#1#2#3#4#5#6#7\relax{\def\x{#1#2#3#4#5#6}}%
\expandafter\x\fmtname xxxxxx\relax \def\y{splain}%
\ifx\x\y   
\gdef\SetFigFont#1#2#3{%
  \ifnum #1<17\tiny\else \ifnum #1<20\small\else
  \ifnum #1<24\normalsize\else \ifnum #1<29\large\else
  \ifnum #1<34\Large\else \ifnum #1<41\LARGE\else
     \huge\fi\fi\fi\fi\fi\fi
  \csname #3\endcsname}%
\else
\gdef\SetFigFont#1#2#3{\begingroup
  \count@#1\relax \ifnum 25<\count@\count@25\fi
  \def\x{\endgroup\@setsize\SetFigFont{#2pt}}%
  \expandafter\x
    \csname \romannumeral\the\count@ pt\expandafter\endcsname
    \csname @\romannumeral\the\count@ pt\endcsname
  \csname #3\endcsname}%
\fi
\fi\endgroup
\begin{picture}(4561,1979)(2221,-2805)
\put(5716,-1186){\makebox(0,0)[lb]{\smash{\SetFigFont{10}{14.4}{rm}$a_2$}}}
\put(3106,-1231){\makebox(0,0)[lb]{\smash{\SetFigFont{10}{14.4}{rm}$a_0$}}}
\put(4276,-1186){\makebox(0,0)[lb]{\smash{\SetFigFont{10}{14.4}{rm}$a_1$}}}
\put(5041,-1726){\makebox(0,0)[lb]{\smash{\SetFigFont{10}{14.4}{rm}$\tau$}}}
\put(3961,-2266){\makebox(0,0)[lb]{\smash{\SetFigFont{10}{14.4}{rm}$P_1$}}}
\put(5356,-2266){\makebox(0,0)[lb]{\smash{\SetFigFont{10}{14.4}{rm}$P_2$}}}
\end{picture}

%% file: pmc_14.tex
\begin{picture}(0,0)%
\special{psfile=pmc_14.ps}%
\end{picture}%
\setlength{\unitlength}{1973sp}%
\begingroup\makeatletter\ifx\SetFigFont\undefined
\def\x#1#2#3#4#5#6#7\relax{\def\x{#1#2#3#4#5#6}}%
\expandafter\x\fmtname xxxxxx\relax \def\y{splain}%
\ifx\x\y   
\gdef\SetFigFont#1#2#3{%
  \ifnum #1<17\tiny\else \ifnum #1<20\small\else
  \ifnum #1<24\normalsize\else \ifnum #1<29\large\else
  \ifnum #1<34\Large\else \ifnum #1<41\LARGE\else
     \huge\fi\fi\fi\fi\fi\fi
  \csname #3\endcsname}%
\else
\gdef\SetFigFont#1#2#3{\begingroup
  \count@#1\relax \ifnum 25<\count@\count@25\fi
  \def\x{\endgroup\@setsize\SetFigFont{#2pt}}%
  \expandafter\x
    \csname \romannumeral\the\count@ pt\expandafter\endcsname
    \csname @\romannumeral\the\count@ pt\endcsname
  \csname #3\endcsname}%
\fi
\fi\endgroup
\begin{picture}(6825,1650)(451,-1486)
\put(4816,-1411){\makebox(0,0)[lb]{\smash{\SetFigFont{10}{14.4}{rm}$z$}}}
\put(1801,-196){\makebox(0,0)[lb]{\smash{\SetFigFont{10}{14.4}{rm}$x_0$}}}
\put(2611,-196){\makebox(0,0)[lb]{\smash{\SetFigFont{10}{14.4}{rm}$y_1$}}}
\put(3511,-196){\makebox(0,0)[lb]{\smash{\SetFigFont{10}{14.4}{rm}$y_2$}}}
\put(4411,-196){\makebox(0,0)[lb]{\smash{\SetFigFont{10}{14.4}{rm}$y_3$}}}
\put(5311,-196){\makebox(0,0)[lb]{\smash{\SetFigFont{10}{14.4}{rm}$y_4$}}}
\put(7111,-196){\makebox(0,0)[lb]{\smash{\SetFigFont{10}{14.4}{rm}$y_{2g-1}$}}}
\put(451,-961){\makebox(0,0)[lb]{\smash{\SetFigFont{10}{14.4}{rm}$\Gamma_g$}}}
\end{picture}

%% file: pmc_15.tex
\begin{picture}(0,0)%
\special{psfile=pmc_15.ps}%
\end{picture}%
\setlength{\unitlength}{1973sp}%
\begingroup\makeatletter\ifx\SetFigFont\undefined
\def\x#1#2#3#4#5#6#7\relax{\def\x{#1#2#3#4#5#6}}%
\expandafter\x\fmtname xxxxxx\relax \def\y{splain}%
\ifx\x\y   
\gdef\SetFigFont#1#2#3{%
  \ifnum #1<17\tiny\else \ifnum #1<20\small\else
  \ifnum #1<24\normalsize\else \ifnum #1<29\large\else
  \ifnum #1<34\Large\else \ifnum #1<41\LARGE\else
     \huge\fi\fi\fi\fi\fi\fi
  \csname #3\endcsname}%
\else
\gdef\SetFigFont#1#2#3{\begingroup
  \count@#1\relax \ifnum 25<\count@\count@25\fi
  \def\x{\endgroup\@setsize\SetFigFont{#2pt}}%
  \expandafter\x
    \csname \romannumeral\the\count@ pt\expandafter\endcsname
    \csname @\romannumeral\the\count@ pt\endcsname
  \csname #3\endcsname}%
\fi
\fi\endgroup
\begin{picture}(9383,2256)(2123,-3739)
\put(8101,-1771){\makebox(0,0)[lb]{\smash{\SetFigFont{10}{14.4}{rm}$\sigma_2$}}}
\put(6301,-1816){\makebox(0,0)[lb]{\smash{\SetFigFont{10}{14.4}{rm}$\sigma_1$}}}
\end{picture}

%% file: pmc_16.tex
\begin{picture}(0,0)%
\special{psfile=pmc_16.ps}%
\end{picture}%
\setlength{\unitlength}{1973sp}%
\begingroup\makeatletter\ifx\SetFigFont\undefined
\def\x#1#2#3#4#5#6#7\relax{\def\x{#1#2#3#4#5#6}}%
\expandafter\x\fmtname xxxxxx\relax \def\y{splain}%
\ifx\x\y   
\gdef\SetFigFont#1#2#3{%
  \ifnum #1<17\tiny\else \ifnum #1<20\small\else
  \ifnum #1<24\normalsize\else \ifnum #1<29\large\else
  \ifnum #1<34\Large\else \ifnum #1<41\LARGE\else
     \huge\fi\fi\fi\fi\fi\fi
  \csname #3\endcsname}%
\else
\gdef\SetFigFont#1#2#3{\begingroup
  \count@#1\relax \ifnum 25<\count@\count@25\fi
  \def\x{\endgroup\@setsize\SetFigFont{#2pt}}%
  \expandafter\x
    \csname \romannumeral\the\count@ pt\expandafter\endcsname
    \csname @\romannumeral\the\count@ pt\endcsname
  \csname #3\endcsname}%
\fi
\fi\endgroup
\begin{picture}(10215,2790)(3151,-3211)
\put(9676,-1321){\makebox(0,0)[lb]{\smash{\SetFigFont{10}{14.4}{rm}$v_1$}}}
\put(8551,-3211){\makebox(0,0)[lb]{\smash{\SetFigFont{10}{14.4}{rm}$x_0$}}}
\put(7921,-2941){\makebox(0,0)[lb]{\smash{\SetFigFont{10}{14.4}{rm}$x_1$}}}
\put(3151,-2356){\makebox(0,0)[lb]{\smash{\SetFigFont{10}{14.4}{rm}$\Gamma_{g,r,n}$}}}
\put(4411,-2221){\makebox(0,0)[lb]{\smash{\SetFigFont{10}{14.4}{rm}$u_1$}}}
\put(5311,-2221){\makebox(0,0)[lb]{\smash{\SetFigFont{10}{14.4}{rm}$u_2$}}}
\put(6751,-2221){\makebox(0,0)[lb]{\smash{\SetFigFont{10}{14.4}{rm}$u_r$}}}
\put(8911,-2221){\makebox(0,0)[lb]{\smash{\SetFigFont{10}{14.4}{rm}$y_1$}}}
\put(9676,-2221){\makebox(0,0)[lb]{\smash{\SetFigFont{10}{14.4}{rm}$y_2$}}}
\put(10756,-2221){\makebox(0,0)[lb]{\smash{\SetFigFont{10}{14.4}{rm}$y_3$}}}
\put(13366,-2221){\makebox(0,0)[lb]{\smash{\SetFigFont{10}{14.4}{rm}$y_{2g-1}$}}}
\put(10801,-2851){\makebox(0,0)[lb]{\smash{\SetFigFont{10}{14.4}{rm}$z$}}}
\put(10576,-1321){\makebox(0,0)[lb]{\smash{\SetFigFont{10}{14.4}{rm}$v_2$}}}
\put(12421,-1276){\makebox(0,0)[lb]{\smash{\SetFigFont{10}{14.4}{rm}$v_{n-1}$}}}
\put(7876,-916){\makebox(0,0)[lb]{\smash{\SetFigFont{10}{14.4}{rm}$x_r$}}}
\put(11476,-2221){\makebox(0,0)[lb]{\smash{\SetFigFont{10}{14.4}{rm}$y_4$}}}
\put(9091,-736){\makebox(0,0)[lb]{\smash{\SetFigFont{10}{14.4}{rm}$4$}}}
\put(8911,-1276){\makebox(0,0)[lb]{\smash{\SetFigFont{10}{14.4}{rm}$x_{r+1}$}}}
\end{picture}

%% file: pmc_17.tex
\begin{picture}(0,0)%
\special{psfile=pmc_17.ps}%
\end{picture}%
\setlength{\unitlength}{1973sp}%
\begingroup\makeatletter\ifx\SetFigFont\undefined
\def\x#1#2#3#4#5#6#7\relax{\def\x{#1#2#3#4#5#6}}%
\expandafter\x\fmtname xxxxxx\relax \def\y{splain}%
\ifx\x\y   
\gdef\SetFigFont#1#2#3{%
  \ifnum #1<17\tiny\else \ifnum #1<20\small\else
  \ifnum #1<24\normalsize\else \ifnum #1<29\large\else
  \ifnum #1<34\Large\else \ifnum #1<41\LARGE\else
     \huge\fi\fi\fi\fi\fi\fi
  \csname #3\endcsname}%
\else
\gdef\SetFigFont#1#2#3{\begingroup
  \count@#1\relax \ifnum 25<\count@\count@25\fi
  \def\x{\endgroup\@setsize\SetFigFont{#2pt}}%
  \expandafter\x
    \csname \romannumeral\the\count@ pt\expandafter\endcsname
    \csname @\romannumeral\the\count@ pt\endcsname
  \csname #3\endcsname}%
\fi
\fi\endgroup
\begin{picture}(7965,2655)(5401,-3211)
\put(9001,-916){\makebox(0,0)[lb]{\smash{\SetFigFont{10}{14.4}{rm}$x_n$}}}
\put(8551,-3211){\makebox(0,0)[lb]{\smash{\SetFigFont{10}{14.4}{rm}$x_0$}}}
\put(7921,-2941){\makebox(0,0)[lb]{\smash{\SetFigFont{10}{14.4}{rm}$x_1$}}}
\put(8911,-2221){\makebox(0,0)[lb]{\smash{\SetFigFont{10}{14.4}{rm}$y_1$}}}
\put(9676,-2221){\makebox(0,0)[lb]{\smash{\SetFigFont{10}{14.4}{rm}$y_2$}}}
\put(10756,-2221){\makebox(0,0)[lb]{\smash{\SetFigFont{10}{14.4}{rm}$y_3$}}}
\put(13366,-2221){\makebox(0,0)[lb]{\smash{\SetFigFont{10}{14.4}{rm}$y_{2g-1}$}}}
\put(10801,-2851){\makebox(0,0)[lb]{\smash{\SetFigFont{10}{14.4}{rm}$z$}}}
\put(11476,-2221){\makebox(0,0)[lb]{\smash{\SetFigFont{10}{14.4}{rm}$y_4$}}}
\put(7876,-916){\makebox(0,0)[lb]{\smash{\SetFigFont{10}{14.4}{rm}$x_{n-1}$}}}
\put(5401,-2311){\makebox(0,0)[lb]{\smash{\SetFigFont{10}{14.4}{rm}$P\Gamma_{g,0,n}$}}}
\end{picture}

%% file: pmc_18.tex
\begin{picture}(0,0)%
\special{psfile=pmc_18.ps}%
\end{picture}%
\setlength{\unitlength}{1973sp}%
\begingroup\makeatletter\ifx\SetFigFont\undefined
\def\x#1#2#3#4#5#6#7\relax{\def\x{#1#2#3#4#5#6}}%
\expandafter\x\fmtname xxxxxx\relax \def\y{splain}%
\ifx\x\y   
\gdef\SetFigFont#1#2#3{%
  \ifnum #1<17\tiny\else \ifnum #1<20\small\else
  \ifnum #1<24\normalsize\else \ifnum #1<29\large\else
  \ifnum #1<34\Large\else \ifnum #1<41\LARGE\else
     \huge\fi\fi\fi\fi\fi\fi
  \csname #3\endcsname}%
\else
\gdef\SetFigFont#1#2#3{\begingroup
  \count@#1\relax \ifnum 25<\count@\count@25\fi
  \def\x{\endgroup\@setsize\SetFigFont{#2pt}}%
  \expandafter\x
    \csname \romannumeral\the\count@ pt\expandafter\endcsname
    \csname @\romannumeral\the\count@ pt\endcsname
  \csname #3\endcsname}%
\fi
\fi\endgroup
\begin{picture}(1950,2115)(7126,-3346)
\put(9001,-2716){\makebox(0,0)[lb]{\smash{\SetFigFont{10}{14.4}{rm}$x_4$}}}
\put(8326,-3346){\makebox(0,0)[lb]{\smash{\SetFigFont{10}{14.4}{rm}$x_1$}}}
\put(8281,-1591){\makebox(0,0)[lb]{\smash{\SetFigFont{10}{14.4}{rm}$x_3$}}}
\put(7201,-2716){\makebox(0,0)[lb]{\smash{\SetFigFont{10}{14.4}{rm}$x_2$}}}
\put(8281,-2716){\makebox(0,0)[lb]{\smash{\SetFigFont{10}{14.4}{rm}$y$}}}
\end{picture}

%% file: pmc_19.tex
\begin{picture}(0,0)%
\special{psfile=pmc_19.ps}%
\end{picture}%
\setlength{\unitlength}{1973sp}%
\begingroup\makeatletter\ifx\SetFigFont\undefined
\def\x#1#2#3#4#5#6#7\relax{\def\x{#1#2#3#4#5#6}}%
\expandafter\x\fmtname xxxxxx\relax \def\y{splain}%
\ifx\x\y   
\gdef\SetFigFont#1#2#3{%
  \ifnum #1<17\tiny\else \ifnum #1<20\small\else
  \ifnum #1<24\normalsize\else \ifnum #1<29\large\else
  \ifnum #1<34\Large\else \ifnum #1<41\LARGE\else
     \huge\fi\fi\fi\fi\fi\fi
  \csname #3\endcsname}%
\else
\gdef\SetFigFont#1#2#3{\begingroup
  \count@#1\relax \ifnum 25<\count@\count@25\fi
  \def\x{\endgroup\@setsize\SetFigFont{#2pt}}%
  \expandafter\x
    \csname \romannumeral\the\count@ pt\expandafter\endcsname
    \csname @\romannumeral\the\count@ pt\endcsname
  \csname #3\endcsname}%
\fi
\fi\endgroup
\begin{picture}(10485,2655)(2881,-3211)
\put(10891,-2851){\makebox(0,0)[lb]{\smash{\SetFigFont{10}{14.4}{rm}$z$}}}
\put(8551,-3211){\makebox(0,0)[lb]{\smash{\SetFigFont{10}{14.4}{rm}$x_0$}}}
\put(7921,-2941){\makebox(0,0)[lb]{\smash{\SetFigFont{10}{14.4}{rm}$x_1$}}}
\put(4411,-2221){\makebox(0,0)[lb]{\smash{\SetFigFont{10}{14.4}{rm}$u_1$}}}
\put(5311,-2221){\makebox(0,0)[lb]{\smash{\SetFigFont{10}{14.4}{rm}$u_2$}}}
\put(6751,-2221){\makebox(0,0)[lb]{\smash{\SetFigFont{10}{14.4}{rm}$u_r$}}}
\put(8911,-2221){\makebox(0,0)[lb]{\smash{\SetFigFont{10}{14.4}{rm}$y_1$}}}
\put(9676,-2221){\makebox(0,0)[lb]{\smash{\SetFigFont{10}{14.4}{rm}$y_2$}}}
\put(13366,-2221){\makebox(0,0)[lb]{\smash{\SetFigFont{10}{14.4}{rm}$y_{2g-1}$}}}
\put(11476,-2221){\makebox(0,0)[lb]{\smash{\SetFigFont{10}{14.4}{rm}$y_4$}}}
\put(2881,-2311){\makebox(0,0)[lb]{\smash{\SetFigFont{10}{14.4}{rm}$P\Gamma_{g,r,n}$}}}
\put(8956,-916){\makebox(0,0)[lb]{\smash{\SetFigFont{10}{14.4}{rm}$x_{n+r}$}}}
\put(10846,-2221){\makebox(0,0)[lb]{\smash{\SetFigFont{10}{14.4}{rm}$y_3$}}}
\end{picture}

%% file: pmc_20.tex
\begin{picture}(0,0)%
\special{psfile=pmc_20.ps}%
\end{picture}%
\setlength{\unitlength}{1579sp}%
\begingroup\makeatletter\ifx\SetFigFont\undefined
\def\x#1#2#3#4#5#6#7\relax{\def\x{#1#2#3#4#5#6}}%
\expandafter\x\fmtname xxxxxx\relax \def\y{splain}%
\ifx\x\y   
\gdef\SetFigFont#1#2#3{%
  \ifnum #1<17\tiny\else \ifnum #1<20\small\else
  \ifnum #1<24\normalsize\else \ifnum #1<29\large\else
  \ifnum #1<34\Large\else \ifnum #1<41\LARGE\else
     \huge\fi\fi\fi\fi\fi\fi
  \csname #3\endcsname}%
\else
\gdef\SetFigFont#1#2#3{\begingroup
  \count@#1\relax \ifnum 25<\count@\count@25\fi
  \def\x{\endgroup\@setsize\SetFigFont{#2pt}}%
  \expandafter\x
    \csname \romannumeral\the\count@ pt\expandafter\endcsname
    \csname @\romannumeral\the\count@ pt\endcsname
  \csname #3\endcsname}%
\fi
\fi\endgroup
\begin{picture}(12801,2938)(3880,-3780)
\put(15571,-2986){\makebox(0,0)[lb]{\smash{\SetFigFont{10}{12.0}{rm}$e'$}}}
\put(7561,-1276){\makebox(0,0)[lb]{\smash{\SetFigFont{10}{12.0}{rm}$a_n$}}}
\put(10576,-3121){\makebox(0,0)[lb]{\smash{\SetFigFont{10}{12.0}{rm}$a_n'$}}}
\put(15841,-1411){\makebox(0,0)[lb]{\smash{\SetFigFont{10}{12.0}{rm}$e$}}}
\end{picture}

%% file: pmc_21.tex
\begin{picture}(0,0)%
\special{psfile=pmc_21.ps}%
\end{picture}%
\setlength{\unitlength}{1579sp}%
\begingroup\makeatletter\ifx\SetFigFont\undefined
\def\x#1#2#3#4#5#6#7\relax{\def\x{#1#2#3#4#5#6}}%
\expandafter\x\fmtname xxxxxx\relax \def\y{splain}%
\ifx\x\y   
\gdef\SetFigFont#1#2#3{%
  \ifnum #1<17\tiny\else \ifnum #1<20\small\else
  \ifnum #1<24\normalsize\else \ifnum #1<29\large\else
  \ifnum #1<34\Large\else \ifnum #1<41\LARGE\else
     \huge\fi\fi\fi\fi\fi\fi
  \csname #3\endcsname}%
\else
\gdef\SetFigFont#1#2#3{\begingroup
  \count@#1\relax \ifnum 25<\count@\count@25\fi
  \def\x{\endgroup\@setsize\SetFigFont{#2pt}}%
  \expandafter\x
    \csname \romannumeral\the\count@ pt\expandafter\endcsname
    \csname @\romannumeral\the\count@ pt\endcsname
  \csname #3\endcsname}%
\fi
\fi\endgroup
\begin{picture}(12449,10553)(2080,-10974)
\put(7651,-826){\makebox(0,0)[lb]{\smash{\SetFigFont{10}{12.0}{rm}$Q$}}}
\put(8731,-9511){\makebox(0,0)[lb]{\smash{\SetFigFont{10}{12.0}{rm}$p$}}}
\put(10711,-5776){\makebox(0,0)[lb]{\smash{\SetFigFont{10}{12.0}{rm}$p$}}}
\put(10531,-9511){\makebox(0,0)[lb]{\smash{\SetFigFont{7}{8.4}{rm}$p+1$}}}
\put(5941,-736){\makebox(0,0)[lb]{\smash{\SetFigFont{10}{12.0}{rm}$i+1$}}}
\put(10171,-9016){\makebox(0,0)[lb]{\smash{\SetFigFont{10}{12.0}{rm}$\bar\beta_{2p}$}}}
\put(4996,-736){\makebox(0,0)[lb]{\smash{\SetFigFont{10}{12.0}{rm}$i$}}}
\put(7291,-1321){\makebox(0,0)[lb]{\smash{\SetFigFont{10}{12.0}{rm}$\bar\alpha_i$}}}
\put(5806,-5191){\makebox(0,0)[lb]{\smash{\SetFigFont{10}{12.0}{rm}$Q$}}}
\put(11656,-5416){\makebox(0,0)[lb]{\smash{\SetFigFont{10}{12.0}{rm}$\bar\beta_{2p-1}$}}}
\put(5581,-8971){\makebox(0,0)[lb]{\smash{\SetFigFont{10}{12.0}{rm}$Q$}}}
\end{picture}

%% file: pmc_22.tex
\begin{picture}(0,0)%
\special{psfile=pmc_22.ps}%
\end{picture}%
\setlength{\unitlength}{1973sp}%
\begingroup\makeatletter\ifx\SetFigFont\undefined
\def\x#1#2#3#4#5#6#7\relax{\def\x{#1#2#3#4#5#6}}%
\expandafter\x\fmtname xxxxxx\relax \def\y{splain}%
\ifx\x\y   
\gdef\SetFigFont#1#2#3{%
  \ifnum #1<17\tiny\else \ifnum #1<20\small\else
  \ifnum #1<24\normalsize\else \ifnum #1<29\large\else
  \ifnum #1<34\Large\else \ifnum #1<41\LARGE\else
     \huge\fi\fi\fi\fi\fi\fi
  \csname #3\endcsname}%
\else
\gdef\SetFigFont#1#2#3{\begingroup
  \count@#1\relax \ifnum 25<\count@\count@25\fi
  \def\x{\endgroup\@setsize\SetFigFont{#2pt}}%
  \expandafter\x
    \csname \romannumeral\the\count@ pt\expandafter\endcsname
    \csname @\romannumeral\the\count@ pt\endcsname
  \csname #3\endcsname}%
\fi
\fi\endgroup
\begin{picture}(11451,3825)(2530,-4516)
\put(13006,-1411){\makebox(0,0)[lb]{\smash{\SetFigFont{10}{14.4}{rm}$e_2'$}}}
\put(9856,-2221){\makebox(0,0)[lb]{\smash{\SetFigFont{10}{14.4}{rm}$e_3$}}}
\put(8776,-4471){\makebox(0,0)[lb]{\smash{\SetFigFont{10}{14.4}{rm}$e_2$}}}
\put(11611,-4516){\makebox(0,0)[lb]{\smash{\SetFigFont{10}{14.4}{rm}$e_1'$}}}
\put(13321,-4021){\makebox(0,0)[lb]{\smash{\SetFigFont{10}{14.4}{rm}$e_3'$}}}
\put(11656,-1051){\makebox(0,0)[lb]{\smash{\SetFigFont{10}{14.4}{rm}$e_1$}}}
\end{picture}